\pdfoutput=1
\documentclass[11pt]{article}

\usepackage[letterpaper,margin=1in]{geometry}
\usepackage{amsmath,amssymb,amsthm,mathtools}
\usepackage{mathrsfs}
\usepackage{enumitem}
\usepackage{xr}
\usepackage{hyperref}
\hypersetup{hidelinks,pdftitle={Renting, Owning, and Moving: Housing Tenure Choice with Transaction Costs},pdfsubject={Finite-horizon complete-market housing tenure},pdfkeywords={housing tenure, transaction costs, optimal switching, timing windows, free boundaries, duality}}
\usepackage[round]{natbib}
\usepackage{setspace}
\usepackage[T1]{fontenc}
\usepackage{txfonts}
\usepackage{microtype}
\usepackage{xcolor,graphicx,booktabs,array}
\usepackage[font={small,color=black},labelfont=bf]{caption}
\usepackage{float}
\usepackage{needspace}
\numberwithin{equation}{section}
\graphicspath{{./}}
\newtheorem{assumption}{Assumption}

\newtheorem{theorem}{Theorem}[section]
\newtheorem{lemma}[theorem]{Lemma}
\newtheorem{proposition}[theorem]{Proposition}
\newtheorem{corollary}[theorem]{Corollary}

\newtheorem{remark}{Remark}[section]

\newcommand{\R}{\mathbb{R}}
\newcommand{\E}{\mathbb{E}}
\newcommand{\F}{\mathbb{F}}
\newcommand{\dd}{\,\mathrm{d}}
\newcommand{\cI}{\mathcal{I}}
\newcommand{\cJ}{\mathcal{J}}
\newcommand{\cL}{\mathcal{L}}
\newcommand{\cO}{\mathcal{O}}

\newcommand{\wt}{\widetilde}
\newcommand{\wh}{\widehat}
\newcommand{\one}{\mathbf{1}}
\newcommand{\tr}{\mathsf T}

\newcommand{\vns}{{v_{\rm ns}}}

\usepackage{placeins}

\begin{document}

\begin{center}
{\LARGE\bfseries Renting, Owning, and Moving: Housing Tenure Choice\\[3pt] with Transaction Costs\par}
\vspace{14pt}
{\large Kexin Chen$^{1}$\quad Gugyum Ha$^{2}$\quad Junkee Jeon$^{3}$\par}
\vspace{6pt}
{\small $^{1}$Department of Applied Mathematics, The Hong Kong Polytechnic University, Hong Kong\par}
{\small $^{2}$Department of Mathematics, Sogang University, Seoul 04107, Republic of Korea\par}
{\small $^{3}$Department of Applied Mathematics, Kyung Hee University, Yongin 17104, Republic of Korea\par}
\vspace{6pt}
{\small October 2026\par}
\end{center}
\vspace{4pt}

\begin{abstract}
\noindent
The wealth thresholds at which a household switches between renting and
owning depend on its preferences, but the dates at which some owner sells or
every renter buys do not. In our finite-horizon model, a household chooses
consumption, a portfolio, and costly tenure changes while rents, ownership
expenses, and transaction costs vary over time on announced dates. When financial markets span house-price
risk, these dates are those of a deterministic minimum cost schedule for
housing, discounted at the risk-adjusted housing carry and computed from a
scalar reflected equation. Preferences, housing services, and labor income
only locate the thresholds inside these timing windows, and house-price
coefficients enter the dates only through the carry, because financial
hedging offsets the exposure of the home. When owning costs at least as much
as renting, the windows concern selling, and an anticipated selling
surcharge closes the selling window before the surcharge peaks and keeps it
closed for years afterwards. We also establish strong duality and the
regularity of the thresholds and compute them by integral equations.
Numerically, the windows disappear when short sales or borrowing are
restricted but do not move with the preferences we vary under an
uninsurable increase in the ownership cost.
\end{abstract}

\vspace{6pt}


\section{Introduction}\label{sec:introduction}

Changing housing tenure is among the largest financial decisions a household
makes. A move changes the housing service bundle and consumes resources, and
the cost of moving changes over time. Transaction taxes are raised, cut, and
suspended on announced dates \citep{BestKleven2018,HilberLyytikainen2017},
and rents and ownership expenses move relative to each other. Housing
commitments also shape the rest of the household portfolio
\citep{FlavinYamashita2002,ChettySzeidl2007}. The decision is therefore a
stochastic control problem. A household chooses consumption, a financial
portfolio, and the dates at which it moves between a rented and an owned
home, and in the optimal policy these moves occur at wealth thresholds that
depend on preferences, housing services, and income.

We study this decision in a finite-horizon model of consumption, investment,
and repeated switches between a rented and an owned housing bundle. Rents,
owner expenses, and transaction costs follow general deterministic paths, and
ownership may carry a terminal adjustment charge. We model anticipated cost
changes with smooth schedules, allowing steep transitions as approximations
to announced step changes. Financial markets span
house-price risk and allow unrestricted positions, so the price risk of a
home can be traded but changing the home still requires a costly move. With
time-dependent costs the optimal policy need not consist of a buying and a
selling threshold at every date. At some dates no owner sells at any
admissible wealth, at others every renter buys, and anticipated cost changes
can produce several separate trading windows or isolated contact dates.
Since a sufficiently wealthy renter buys at every date, the first question is
at which dates selling is optimal for some owner and at which dates buying is
optimal for every renter.

Our main result is that these dates are determined by a deterministic
problem. Let $m_0(t)$ and $m_1(t)$ be the minimum cost, per unit of the
house-price index, of financing housing from date $t$ to the horizon for a
renter and for an owner, over deterministic switching schedules and with
costs discounted at the risk-adjusted housing carry.
Theorem~\ref{thm:timing-windows} shows that some owner sells at date $t$
exactly when a minimum cost schedule that starts from owning sells
immediately throughout a neighborhood of $t$ and the local sale incentive is
strictly positive. Every renter buys at $t$ exactly when a minimum cost
schedule that starts from renting buys immediately. At every other date a renter buys only above a
finite wealth threshold and no owner sells. The functions $m_i$ are computed
from a scalar reflected equation and depend only on the interest rate, the
market price of risk, the house-price coefficients, and the cost schedules.
Risk aversion, the consumption share, time preference, the housing service
levels, and labor income determine who moves inside a window, but not when a
window is open.

Ownership comes with a larger service bundle, and this service
advantage is what makes tenure choice depend on wealth. With equal services,
a case our assumptions exclude, a deterministic minimum cost schedule would
be optimal at every admissible wealth level. Financial hedging offsets the
house-price exposure of the home, and house-price coefficients affect the
timing windows only through the risk-adjusted carry. The dates at which
physical flexibility is used are thus set by a deterministic comparison of
costs.

Part of this statement is expected. A household whose wealth is close to its
financing floor has almost no resources beyond those needed to pay for
housing, so its tenure changes must approach the minimum cost schedule. The
content of the theorem is that this limiting behavior decides the question
for every wealth level. No owner, however wealthy or risk averse, sells at a
date at which the schedule does not sell, and when the schedule buys
immediately, every renter buys. The proof identifies the limit of a double
obstacle problem for large dual states with $m_1-m_0$
(Proposition~\ref{prop:minimum-financing}).

The timing windows have a direct reading. When owning costs at least as much as
renting at every later date, as in our baseline, no date is one at which
every renter buys, and the predictions that do not depend on preferences
concern selling. The schedule sells at a date only if no later date before
the next all-buy date offers a larger discounted saving from selling, and
the symmetric statement holds for purchases
(Lemma~\ref{lem:running-record}). In our examples this condition
implies that an anticipated surcharge on selling closes the selling window
about 0.8 years before the surcharge peaks and keeps it closed for 2.4 years
afterwards. It also explains why the terminal treatment of an owned home can
remove sales in the final years, with an explicit formula for the length of
that period. Because these dates follow from a scalar equation, the effect
of an announced cost schedule on whether anyone moves can be computed before
any household parameter is specified.

Within the windows, the thresholds are free boundaries of a parabolic double
obstacle problem with time-dependent obstacles. We establish their local
$C^{1,\alpha}$ regularity and their behavior at internal and terminal
endpoints, without assuming monotonicity in time or finitely many
components. Strong duality holds on the natural financing domain
$x>p\,m_i(t)$, so the functions $m_i$ that determine the windows are also the
financing floors of the primal problem. Integral equations convert the
boundaries into wealth thresholds, consumption, financial hedges, and the
value of tenure flexibility. In the baseline example, halving transaction
costs moves the buying and selling thresholds with ten years remaining from
1.3940 and 0.7852 to 1.2495 and 0.8142 property values, and for a renter
with wealth of 2.5 property values the option to change tenure is worth 11.2
percent of wealth.

Theorem~\ref{thm:timing-windows} rests on two features of the market.
Positions are unrestricted, so the house-price exposure of the home can be
offset at every wealth level, and every cost schedule is known in advance.
Two numerical experiments in the e-companion remove each feature in turn.
When short sales or borrowing are restricted, an owner with little net worth
sells to stay solvent rather than to save on housing costs, and the windows
disappear. When the ownership cost is instead subject to an uninsurable
permanent increase, the dates at which some owner sells remain deterministic
in our computations and are read from the worst-case schedule.

\subsection{Related literature}
The consumption and investment framework is due to
\citet{Merton1969,Merton1971}. Illiquid durables create inaction and affect
portfolios \citep{GrossmanLaroque1990,DamgaardFuglsbjergMunk2003}, and
\citet{BankRiedel2003} obtain a backward integral equation for irreversible
durable accumulation. Housing models incorporate labor income, tenure choice,
borrowing constraints, and predictable price growth
\citep{Cocco2005,YaoZhang2005,CorradinFillatVergara2014}.
\citet{KraftMunk2011} provide the closest complete-market benchmark, with a
continuous housing choice, and we add costly discrete tenure changes between
fixed service bundles. Housing habits, the hedging of housing expenses, and
the portfolio effects of housing are studied by
\citet{FlavinNakagawa2008,Yang2009,KraftMunkWagner2018} and
\citet{SinaiSouleles2005,Han2010,ChettySandorSzeidl2017}. Transaction costs
and the expected length of stay influence tenure \citep{HaurinGill2002},
ownership often persists late in life \citep{VentiWise2004}, and taxation
affects tenure and mobility
\citep{RosenRosen1980,BestKleven2018,HilberLyytikainen2017}.
\citet{HanNgaiSheedy2026} study transaction taxes in housing market
equilibrium, and the holding-period tax in \citet{AgarwalEtAl2026} requires
the duration since purchase as an additional state. We take prices and
calendar-time costs as given and characterize exactly how a known cost path
determines the dates at which households can optimally move.

Costly reversibility is central to entry and exit \citep{Dixit1989} and to
factor adjustment \citep{EberlyVanMieghem1997}, and the operations and
finance literature studies how physical flexibility interacts with financial
hedging \citep{CaldenteyHaugh2006,CaldenteyHaugh2009,ChodRudiVanMieghem2010}.
In our model financial spanning removes house-price exposure, but tenure
flexibility keeps its value (Section~\ref{sec:numerical}), and the dates at
which it can be exercised are decided by a deterministic cost comparison.

The reduction uses martingale duality
\citep{KaratzasLehoczkyShreve1987,CoxHuang1989,KaratzasWang2000}.
Integral representations and free boundary methods for finite horizons are
developed by \citet{Kim1990,Jacka1991,Peskir2005LocalTime,Peskir2005American,
PeskirShiryaev2006}. Related two-boundary problems include investment with
transaction costs \citep{DaiYi2009,DaiJiangLiYi2009}, trend following
\citep{DaiYangZhangZhu2016}, and partially reversible investment
\citep{DeAngelisFerrari2014,HanYiZhang2024}, and continuity and local
regularity of boundaries without time monotonicity are studied by
\citet{DeAngelis2015,DeAngelisStabile2019}. \citet{JeonYang} combine
finite-horizon consumption and investment with costly switching, duality,
and integral representations. \citet{HaJeonOk2026} develop local derivative
and level-curve estimates for time-dependent switching costs under global
incentive inequalities that prescribe the switching interval, with a
terminal datum common to both regimes. Here those inequalities may fail and
the terminal data differ by tenure, so the active time domains are part of
the solution. With constant costs and incomes, \citet{JeonYang} obtain a
single critical remaining time before which one switching region is empty,
their equation (2.4), and \eqref{eq:sale-onset-economic} below is its
analogue for constant housing costs. With time-dependent costs the active
time sets can open and close repeatedly, and their identification through a
reflected scalar equation and a deterministic minimum cost schedule is the
new element of our analysis.

Section~\ref{sec:model} presents the model and its dual reduction,
Section~\ref{sec:dob} contains the main analysis,
Section~\ref{sec:policies} recovers the optimal policies, and
Section~\ref{sec:numerical} presents numerical results. The e-companion
contains auxiliary results, all proofs, the integral representation used for
computation, and additional examples.

\section{Model and dual reduction}\label{sec:model}

On $[0,T]$, a household chooses consumption, financial investment, and costly
switches between rented and owned housing. The horizon may represent retirement
or planned relocation. Financial spanning and unrestricted portfolios provide
a benchmark in which physical adjustment remains costly while house-price
risk can be hedged.
\subsection{Financial market and house-price index}

Let $(\Omega,\mathcal F,\F=(\mathcal F_t)_{t\ge0},\mathbb P)$ satisfy the usual
conditions, and let $\F$ be the usual augmentation of the filtration
generated by a $d$-dimensional Brownian motion $B=(B^1,\ldots,B^d)$. The
risk-free asset satisfies
\begin{equation}\label{eq:S0}
    \dd S^0_t=rS^0_t\,\dd t,
    \qquad r>0.
\end{equation}
There are $d$ traded risky securities $S^1,\ldots,S^d$ satisfying
\begin{equation}\label{eq:risky}
    \frac{\dd S^m_t}{S^m_t}
    =\mu_m(t)\,\dd t+\sigma_m(t)^{\tr}\dd B_t,
    \qquad m=1,\ldots,d.
\end{equation}
Let
\begin{equation}\label{eq:Sigma}
    \Sigma(t):=
    \begin{pmatrix}
        \sigma_1(t)^{\tr}\\
        \vdots\\
        \sigma_d(t)^{\tr}
    \end{pmatrix},
    \qquad
    \mu(t):=(\mu_1(t),\ldots,\mu_d(t))^{\tr}.
\end{equation}
The local house-price index $P$ follows
\begin{equation}\label{eq:P}
    \frac{\dd P_t}{P_t}=\mu_P(t)\,\dd t+\sigma_P(t)^{\tr}\dd B_t,
    \qquad P_0=p>0.
\end{equation}
All coefficient functions are deterministic and bounded.

The index $P$ prices one standardized housing unit. It affects monetary
outlays and the owner's balance sheet. Physical housing services are specified
separately below.

\begin{assumption}\label{ass:market}
The matrix $\Sigma(t)$ is invertible for every $t\in[0,T]$, and there exists
$\lambda_\Sigma>0$ such that
\[
    \Sigma(t)\Sigma(t)^{\tr}\ge \lambda_\Sigma I_d,
    \qquad t\in[0,T].
\]
Moreover, the house-price index $P$ is driven by the same Brownian motion as the
traded securities. One of the risky securities may be interpreted as a REIT, a
housing-index security, or another house-price-linked traded asset.
\end{assumption}

Spanning prices house-price diffusion through traded assets, but changing the
physical home still requires a costly tenure switch.

\begin{assumption}\label{ass:portfolio}
The household may take any progressively measurable dollar position
\(\pi_t\in\R^d\) satisfying the required square-integrability condition.  In
particular, there are no short-sale, leverage, portfolio-cone, mortgage, or
down-payment constraints in the benchmark model, and the traded portfolio may be
rebalanced frictionlessly at a tenure-switching time after the switching cost is
paid.
\end{assumption}

Assumption \ref{ass:portfolio} is essential for the neutralization result below.
Spanning of Brownian risk alone is not sufficient if the household is prevented from
holding the offsetting traded position, a case computed in
Section~\ref{ec:pc} of the e-companion.

The unique market price of Brownian risk is
\begin{equation}\label{eq:theta}
    \theta(t):=\Sigma(t)^{-1}\big(\mu(t)-r\one\big).
\end{equation}

\subsection{Housing tenure, housing services, and owner exposure}

The household can be in one of two tenure states:
\[
    \cI:=\{0,1\},
    \qquad 0=\text{renter},
    \qquad 1=\text{owner}.
\]
Let $I_t\in\cI$ denote the tenure process. In state $i$, the household occupies a
fixed housing-service bundle of size $H_i>0$, with
\[
    0<H_0< H_1.
\]

The strict ordering represents the service advantage of ownership and is used
in the source monotonicity below. Housing services in regime $i$ are
\begin{equation}\label{eq:housing-service}
 h=H_i.
\end{equation}
Thus a price increase changes housing costs without directly changing the
physical service flow.

The renter does not physically hold the house-price asset. The owner does. Let
\begin{equation}\label{eq:q}
    q_0:=0,
    \qquad q_1:=\bar q>0.
\end{equation}

One may set $\bar q=H_1$. We keep owned units $\bar q$ separate from the
utility shifter $H_1$ to distinguish investment exposure from services.

The household pays a tenure-dependent housing flow cost
\begin{equation}\label{eq:flow-cost}
    P_ta_i(t),
    \qquad i=0,1.
\end{equation}
Here $a_0(t)$ represents the rental cost per unit of the house-price index, and
$a_1(t)$ represents maintenance, property taxes, insurance, and other ownership
costs per unit of the house-price index. The functions $a_i$ are deterministic and
bounded.

Switching from state $i$ to state $j\ne i$ at time $t<T$ costs
\begin{equation}\label{eq:switch-cost}
    K_{ij}(t,P_t)=P_tk_{ij}(t).
\end{equation}
Here $k_{01}$ is the buying and moving cost, and $k_{10}$ is the selling and moving
cost. The proportional specification captures value-related fees and taxes
and provides a homogeneous approximation to other moving expenses. It does
not require every component of actual moving costs to be proportional to value.

\begin{assumption}\label{ass:k}
The functions $k_{01},k_{10}\in C^1([0,T])$ are strictly positive. By continuity on the compact horizon,
these functions are bounded and
\[
    \inf_{t\in[0,T]}\{k_{01}(t)+k_{10}(t)\}>0.
\]
Switches are allowed only before $T$. The terminal housing adjustment is handled
separately by the terminal payoff.
\end{assumption}

The positive round-trip cost rules out costless cycling. Unlike the continuous
housing controls in \citet{KraftMunk2011,KraftMunkWagner2018}, the service
bundle here changes only through a tenure switch.

\subsection{Labor income and wealth dynamics}

Labor income is deterministic, time-dependent, and independent of tenure:
\begin{equation}\label{eq:labor}
    \ell\in C([0,T];\R_+).
\end{equation}

Define human wealth by
\begin{equation}\label{eq:human}
    \mathcal H(t):=\int_t^T e^{-r(s-t)}\ell(s)\,\dd s.
\end{equation}
Then $\mathcal H'(t)=r\mathcal H(t)-\ell(t)$ and $\mathcal H(T)=0$.

Let $W_t$ be total net worth, including the market value $q_{I_t}P_t$ of the
physical housing position when the household is an owner. Let
$\pi_t=(\pi_t^1,\ldots,\pi_t^d)^{\tr}$ be the vector of dollar amounts invested in
the traded risky securities. Between switching times, wealth evolves as
\begin{align}\label{eq:W-dyn}
    \dd W_t
    ={}&
    \Big[
        rW_t+\pi_t^{\tr}(\mu(t)-r\one)
        +q_{I_t}P_t\big(\mu_P(t)-r\big)
        -c_t-P_ta_{I_t}(t)+\ell(t)
    \Big]\dd t  \notag\\
    &\quad
    +\Big[\pi_t^{\tr}\Sigma(t)+q_{I_t}P_t\sigma_P(t)^{\tr}\Big]\dd B_t.
\end{align}

The terms involving $q_{I_t}$ record the physical home's return net of its
risk-free financing cost.

Define the traded portfolio that reproduces one unit of house-price diffusion by
\begin{equation}\label{eq:zetaP}
    \zeta_P(t):=\Sigma(t)^{-\tr}\sigma_P(t),
    \qquad
    \zeta_P(t)^{\tr}\Sigma(t)=\sigma_P(t)^{\tr}.
\end{equation}
For a fixed tenure state \(i\), introduce the effective traded position
\begin{equation}\label{eq:pihat}
    \wh\pi_t:=\pi_t+q_iP_t\zeta_P(t).
\end{equation}

\begin{proposition}
\label{prop:neutralization}
Under Assumptions \ref{ass:market} and \ref{ass:portfolio}, between switching
times the wealth equation in regime \(i\) is equivalently
\begin{align}\label{eq:W-neutral}
    \dd W_t
    ={}&
    \Big[
        rW_t+\wh\pi_t^{\tr}(\mu(t)-r\one)
        -c_t-P_t\wh a_i(t)+\ell(t)
    \Big]\dd t
    +\wh\pi_t^{\tr}\Sigma(t)\dd B_t,
\end{align}
where
\begin{equation}\label{eq:kappaH-early}
    \kappa_H(t):=\mu_P(t)-r-\sigma_P(t)^{\tr}\theta(t),
    \qquad
    \wh a_i(t):=a_i(t)-q_i\kappa_H(t).
\end{equation}
Consequently, the owner's physical diffusion does not create a new control or
state variable in the complete-market benchmark.  Its only residual budget effect
is the risk-adjusted carry \(q_i\kappa_H\).
\end{proposition}

\begin{remark}\label{rem:neutrality}
If $P$ is a continuously traded ex-dividend asset with no omitted carry,
no-arbitrage gives $\kappa_H=0$. A physical home may instead have nontraded
service or convenience yield, so spanning its diffusion need not eliminate
the residual carry. Portfolio constraints can also prevent neutralization.
\end{remark}

At a switching time $\tau_n<T$, if the state changes from $\iota_{n-1}$ to
$\iota_n$, then
\begin{equation}\label{eq:W-jump}
    W_{\tau_n}=W_{\tau_n-}-P_{\tau_n}k_{\iota_{n-1},\iota_n}(\tau_n).
\end{equation}

Since $W$ includes physical housing, a purchase exchanges financial for housing
wealth, and only the transaction fee reduces total net worth. The regimes continue
to differ in services, effective flow costs, and switching costs.

\subsection{Preferences, admissibility, and the primal problem}

The household has Cobb--Douglas CRRA utility over perishable consumption and
housing services:
\begin{equation}\label{eq:utility}
    U(c,h)=\frac{c^{1-\gamma_1}h^{\gamma_1-\gamma}}{1-\gamma_1},
    \qquad c>0,\quad h>0,
\end{equation}
where
\begin{equation}\label{eq:gamma1}
    \gamma_1:=1-\alpha(1-\gamma),
    \qquad \alpha\in(0,1),\quad \gamma>0.
\end{equation}
We exclude $\gamma=1$, equivalently $\gamma_1=1$.
Positivity of $\gamma_1=1-\alpha+\alpha\gamma$ is automatic. In tenure state $i$, the housing-service
argument is the fixed physical service level
\[
    h=H_i.
\]

Conditional on tenure, utility is CRRA in consumption with curvature
$\gamma_1$, and prices enter housing outlays rather than the service argument.

At $T$, the household receives terminal utility from total net worth after a
tenure-dependent terminal housing adjustment.  To preserve the same monetary
homogeneity as period consumption utility, terminal wealth has effective curvature
\(\gamma_1\).  Let $\eta>0$, set $\chi_0:=0$, and let the owner's terminal
liquidation parameter satisfy
\[
    0\le \chi_1\le k_{10}(T).
\]
Define
\begin{equation}\label{eq:terminal}
    \Phi_i(w,p)
    =
    \eta\frac{(w-p\chi_i)^{1-\gamma_1}}{1-\gamma_1},
    \qquad i=0,1,
\end{equation}
for $w>p\chi_i$. The cases $\chi_1=0$, $0<\chi_1<k_{10}(T)$, and
$\chi_1=k_{10}(T)$ represent retention, partial adjustment, and full terminal
liquidation. The charge is a resource cost or haircut on net worth, which
already includes the house. The common curvature $\gamma_1$ of monetary
consumption and terminal surplus permits the normalization $z=yp^{\gamma_1}$.

A switching strategy is a sequence
\(\vartheta=((\tau_n,\iota_n))_{n\ge0}\), where \(\tau_0=t\),
\(\iota_0=i\), and each \(\tau_n\) is an \(\F\)-stopping time.  On every
active path the switching times are strictly increasing before the horizon,
\[
    t=\tau_0\le\tau_1<\cdots<\tau_N<T,
    \qquad N<\infty\quad\text{a.s.},
\]
and \(\iota_n\in\{0,1\}\) is \(\mathcal F_{\tau_n}\)-measurable with
\(\iota_n\ne\iota_{n-1}\).  After the last active switch we set all remaining
\(\tau_n=T\) as cemetery times and do not charge a switching cost at \(T\).
An immediate first switch at $t$ is permitted, and $i$ is the tenure just before
that decision. All switching-cost sums run over active indices $n\ge1$,
never the bookkeeping index $n=0$. The tenure process is right-continuous
and satisfies $I_T=I_{T-}$.

A consumption--portfolio pair \((c,\pi)\) consists of a nonnegative progressively
measurable process \(c\) and an \(\R^d\)-valued progressively measurable process
\(\pi\).  We require on \([t,T]\),
\[
 \int_t^T c_s\,\dd s<\infty,
 \qquad
 \int_t^T\left\|\pi_s^{\tr}\Sigma(s)
       +q_{I_s}P_s\sigma_P(s)^{\tr}\right\|^2\dd s<\infty
 \quad\text{a.s.}
\]
For such a control, define the deflated gain process
\[
 M_u:=\int_t^u\Gamma_{t,s}
 \bigl[\Sigma(s)^{\tr}\pi_s+q_{I_{s-}}P_s\sigma_P(s)
                  -W_{s-}\theta(s)\bigr]^{\tr}\,\dd B_s,
 \qquad u\in[t,T],
\]
where $\Gamma_{t,s}$ is the state-price density ratio defined in
\eqref{eq:Gamma}. A control \((c,\pi,\vartheta)\) is admissible if the wealth
equation \eqref{eq:W-dyn}--\eqref{eq:W-jump} is well defined, augmented
wealth satisfies $W_s+\mathcal H(s)>P_s m_{I_s}(s)$ for $s<T$ including
immediately after switches, $M$ is a supermartingale, and
\(W_T-P_T\chi_{I_T}>0\) a.s. Here $m_i(s)$ is the deterministic minimum
remaining housing cost per unit of the current house-price index, defined
in \eqref{eq:minimum-housing-cost} with discount rate
$r_Z=r-\mu_P+\theta^{\tr}\sigma_P=-\kappa_H$. Thus every active switch
must remain financeable after its transaction cost under unconstrained
borrowing. A sufficient condition for the local martingale $M$ to be a
supermartingale is that its negative part is of class~$D$.
For a candidate optimizer we impose the stronger condition that $M$ is a
uniformly integrable martingale, so that the budget inequality is attained
with equality. The admissible set is denoted by
\(\mathfrak A_i(t,w,p)\).

\begin{lemma}\label{lem:no-terminal-switch}
Under \(k_{01}(T)>0\) and \(0\le\chi_1\le k_{10}(T)\), excluding a switch
exactly at \(T\) is without loss of optimality.
\end{lemma}

For initial state $(t,w,p,i)$, define
\begin{align}\label{eq:V-original}
    V_i(t,w,p):=
    \sup_{(c,\pi,\vartheta)\in\mathfrak A_i(t,w,p)}
    \E_{t,w,p,i}\bigg[
        \int_t^T e^{-\beta(s-t)}U(c_s,H_{I_s})\,\dd s
        +e^{-\beta(T-t)}\Phi_{I_T}(W_T,P_T)
    \bigg].
\end{align}
Equivalently, after adding human wealth $x=w+\mathcal H(t)$, we write
\begin{equation}\label{eq:V-augmented}
    V_i(t,w,p)=\wt V_i(t,x,p),
    \qquad x=w+\mathcal H(t).
\end{equation}

\subsection{State prices and convex conjugates}\label{sec:dual}

The reduction uses the martingale approach of
\citet{KaratzasLehoczkyShreve1987,CoxHuang1989} and the finite-horizon
switching framework of \citet{JeonYang}. The additional step here is the
joint homogeneity in the multiplier and the stochastic house-price index.

We now introduce the state-price density used in the martingale reduction. Under
Assumption \ref{ass:market}, the market price of Brownian risk is the process
\(\theta\) defined in \eqref{eq:theta}, and the state-price density
$\Upsilon$ satisfies
\begin{equation}\label{eq:sdf}
    \frac{\dd\Upsilon_t}{\Upsilon_t}
    =-r\,\dd t-\theta(t)^{\tr}\dd B_t,
    \qquad \Upsilon_0=1.
\end{equation}
For $t\le s\le T$, define
\begin{equation}\label{eq:Gamma}
    \Gamma_{t,s}:=\frac{\Upsilon_s}{\Upsilon_t}.
\end{equation}

Recall from \eqref{eq:kappaH-early} the risk-adjusted housing carry
\begin{equation}\label{eq:kappaH}
    \kappa_H(t)=\mu_P(t)-r-\sigma_P(t)^{\tr}\theta(t)
\end{equation}
and the effective net housing cost
\begin{equation}\label{eq:ahat}
    \wh a_i(t)=a_i(t)-q_i\kappa_H(t),
    \qquad i=0,1.
\end{equation}

For $y>0$, define the dual state process
\begin{equation}\label{eq:Y}
    Y_s^{t,y}:=y e^{\beta(s-t)}\Gamma_{t,s},
    \qquad s\in[t,T].
\end{equation}
Then
\begin{equation}\label{eq:Y-dyn}
    \frac{\dd Y_s^{t,y}}{Y_s^{t,y}}
    =
    (\beta-r)\dd s-\theta(s)^{\tr}\dd B_s,
    \qquad Y_t^{t,y}=y.
\end{equation}

For each tenure state $i$, define the convex conjugate
\[
    \wt U_i(y):=\sup_{c>0}\{U(c,H_i)-yc\}.
\]
A direct calculation gives
\begin{equation}\label{eq:Uconj}
    \wt U_i(y)
    =
    \frac{\gamma_1}{1-\gamma_1}
    H_i^{\frac{\gamma_1-\gamma}{\gamma_1}}
    y^{-\frac{1-\gamma_1}{\gamma_1}},
\end{equation}
and the optimizer is
\begin{equation}\label{eq:cstar}
    c_i^*(y)
    =
    H_i^{\frac{\gamma_1-\gamma}{\gamma_1}}
    y^{-1/\gamma_1}.
\end{equation}
The terminal conjugate is
\begin{equation}\label{eq:Phiconj}
    \wt\Phi(y)
    :=
    \sup_{x>0}
    \left\{
        \eta\frac{x^{1-\gamma_1}}{1-\gamma_1}-yx
    \right\}
    =
    \frac{\gamma_1}{1-\gamma_1}\eta^{1/\gamma_1}
    y^{-\frac{1-\gamma_1}{\gamma_1}},
\end{equation}
and the shifted terminal conjugate in state $i$ is
\begin{equation}\label{eq:Phishift}
    \wt\Phi_i(y,p)=\wt\Phi(y)-yp\chi_i.
\end{equation}

Define the normalized dual state
\begin{equation}\label{eq:z}
    z:=yp^{\gamma_1}.
\end{equation}
Then
\begin{align}
    \wt U_i(y)&=p^{1-\gamma_1}\wh u_i(z),\label{eq:homU}\\
    yp\,\wh a_i(t)&=p^{1-\gamma_1}z\,\wh a_i(t),\label{eq:homA}\\
    yp\,k_{ij}(t)&=p^{1-\gamma_1}z\,k_{ij}(t),\label{eq:homK}\\
    \wt\Phi_i(y,p)&=p^{1-\gamma_1}
       \big(\wh\Phi(z)-z\chi_i\big),\label{eq:homPhi}
\end{align}
where
\begin{align}
    \wh u_i(z)
    &:=\frac{\gamma_1}{1-\gamma_1}
    H_i^{\frac{\gamma_1-\gamma}{\gamma_1}}
    z^{-\frac{1-\gamma_1}{\gamma_1}},\label{eq:uhat}\\
    \wh\Phi(z)
    &:=\frac{\gamma_1}{1-\gamma_1}\eta^{1/\gamma_1}
    z^{-\frac{1-\gamma_1}{\gamma_1}}.\label{eq:Phihat}
\end{align}

Spanning and homogeneity play distinct roles: the former permits financial
replication, and the latter removes $p$ from the normalized value. The common
factor $p^{1-\gamma_1}$ uses fixed services, price-proportional monetary
costs, and the same consumption and terminal curvature.

\subsection{Pure finite-horizon switching problem}

Using Lemma \ref{lem:budget} and the Fenchel inequalities, the augmented primal
value satisfies
\begin{equation}\label{eq:duality-bound}
    \wt V_i(t,x,p)\le
    \inf_{y>0}\{\cJ_i(t,y,p)+yx\},
\end{equation}
where the dual value is the pure switching problem
\begin{align}\label{eq:Jyp}
    \cJ_i(t,y,p)
    :=
    \sup_{\vartheta\in\cO_i(t)}
    \E_{t,y,p,i}\bigg[
        &\int_t^T e^{-\beta(s-t)}
        \Big(
            \wt U_{I_s}(Y_s)
            -Y_sP_s\wh a_{I_s}(s)
        \Big)\dd s
        \notag\\
        &+e^{-\beta(T-t)}
        \Big(
            \wt\Phi(Y_T)-Y_TP_T\chi_{I_T}
        \Big)
        \notag\\
        &-\sum_{n\ge1:\,\tau_n<T}
        e^{-\beta(\tau_n-t)}
        Y_{\tau_n}P_{\tau_n}
        k_{\iota_{n-1},\iota_n}(\tau_n)
    \bigg].
\end{align}
Here $\cO_i(t)$ contains all adapted finite-activity switching strategies
starting from tenure state $i$, with the timing conventions of Section~\ref{sec:model}.
No wealth-dependent primal financeability restriction is imposed in this
dual supremum. A strategy with infinite expected nonnegative discounted
switching costs is assigned value $-\infty$. All other payoff terms have
finite expectations by the lognormal moment bounds.

The standard martingale verification argument of \citet{JeonYang} turns
\eqref{eq:duality-bound} into a strong duality criterion, stated as
Proposition~\ref{prop:duality} in the e-companion. A minimizing multiplier
satisfies the first-order condition \eqref{eq:foc}, and
Theorem~\ref{thm:complete-verification} establishes attainment and
financeability on the natural wealth domain.

\subsection{One-dimensional reduction and variational inequalities}

Let
\begin{equation}\label{eq:Z}
    Z_s^{t,z}:=Y_s^{t,y}P_s^{\gamma_1},
    \qquad z=yp^{\gamma_1}.
\end{equation}
By It\^o's formula,
\begin{equation}\label{eq:Zdyn}
    \frac{\dd Z_s}{Z_s}
    =
    b_Z(s)\dd s+\nu_Z(s)^{\tr}\dd B_s,
    \qquad Z_t=z,
\end{equation}
where
\begin{align}
    \nu_Z(s)&:=\gamma_1\sigma_P(s)-\theta(s),\label{eq:nuZ}\\
    b_Z(s)&:=
    \beta-r+\gamma_1\mu_P(s)
    +\frac12\gamma_1(\gamma_1-1)\|\sigma_P(s)\|^2
    -\gamma_1\theta(s)^{\tr}\sigma_P(s).\label{eq:bZ}
\end{align}
Set
\begin{equation}\label{eq:sigmaZ}
    \sigma_Z^2(s):=\|\nu_Z(s)\|^2.
\end{equation}
For the uniform parabolicity, we impose
\begin{equation}\label{eq:uniform-parabolicity}
    \inf_{s\in[0,T]}\sigma_Z(s)>0.
\end{equation}

The homogeneity identities \eqref{eq:homU}--\eqref{eq:homPhi} imply that
\begin{equation}\label{eq:Jhom}
    \cJ_i(t,y,p)=p^{1-\gamma_1}j_i(t,z),
    \qquad z=yp^{\gamma_1}.
\end{equation}
Define the effective discount rate in the normalized dual state by
\begin{equation}\label{eq:betaZ}
    \beta_Z(t)
    :=
    \beta-(1-\gamma_1)\mu_P(t)
    +\frac12\gamma_1(1-\gamma_1)\|\sigma_P(t)\|^2.
\end{equation}
Also define the effective normalized interest rate
\begin{equation}\label{eq:rZ}
    r_Z(t):=
    r-\mu_P(t)+\theta(t)^{\tr}\sigma_P(t).
\end{equation}
Equivalently, by \eqref{eq:kappaH},
\[
    r_Z(t)=-\kappa_H(t).
\]
Thus \(r_Z\) measures the risk-adjusted carrying rate generated by using the
house-price index as the normalization variable. The drift coefficient in the
normalized operator can then be written as
\begin{equation}\label{eq:b}
    \mathfrak b_Z(t):=\beta_Z(t)-r_Z(t)
    =
    \beta-r+\gamma_1\mu_P(t)
    +\frac12\gamma_1(1-\gamma_1)\|\sigma_P(t)\|^2
    -\theta(t)^{\tr}\sigma_P(t).
\end{equation}
The effective one-dimensional operator is
\begin{equation}\label{eq:Lz}
    (\cL_t f)(z)
    :=
    \frac12\sigma_Z^2(t)z^2f_{zz}(z)
    +\big(\beta_Z(t)-r_Z(t)\big)zf_z(z)-\beta_Z(t)f(z).
\end{equation}

The drift $\mathfrak b_Z$ differs from the physical drift $b_Z$ because
applying the joint generator to $p^{1-\gamma_1}j_i(t,yp^{\gamma_1})$
also differentiates the prefactor.

The normalized values satisfy the switching system below. Lemma
\ref{lem:switching-attainment} identifies its solution with
\eqref{eq:Jyp} by reconstructing the individual values from the scalar
obstacle solution and applying switching verification.

\begin{equation}\label{eq:VI}
\begin{cases}
\displaystyle
\max\Big\{
    \partial_tj_0+\cL_tj_0+\wh u_0(z)-\wh a_0(t)z,
    j_1-j_0-k_{01}(t)z
\Big\}=0,\\[0.8em]
\displaystyle
\max\Big\{
    \partial_tj_1+\cL_tj_1+\wh u_1(z)-\wh a_1(t)z,
    j_0-j_1-k_{10}(t)z
\Big\}=0,
\end{cases}
\end{equation}
for $(t,z)\in[0,T)\times(0,\infty)$, with terminal conditions
\begin{equation}\label{eq:jterminal}
    j_i(T,z)=\wh\Phi(z)-\chi_i z,
    \qquad
    \chi_0=0,\quad 0\le \chi_1\le k_{10}(T).
\end{equation}

\section{Timing windows for tenure changes}\label{sec:dob}

This section contains the main analysis. We write the difference of the two
normalized values as a double obstacle problem with time-dependent obstacles.
Its limit for large dual states solves a reflected scalar equation, and this
limit equals the difference between the minimum costs of financing housing in
the two tenures. Theorem~\ref{thm:timing-windows} uses this identity to
determine the dates at which some owner sells and the dates at which every
renter buys. We then describe
the windows by running records and establish the regularity of the
boundaries inside them.

The double-obstacle reduction and the local regularity
arguments build on \citet{JeonYang,HaJeonOk2026}. We do not impose the
global incentive inequalities in Assumption~2.2 of \citet{HaJeonOk2026}.
Instead, the limiting obstacle ODE below determines the active time sets,
including their possible multiple components and the effect of the
regime-dependent terminal datum. Local derivative and boundary estimates
are applied only inside those active components.

The analysis of this section and the policy verification of
Section~\ref{sec:policies} use the following regularity.
\begin{assumption}\label{ass:fb-regularity}
In addition to the standing assumptions above, suppose that
\[
    \sigma_Z^2,\beta_Z,r_Z,a_0,a_1\in C^{0,1}([0,T]),
    \qquad
    k_{01},k_{10}\in C^{1,1}([0,T]).
\]
\end{assumption}
Since $\kappa_H=-r_Z$, its Lipschitz regularity follows automatically.
\subsection{Value-difference formulation}

Define
\begin{equation}\label{eq:Q}
    Q(t,z):=j_0(t,z)-j_1(t,z).
\end{equation}
The switching constraints in \eqref{eq:VI} imply
\begin{equation}\label{eq:obstacles}
    -k_{01}(t)z\le Q(t,z)\le k_{10}(t)z.
\end{equation}
The lower obstacle is reached when a renter optimally buys and becomes an owner.
The upper obstacle is reached when an owner optimally sells and becomes a renter.

Let
\begin{equation}\label{eq:Delta}
    \Delta(t,z):=
    \wh u_0(z)-\wh u_1(z)
    -\big(\wh a_0(t)-\wh a_1(t)\big)z.
\end{equation}
Using the definitions of \(\wh u_i\) and \(\wh a_i\), we can write this source term
more explicitly as
\begin{equation}\label{eq:Delta-expanded}
\begin{aligned}
    \Delta(t,z)
    &=
    \frac{\gamma_1}{1-\gamma_1}
    \left(
        H_0^{\frac{\gamma_1-\gamma}{\gamma_1}}
        -
        H_1^{\frac{\gamma_1-\gamma}{\gamma_1}}
    \right)
    z^{-\frac{1-\gamma_1}{\gamma_1}}  \\
    &\quad
    -
    \left(
        a_0(t)-a_1(t)+q_1\kappa_H(t)
    \right)z .
\end{aligned}
\end{equation}

The source compares the service benefit and effective cost of renting
relative to owning, including the owner's risk-adjusted carry.

Then $Q$ solves
\begin{equation}\label{eq:Qdob}
\begin{cases}
\partial_tQ+\cL_tQ+\Delta(t,z)=0,
& -k_{01}(t)z<Q(t,z)<k_{10}(t)z,\\
\partial_tQ+\cL_tQ+\Delta(t,z)\le0,
& Q(t,z)=-k_{01}(t)z,\\
\partial_tQ+\cL_tQ+\Delta(t,z)\ge0,
& Q(t,z)=k_{10}(t)z,\\
-k_{01}(t)z\le Q(t,z)\le k_{10}(t)z,
& (t,z)\in[0,T)\times(0,\infty),\\
Q(T,z)=\chi_1 z,
& z>0.
\end{cases}
\end{equation}

\subsection{Logarithmic transformation}

Set
\begin{equation}\label{eq:log}
    \tau:=T-t,
    \qquad \xi:=\log z,
    \qquad
    v(\tau,\xi):=e^{-\xi}Q(T-\tau,e^\xi).
\end{equation}
Then
\begin{equation}\label{eq:vobstacles}
    -\psi_0(\tau)\le v(\tau,\xi)\le\psi_1(\tau),
    \qquad
    \psi_0(\tau):=k_{01}(T-\tau),
    \qquad
    \psi_1(\tau):=k_{10}(T-\tau),
\end{equation}
and
\begin{equation}\label{eq:vinitial}
    v(0,\xi)=\chi_1.
\end{equation}
For notational convenience, define the time-reversed coefficients
\begin{equation}\label{eq:time-reversed-coefficients}
\begin{gathered}
    \widetilde\sigma(\tau):=\sigma_Z(T-\tau),
    \qquad
    \widetilde\beta(\tau):=\beta_Z(T-\tau),
    \qquad
    \widetilde r(\tau):=r_Z(T-\tau),
    \\
    \widetilde\kappa_H(\tau):=\kappa_H(T-\tau),
    \qquad
    \widetilde a_i(\tau):=a_i(T-\tau),
    \quad i=0,1.
\end{gathered}
\end{equation}
Moreover, set
\begin{equation}
    \mathfrak a(\tau)
    :=
    \widetilde a_0(\tau)-\widetilde a_1(\tau)
    +q_1\widetilde\kappa_H(\tau).
\end{equation}
Also define
\begin{equation}\label{eq:CH}
    C_H
    :=
    \frac{\gamma_1}{1-\gamma_1}
    \left(
        H_1^{\frac{\gamma_1-\gamma}{\gamma_1}}
        -
        H_0^{\frac{\gamma_1-\gamma}{\gamma_1}}
    \right)
    {>0}.
\end{equation}
The sign follows from $H_1>H_0$ and
$\gamma_1-\gamma=(1-\alpha)(1-\gamma)$,
$1-\gamma_1=\alpha(1-\gamma)$: the exponent in the difference
and the denominator have the same sign. Thus $C_H>0$ for both
$\gamma_1<1$ and $\gamma_1>1$.
Then the source term is given by
\begin{equation}\label{eq:Rsource}
    \begin{aligned}
        \mathcal R(\tau,\xi)
        &:=
        e^{-\xi}
        \bigl(
            \widehat u_0(e^\xi)-\widehat u_1(e^\xi)
        \bigr)
        -\mathfrak a(\tau)\\
        &=
        -C_H e^{-\xi/\gamma_1}
        -\mathfrak a(\tau).
    \end{aligned}
\end{equation}
Then in the domain $\Omega_T:=(0,T)\times \mathbb{R}$,
$v$ satisfies
\begin{equation}\label{eq:vdob}
\begin{cases}
\partial_\tau v-\mathcal L_\xi v=\mathcal R(\tau,\xi),
& -\psi_0(\tau)<v(\tau,\xi)<\psi_1(\tau),\\
\partial_\tau v-\mathcal L_\xi v\ge\mathcal R(\tau,\xi),
& v(\tau,\xi)=-\psi_0(\tau),\\
\partial_\tau v-\mathcal L_\xi v\le\mathcal R(\tau,\xi),
& v(\tau,\xi)=\psi_1(\tau),\\
-\psi_0(\tau)\le v(\tau,\xi)\le\psi_1(\tau),
& (\tau,\xi)\in \Omega_T,\\
v(0,\xi)=\chi_1,
& \xi\in\R.
\end{cases}
\end{equation}
where
\begin{equation}\label{eq:Lx}
    \mathcal L_\xi  v
    :=
    \frac{1}{2}\widetilde\sigma^2(\tau)\partial_{\xi\xi}v
    +
    \left(
        \frac{1}{2}\widetilde\sigma^2(\tau)
        +\widetilde\beta(\tau)-\widetilde r(\tau)
    \right)\partial_\xi v
    -
    \widetilde r(\tau)v.
\end{equation}

The logarithmic normalization makes both obstacles depend only on time.

\subsection{Buying and selling regions}

The buy and sell regions in normalized dual variables are
\begin{align}
    \mathcal B
    &:=\{(t,z):Q(t,z)=-k_{01}(t)z\},\label{eq:buy-region}\\
    \mathcal S
    &:=\{(t,z):Q(t,z)=k_{10}(t)z\}.\label{eq:sell-region}
\end{align}
A renter buys when the lower obstacle is reached, while an owner sells when the
upper obstacle is reached. Through the dual first-order condition \eqref{eq:foc},
these dual boundaries correspond to time-dependent wealth-to-house-price thresholds
in the primal problem.

Lemma~\ref{lem:dob-existence} in the e-companion shows that \eqref{eq:vdob}
has a unique bounded strong solution
$v\in W^{1,2}_{p,\mathrm{loc}}(\Omega_T)\cap C([0,T]\times\R)$ for every
$1<p<\infty$, and Lemma~\ref{lem:bounded-comparison} gives the comparison
principle used throughout. Smooth fit follows from the resulting spatial
regularity.
\begin{proposition}
\label{prop:smooth-fit-complete}
The strong-solution estimates of Lemma~\ref{lem:dob-existence} imply that $Q(t,\cdot)$ is
continuously differentiable for $t<T$. At any finite buy or sell boundary, this spatial $C^1$ regularity,
together with contact with the one-sided obstacle, implies
\[
    Q_z(t,z_b(t))=-k_{01}(t),
    \qquad
    Q_z(t,z_s(t))=k_{10}(t).
\]
Equivalently, for the logarithmically normalized value \(v\), one has
\(v_\xi=0\) at every finite spatial contact point with either time-dependent obstacle.
\end{proposition}

\subsection{Contact regions and the limiting problem}
Assumption~\ref{ass:fb-regularity} makes the transformed coefficients
Lipschitz and the obstacles $C^{1,1}$. Spatial monotonicity first identifies
the contact half-lines, and the limiting cost equation then identifies their
active time domains.

\begin{lemma}\label{lem:state-monotonicity}
Let $v$ be the solution to \eqref{eq:vdob}. Then
\[
    \partial_\xi v\ge0
    \qquad\text{in }\Omega_T.
\]
Moreover, the strict inequality $\partial_\xi v>0$ holds in
$\{-\psi_0<v<\psi_1\}$.
\end{lemma}

Thus lower and upper contact sets are respectively left and right half-lines.
Their locations satisfy the following strict restrictions.

\begin{lemma}\label{lem:contact-location}
For $\tau\in[0,T]$, define
\[
\begin{aligned}
    G_0(\tau)
    &:=
    \psi_0'(\tau)
    +\widetilde r(\tau)\psi_0(\tau)
    -\mathfrak a(\tau),\\
    G_1(\tau)
    &:=
    -\psi_1'(\tau)
    -\widetilde r(\tau)\psi_1(\tau)
    -\mathfrak a(\tau),
\end{aligned}
\]
and, for $i=0,1$,
\[
    \Lambda_i(\tau)
    :=
    \begin{cases}
    \displaystyle
    \gamma_1\log\left(\frac{C_H}{G_i(\tau)}\right),
    &G_i(\tau)>0,\\[1ex]
    +\infty,
    &G_i(\tau)\le0.
    \end{cases}
\]
Then
\[
    \{v=-\psi_0\}
    \subset
    \{(\tau,\xi)\in\Omega_T:\xi<\Lambda_0(\tau)\},
\]
and
\[
    \{v=\psi_1\}
    \subset
    \{(\tau,\xi)\in\Omega_T:
        \xi>\Lambda_1(\tau),\ \Lambda_1(\tau)<\infty\}.
\]
\end{lemma}

The increasing bounded function $v(\tau,\cdot)$ has a limit at infinity.
Its scalar obstacle problem determines the temporal contact structure.

\begin{lemma}\label{lem:limit-v}
Define
\[
    v_\infty(\tau):=\lim_{\xi\to\infty}v(\tau,\xi),
    \qquad 0\le\tau\le T.
\]
Then $v_\infty$ is the unique absolutely continuous solution on $[0,T]$
(with the differential relations understood a.e.) of
\begin{equation}\label{v,infty}
\begin{cases}
-\psi_0\le v_\infty\le\psi_1,\\
v_\infty'+\widetilde r(\tau)v_\infty=-\mathfrak a,
&-\psi_0<v_\infty<\psi_1,\\
v_\infty'+\widetilde r(\tau)v_\infty\ge-\mathfrak a,
&v_\infty=-\psi_0,\\
v_\infty'+\widetilde r(\tau)v_\infty\le-\mathfrak a,
&v_\infty=\psi_1,
\end{cases}
\qquad
v_\infty(0)=\chi_1.    
\end{equation}
In fact $v_\infty\in W^{1,\infty}(0,T)$, and
\[
 \lim_{\xi\to\infty}\sup_{0\le\tau\le T}
       |v(\tau,\xi)-v_\infty(\tau)|=0.
\]
Moreover, denote the limiting contact sets by
\[
    \mathcal C_0(\infty)
    :=
    \{\tau\in(0,T):v_\infty(\tau)=-\psi_0(\tau)\},
    \qquad
    \mathcal C_1(\infty)
    :=
    \{\tau\in(0,T):v_\infty(\tau)=\psi_1(\tau)\}.
\]
Then
\[
    v(\tau,\cdot)\equiv-\psi_0(\tau)
    \qquad\text{for all }\tau\in\mathcal C_0(\infty),
\]
and
\[
    v(\tau,\xi)<\psi_1(\tau)
    \qquad\text{for all }(\tau,\xi)
    \in\bigl((0,T)\setminus\mathcal C_1(\infty)\bigr)\times\R.
\]
\end{lemma}

\subsection{The minimum cost schedule and the timing windows}\label{sec:timing}
By Lemma~\ref{lem:state-monotonicity}, each nonempty contact set is a
half-line, so the free boundaries are the spatial endpoints
\[
    \Xi_0(\tau):=\sup\{\xi\in\R:v(\tau,\xi)=-\psi_0(\tau)\},
    \qquad
    \Xi_1(\tau):=\inf\{\xi\in\R:v(\tau,\xi)=\psi_1(\tau)\},
\]
with the usual conventions for empty and unbounded sets.

\begin{lemma}\label{lem:lower-contact-existence}
For every $\varepsilon\in(0,T)$, there exists
$\Xi_\varepsilon\in\R$ such that
\[
    v(\tau,\xi)=-\psi_0(\tau),
    \qquad
    (\tau,\xi)\in[\varepsilon,T]\times(-\infty,\Xi_\varepsilon].
\]
\end{lemma}

We next identify the times at which the lower boundary is finite.

\begin{proposition}\label{prop:lower-contact-times}
Let
\[
    \mathcal T_0
    :=
    \{\tau\in(0,T):\Xi_0(\tau)\in\R\}.
\]
Then
\[
    \mathcal T_0
    =
    (0,T)\setminus\mathcal C_0(\infty)
    =
    \{\tau\in(0,T):v_\infty(\tau)>-\psi_0(\tau)\}.
\]
\end{proposition}

We now turn to the upper boundary. Since the lower contact region is
nonempty at every positive time and $\psi_0+\psi_1>0$, one has
$\Xi_1(\tau)>-\infty$ for every $\tau\in(0,T)$.
Lemma~\ref{lem:finite-contact-interior} in the e-companion shows that
contact with the upper obstacle at a finite spatial point occurs only at an
interior point of $\mathcal C_1(\infty)$, and
Lemma~\ref{lem:recover-upper-contact} gives the converse at interior times
with $\Lambda_1<\infty$. Hence the upper boundary is finite exactly on the
following set.
\begin{proposition}\label{prop:upper-contact-times}
Let
\[
    \mathcal T_1
    :=
    \{\tau\in(0,T):\Xi_1(\tau)\in\R\}.
\]
Then
\[
    \mathcal T_1
    =
    \operatorname{int}\mathcal C_1(\infty)
    \cap
    \{\tau\in(0,T):\Lambda_1(\tau)<\infty\}.
\]
\end{proposition}

The limiting function $v_\infty$ has a direct economic meaning, which we now
make precise.

Write $g=\gamma_1$. Define the minimum cost of
financing housing, excluding perishable consumption and positive terminal
surplus, by
\begin{align}\label{eq:minimum-housing-cost}
 m_i(t)=\inf_{\zeta}\bigg\{
 \int_t^T D_H(t,u)\widehat a_{I_u}(u)\,du
 +D_H(t,T)\chi_{I_T}
 &+\sum_{n}D_H(t,\tau_n)k_{\iota_{n-1},\iota_n}(\tau_n)\bigg\},\notag\\
 D_H(t,u)&=\exp\left(-\int_t^u r_Z(q)\,dq\right).
\end{align}
Here the infimum is over deterministic finite switching schedules, including
an immediate switch at $t$, and excludes a switch at $T$.
All sums are over active switches with index $n\ge1$.
The natural interior domain of the unconstrained complete-market problem is
\begin{equation}\label{eq:financeability-domain}
 \mathcal D_i=\{(t,x,p):t<T,\ p>0,\ x>p m_i(t)\}.
\end{equation}
In particular, $x>0$ alone is not a sufficient solvency condition when
mandatory housing outlays are positive. The quantity $m_i$ can be negative
when the net housing carry is sufficiently favorable.
\begin{proposition}\label{prop:minimum-financing}
The deterministic infimum in \eqref{eq:minimum-housing-cost} is attained,
and $m_i$ is bounded and continuous. Among adapted finite switching
policies with integrable state-price costs, the infimum of the expected
remaining housing cost is $p m_i(t)$. This bound also holds conditionally
at stopping times. Moreover,
\[
 v_\infty(T-t)=m_1(t)-m_0(t).
\]
\end{proposition}

The quantity $m_i(t)$ is the smallest cost, per unit of the current
house-price index, of paying for housing from $t$ to $T$ when the household
starts in tenure $i$ and follows a deterministic schedule. Costs are
discounted at $r_Z=-\kappa_H$, the rate that applies when the house-price
index is used as numeraire (see the proof of
Proposition~\ref{prop:minimum-financing}). Since an immediate switch is
allowed in \eqref{eq:minimum-housing-cost},
\begin{equation}\label{eq:immediate-switch-bounds}
 m_0(t)\le m_1(t)+k_{01}(t),
 \qquad
 m_1(t)\le m_0(t)+k_{10}(t),
\end{equation}
with equality exactly when some minimum cost schedule switches at $t$. The
minimizing schedule need not be unique.
Combining Proposition~\ref{prop:minimum-financing} with
Propositions~\ref{prop:lower-contact-times}
and~\ref{prop:upper-contact-times} gives the main result of the paper.

\begin{theorem}\label{thm:timing-windows}
Let $t\in(0,T)$ and $\tau=T-t$.
\begin{enumerate}[label=(\roman*),leftmargin=*]
\item Every renter with $(t,x,p)\in\mathcal D_0$ buys at $t$ if and only if
$m_0(t)=m_1(t)+k_{01}(t)$, that is, if and only if a minimum cost schedule
starting from renting buys at $t$. At every other $t\in(0,T)$, the buying
region is nonempty and consists of the renters whose wealth lies above a
finite threshold.
\item Some owner with $(t,x,p)\in\mathcal D_1$ sells at $t$ if and only if
$\tau$ is an interior point of
$\{u\in(0,T):m_1(T-u)=m_0(T-u)+k_{10}(T-u)\}$ and $G_1(\tau)>0$. In this
case the selling region consists of the owners whose wealth lies between the
financing floor $p\,m_1(t)$ and a finite threshold. Otherwise no owner sells
at $t$.
\item The sets of dates described in (i) and (ii) are determined by $r_Z$,
$a_0$, $a_1$, $q_1$, $k_{01}$, $k_{10}$ and $\chi_1$. They are the same for
all values of $\beta$, $\gamma$, $\alpha$, $H_0<H_1$, $\eta$ and $\ell$ that
satisfy the standing assumptions.
\end{enumerate}
\end{theorem}

\begin{proof}
By Proposition~\ref{prop:minimum-financing},
$v_\infty(\tau)=m_1(t)-m_0(t)$, while $\psi_0(\tau)=k_{01}(t)$ and
$\psi_1(\tau)=k_{10}(t)$. Hence $\tau\in\mathcal C_0(\infty)$ if and only if
$m_0(t)=m_1(t)+k_{01}(t)$, and $\tau\in\mathcal C_1(\infty)$ if and only if
$m_1(t)=m_0(t)+k_{10}(t)$. Lemma~\ref{lem:limit-v} gives
$v(\tau,\cdot)\equiv-\psi_0(\tau)$ on $\mathcal C_0(\infty)$.
Proposition~\ref{prop:lower-contact-times} and
Lemma~\ref{lem:lower-contact-existence} show that at every other positive
time the lower contact set is a nonempty half-line with a finite endpoint.
Proposition~\ref{prop:upper-contact-times} identifies the times with finite
upper contact, and $\Lambda_1(\tau)<\infty$ if and only if $G_1(\tau)>0$.
Proposition~\ref{prop:wealth-coordinates} maps the dual contact half-lines
onto the stated wealth regions. It maps an all-buy slice onto the entire
admissible renter domain and leaves an absent selling region absent. This
proves (i) and (ii). For (iii), the functions $m_i$ in
\eqref{eq:minimum-housing-cost} involve only $r_Z$,
$\widehat a_i=a_i+q_ir_Z$, $k_{01}$, $k_{10}$ and $\chi_1$, and
$G_1=-\psi_1'-\widetilde r\psi_1-\mathfrak a$ involves only $r_Z$, $k_{10}$
and $\mathfrak a$, which is the time reversal of $a_0-a_1-q_1r_Z$.
\end{proof}

Theorem~\ref{thm:timing-windows} separates two questions. Whether some
owner sells, or every renter buys, at a given date is answered by a
deterministic minimum cost problem. Which households make the change is answered by the full
problem, through thresholds that depend on preferences, services, and labor
income. The service advantage $H_1>H_0$ is what makes the second answer
depend on wealth. If the service levels were equal, then $C_H=0$, the source
$\mathcal R=-\mathfrak a$ would not depend on $\xi$, and the unique bounded
solution of \eqref{eq:vdob} would be $v\equiv v_\infty$. A deterministic
minimum cost schedule would then be optimal at every admissible wealth level.
The assumption $H_0<H_1$
excludes this degenerate case, but it shows that wealth matters for tenure
changes only because the two tenures provide different services.

The theorem also shows which windows can arise. Away from the obstacles,
\eqref{v,infty} reads $v_\infty'=-\widetilde r v_\infty-\mathfrak a$, and
$v_\infty(0)=\chi_1\ge0$. On an interval of remaining time on which
$\mathfrak a\le0$, the function $v_\infty$ therefore cannot become negative,
and it never meets the negative lower obstacle $-\psi_0$. In calendar time,
$\mathfrak a\le0$ means $d:=\widehat a_1-\widehat a_0\ge0$ with
$\widehat a_i=a_i+q_ir_Z$, that is, owning costs at least as much as renting
once the risk-adjusted carry is included. A date $t$ can thus be one at which
every renter buys only if owning is strictly cheaper than renting at some
date in $[t,T]$. When owning costs at least as much as renting throughout, as
in the baseline of Section~\ref{sec:numerical}, the part of
Theorem~\ref{thm:timing-windows} that does not depend on preferences
concerns selling.

For applications, the windows created by a cost schedule can therefore be
evaluated with the scalar equation \eqref{v,infty} alone, as in
Table~\ref{tab:components} below.

\subsection{Running records}\label{sec:running-records}

The scalar problem \eqref{v,infty} also yields explicit conditions for its
contact times, which describe the minimum cost schedule directly. Let $\vns$
solve
\[
    \vns'(\tau)+\widetilde r(\tau)\vns(\tau)=-\mathfrak a(\tau),
    \qquad \vns(0)=\chi_1.
\]
Equivalently,
\[
    \vns(\tau)
    =
    \exp\left(-\int_0^\tau\widetilde r(q)\,dq\right)\chi_1
    -
    \int_0^\tau
    \exp\left(-\int_s^\tau\widetilde r(q)\,dq\right)
    \mathfrak a(s)\,ds.
\]
Since $\tau=T-t$, this representation expresses the monetary-value
difference between the two states, under no further switching, as the
discounted value accumulated from the current time $t$ to the terminal
time $T$, together with the discounted terminal difference $\chi_1$.

Define the transformed switching gains
\[
\begin{aligned}
    F_1(\tau)
    &:=
    \exp\left(\int_0^\tau\widetilde r(q)\,dq\right)
    \bigl(\vns(\tau)-\psi_1(\tau)\bigr),\\
    F_0(\tau)
    &:=
    \exp\left(\int_0^\tau\widetilde r(q)\,dq\right)
    \bigl(-\vns(\tau)-\psi_0(\tau)\bigr).
\end{aligned}
\]

For $\tau\in(0,T)$, define
\begin{equation}\label{eq:ell-contact}
    \ell_i(\tau)
    :=
    \sup\bigl(
        \{0\}\cup
        ({\mathcal C_i(\infty)}\cap[0,\tau))
    \bigr),
    \qquad i=0,1.
\end{equation}
Thus $\ell_i(\tau)$ records the most recent limiting contact with the
corresponding obstacle below $\tau$. If no such contact exists, we set
$\ell_i(\tau)=0$ and, by convention, interpret
$F_i(\ell_i(\tau))$ as zero.

\begin{lemma}\label{lem:running-record}
For every $\tau\in(0,T)$,
    $\tau\in\mathcal C_1(\infty)$
if and only if the following two conditions hold:
\[
\begin{aligned}
    &\text{\rm (i)}\quad
    \begin{cases}
        F_1(\tau)+F_0(\ell_0(\tau))\ge0,
        &\ell_0(\tau)>0,\\
        F_1(\tau)\ge0,
        &\ell_0(\tau)=0,
    \end{cases}\\
    &\text{\rm (ii)}\quad
    F_1(\tau)
    =
    \max_{\ell_0(\tau)\le s\le\tau}F_1(s).
\end{aligned}
\]
Moreover,
\[
    F_1'(\tau)\ge0
    \qquad
    \text{for a.e. }\tau\in\mathcal C_1(\infty).
\]

Similarly, for every $\tau\in(0,T)$,
    $\tau\in\mathcal C_0(\infty)$
if and only if
\[
\begin{aligned}
    &\text{\rm (i)}\quad
    \begin{cases}
        F_0(\tau)+F_1(\ell_1(\tau))\ge0,
        &\ell_1(\tau)>0,\\
        F_0(\tau)\ge0,
        &\ell_1(\tau)=0,
    \end{cases}\\
    &\text{\rm (ii)}\quad
    F_0(\tau)
    =
    \max_{\ell_1(\tau)\le s\le\tau}F_0(s).
\end{aligned}
\]
Moreover,
\[
    F_0'(\tau)\ge0
    \qquad
    \text{for a.e. }\tau\in\mathcal C_0(\infty).
\]
\end{lemma}

The function $\vns(\tau)$ is the difference between the costs of owning and
of renting from date $T-\tau$ to $T$ without further switches, in units of
the current house-price index. Hence $F_1(\tau)$ is the cost saving from selling
at date $T-\tau$ and renting until $T$, compounded to the horizon at the rate
$r_Z$, and $F_0(\tau)$ is the corresponding saving from buying. Since
$s<\tau$ corresponds to a later calendar date, condition (ii) states that no
later date before the next all-buy date offers a larger saving from selling.
Condition (i) states that selling now and buying back at that all-buy date
does not lose money. The lower conditions have the same reading for
purchases. The running records thus describe the minimum cost schedule
explicitly. Corollary~\ref{cor:lower-time-characterization} and
Lemma~\ref{lem:upper-time-characterization} in the e-companion restate both
parts of Theorem~\ref{thm:timing-windows} in these terms. The lower
conditions are weak and the upper conditions are strict. This asymmetry
reflects the service advantage of ownership: limiting upper contact alone
does not make selling optimal, and a finite selling boundary also requires
$F_1'(\tau)>0$, that is, $G_1(\tau)>0$.

\subsection{Regularity and endpoint behavior}\label{sec:regularity}

Local derivative estimates and non-degeneracy make $\Xi_0$ and $\Xi_1$ locally
Lipschitz on $\mathcal T_0$ and $\mathcal T_1$, and a boundary Harnack
argument gives $\Xi_0\in C^{1,\alpha}_{\rm loc}(\mathcal T_0)$ and
$\Xi_1\in C^{1,\alpha}_{\rm loc}(\mathcal T_1)$ for every $\alpha\in(0,1/2)$
(Lemma~\ref{lem:boundary-C1alpha}). The lower boundary satisfies
$\Xi_0(\tau)\to-\infty$ as $\tau\downarrow0$, and $\Xi_0$ and $\Xi_1$ tend to
$+\infty$ at every endpoint of $\mathcal T_0$ and $\mathcal T_1$ in $(0,T)$
(Lemma~\ref{lem:boundary-behavior-free-boundaries}), so their graphs are the
entire free boundaries (Corollary~\ref{cor:free-boundary-graphs}). Near the
horizon, a finite selling boundary exists if $\psi_1(0)=\chi_1$ and
$G_1(0)>0$, and then $\Xi_1(\tau)\to\Lambda_1(0)$ as $\tau\downarrow0$
(Lemma~\ref{lem:initial-upper-boundary}). Selling is absent near the
horizon if $\chi_1<\psi_1(0)$, or if $\chi_1=\psi_1(0)$ and $G_1(0)<0$. If
$\chi_1=\psi_1(0)$ and $G_1(0)=0$, either selling is absent near the horizon
or the selling boundary diverges there
(Lemma~\ref{lem:remaining-terminal-cases}).

\section{Optimal housing and financial policies}\label{sec:policies}
\subsection{Strong duality and attainment}\label{sec:verification-complete}

The minimum costs $m_i$ of Section~\ref{sec:timing} are also the financing
floors of the primal problem. The dual values have the following regularity
and range.
\begin{lemma}\label{lem:switching-attainment}
The dual switching values are attained. For $t<T$, each
$\mathcal J_i(t,\cdot,p)$ is continuously differentiable and strictly
convex, with
\begin{equation}\label{eq:dual-derivative-range}
 \lim_{y\downarrow0}\mathcal J_{i,y}(t,y,p)=-\infty,
 \qquad
 \lim_{y\uparrow\infty}\mathcal J_{i,y}(t,y,p)=-p m_i(t).
\end{equation}
\end{lemma}

\begin{theorem}\label{thm:complete-verification}
Under the complete-market assumptions and Assumption
\ref{ass:fb-regularity}, for every $(t,x,p)\in\mathcal D_i$ there is a
unique multiplier $y^*>0$ satisfying
\begin{equation}\label{eq:unique-multiplier}
 x=-\mathcal J_{i,y}(t,y^*,p).
\end{equation}
The contact switching strategy $\vartheta^*$, the consumption and terminal
wealth in Proposition~\ref{prop:duality}, and a financing traded portfolio
are admissible and optimal. In particular,
\[
 \widetilde V_i(t,x,p)
 =\min_{y>0}\{\mathcal J_i(t,y,p)+yx\},
 \qquad
 V_i(t,w,p)=\widetilde V_i(t,w+\mathcal H(t),p).
\]
This equality concerns the original finite-activity switching problem,
and no randomized or relaxed tenure choice is introduced.
\end{theorem}

\subsection{Wealth thresholds, consumption and financial hedging}
\label{sec:primal-policies}
To express the optimal policies, use the positive dual coordinate
\begin{equation}\label{eq:s-coordinate}
 s=C_H z^{-1/g},\qquad
 h_i(\tau,s)=j_i(T-\tau,z)/z,\qquad
 w_i(\tau,s)=-h_i(\tau,s)+\frac{s}{g}h_{i,s}(\tau,s).
\end{equation}
Then $X/P=w_i(\tau,s)$ and $h_0-h_1=v(\tau,\log z)$.
The variables $z$ and $s$ are dual coordinates, and neither is household wealth.
Strict convexity of $j_i$ makes $w_i$ strictly increasing in $s$,
with range $(m_i(T-\tau),\infty)$. Thus inversion of this single
monotone function determines the optimal multiplier from current wealth.
Writing $A_i=H_i^{(g-\gamma)/g}$, optimal consumption is
\begin{equation}\label{eq:primal-consumption-ratio}
 \frac{c^*}{P}=\frac{A_i}{C_H}s.
\end{equation}
Define $b(\tau)=C_He^{-\Xi_0(\tau)/g}$ and
$a(\tau)=C_He^{-\Xi_1(\tau)/g}$, using $b=0$ on
$\mathcal C_0(\infty)$ and $a=0$ outside $\mathcal T_1$.
Here $a(\tau),b(\tau)$ are free boundaries, while $a_i(t)$ denotes a flow
cost and $\mathfrak b_Z(t)$ the normalized generator drift.
Buying occurs for $s\ge b$, and selling occurs for $0<s\le a$.
The pre-switch wealth thresholds, whenever finite, are
\[
 \frac{w_{\rm buy}^{\rm net}+\mathcal H(t)}{P}=w_0(T-t,b(T-t)),
 \qquad
 \frac{w_{\rm sell}^{\rm net}+\mathcal H(t)}{P}=w_1(T-t,a(T-t)).
\]
The thresholds use their respective tenure values. On each whole
contact region the corresponding value-matching identity holds
identically. Differentiating there and extending to the boundary by
continuity of $j_{i,z}$ gives the financing identities
\begin{equation}\label{eq:wealth-jump-identities}
 w_0(\tau,b)-w_1(\tau,b)=\psi_0(\tau),\qquad
 w_1(\tau,a)-w_0(\tau,a)=\psi_1(\tau).
\end{equation}
Whenever both boundaries are finite and positive, $0<a<b$ and the
strict increase of $w_1(\tau,\cdot)$ imply
\begin{equation}\label{eq:wealth-threshold-order}
 w_1(\tau,a)<w_1(\tau,b)
       =w_0(\tau,b)-\psi_0(\tau)<w_0(\tau,b).
\end{equation}
Thus the pre-purchase wealth threshold exceeds the pre-sale threshold
by strictly more than the buying cost. In net-worth units,
$w_{\rm buy}^{\rm net}-w_{\rm sell}^{\rm net}>Pk_{01}(t)$.
This ordering holds for every parameter set satisfying the assumptions,
but it does not imply time monotonicity or a curvature restriction.
The shadow price is continuous through a switch, while wealth loses the
transaction cost and the consumption coefficient changes from $A_i$ to $A_j$.

At continuation points where the indicated derivatives exist, the financial
portfolio has the explicit form
\begin{equation}\label{eq:portfolio-recovery}
 \frac{\pi_i^*}{P}
 =\Sigma^{-\mathsf T}\left[
    (w_i-q_i)\sigma_P
   +s w_{i,s}\left(\frac{\theta}{g}-\sigma_P\right)\right],
 \quad
 w_{i,s}=-\left(1-\frac1g\right)h_{i,s}+\frac{s}{g}h_{i,ss}.
\end{equation}
The formula holds almost everywhere in time,
which suffices for stochastic integration, and the martingale representation in
\eqref{eq:wealth-replication} defines the portfolio everywhere up to the
usual equivalence of stochastic integrands. Complete spanning therefore
removes unhedgeable housing risk, but it does not eliminate the household's
financial hedge or the effective owner cost.

\begin{proposition}\label{prop:wealth-coordinates}
For each $t<T$, the regime-specific map
$z\mapsto-j_{i,z}(t,z)$ is continuous and strictly decreasing from
$+\infty$ to $m_i(t)$. It maps the dual buy and sell regions to their
pre-switch wealth regions. Finite buy and sell boundaries are evaluated
with $j_0$ and $j_1$, respectively. If a log-dual boundary tends to
$+\infty$ at an internal endpoint, its wealth threshold tends to the
corresponding financing floor. Thus the active time components are
preserved under the regime-specific coordinate changes, but time
monotonicity and curvature of the boundary graphs need not be preserved.
Whenever both boundaries are finite, their pre-switch wealth levels
satisfy the strict separation \eqref{eq:wealth-threshold-order}.
The consumption, wealth-jump and portfolio identities above hold for the
optimal policy.
\end{proposition}

\section{Numerical results}\label{sec:numerical}

The calculations answer three questions: at which dates some owner sells or
every renter buys, what wealth triggers a change, and what consumption,
financial hedge, and flexibility value follow. By Theorem~\ref{thm:timing-windows}, the first is
answered by the scalar equation \eqref{v,infty} alone. The coupled integral
equations then locate the thresholds inside the windows, and these are
converted to wealth with the value of the current tenure. An absent selling
region is therefore never treated as a small positive sale threshold.
\subsection{Parameters}
All rates are annual and $T=20$ years. Table~\ref{tab:parameters} gives
the baseline primitives. These are illustrative, economically interpretable
parameters, not estimates of a particular housing market. The house-price
index is the value of one standard property, and $H_0=0.8$ and $H_1=1$
represent a smaller rental service bundle and a larger owned bundle.
The renter pays 4.5 percent of the property value annually. Owner cash
expenses are 3.2 percent, while the risk-adjusted carrying cost adds
2.5 percent, producing an effective ownership premium of 1.2 percentage
points. All transaction costs are strictly positive, including at $T$.

\begin{table}[htbp]\centering\small
\caption{Baseline parameters.}\label{tab:parameters}
\begin{tabular}{llr@{\qquad}llr}
\toprule
Quantity & Symbol & Value & Quantity & Symbol & Value\\
\midrule
Horizon & $T$ &20 & Interest rate &$r$&0.030\\
Utility discount &$\beta$&0.035 & Consumption share &$\alpha$&0.50\\
Composite curvature &$\gamma$&3 & Effective curvature &$g$&2\\
Rental services &$H_0$&0.8 & Owned services / units &$H_1,q_1$&1\\
Stock volatility &$\sigma_S$&0.20 & Market price of risk &$\theta$&0.40\\
House volatility &$\sigma_P$&0.10 & House-price drift &$\mu_P$&0.045\\
Rental flow cost &$a_0$&0.045 & Owner cash expenses &$a_1$&0.032\\
Buying cost &$k_{01}$&0.05 & Selling cost &$k_{10}$&0.04\\
Terminal adjustment &$\chi_1$&0.015 & Terminal utility weight &$\eta$&1\\
\bottomrule
\end{tabular}
\end{table}

The implied quantities are
\[
 C_H=0.23606798,\quad r_Z=0.025,\quad\beta_Z=0.070,
 \quad\sigma_Z=0.20,\quad\varsigma=0.10,\quad b_s=-0.0275.
\]
All results use the coupled Gaussian integral equations of
Section~\ref{sec:integral-method} of the e-companion, which defines the
coefficients $\varsigma$ and $b_s$ of the integral representation and
reports independent numerical checks.

\subsection{Timing windows}
Let $\mathcal G_{c,w}(\tau)=\exp(-((\tau-c)/w)^2/2)$ and
$d=\widehat a_1-\widehat a_0=-\mathfrak a$.
Eleven scenarios illustrate the contact and endpoint mechanisms
of Section~\ref{sec:dob}. Unspecified parameters are those in
Table~\ref{tab:parameters}. When $d$ varies, keep $a_1=0.032$ and set
$a_0(\tau)=0.032+0.025-d(\tau)$, so all rental yields and transaction costs
remain positive. Anticipated taxes and changing relative housing expenses
motivate the time-varying paths
\citep{BestKleven2018,HilberLyytikainen2017,RosenRosen1980}.
Cases 8, 10 and 11 isolate critical or equality
configurations and are described in Section~\ref{ec:cases} of the
e-companion. The theory does not restrict the number of active components.

\begin{enumerate}[leftmargin=*,itemsep=4pt]
\item \textbf{Partial terminal adjustment.}
The baseline has $d=0.012$ and $\chi_1=0.015<k_{10}(T)$.
Retaining the property avoids part of the sale cost, and selling is unavailable
during the last 2.2105 years. Allowing terminal retention is motivated by
late-life ownership persistence \citep{VentiWise2004}.

\item \textbf{Full terminal liquidation.}
Set $\chi_1=k_{10}(T)=0.04$. Waiting until $T$ no longer avoids the
selling cost. The positive effective ownership premium then permits a
finite sale threshold arbitrarily close to $T$.

\item \textbf{No selling during the horizon.}
Set $d=0.0005$ and $\chi_1=0$. The small ownership cost premium never
justifies paying to leave the superior service bundle, so the selling region
is empty throughout the horizon.

\item \textbf{An all-buy interval.}
Set $d=-0.008$ and $\chi_1=0.015$. Ownership offers better services and
lower effective expenses. With at least 8.6283 years remaining, every
financeable renter buys. The role of expected residence length is consistent
with the mechanism studied by \citet{HaurinGill2002}.

\item \textbf{One interruption of selling.}
Set $\psi_1=0.04+0.05\mathcal G_{7,1.2}$ and $\chi_1=\psi_1(0)$.
The anticipated selling-cost surcharge removes the finite sale boundary
between 4.5686 and 7.7990 remaining years.

\item \textbf{Two interruptions of selling.}
Set $\psi_1=0.04+0.05(\mathcal G_{6,1}+\mathcal G_{13,1})$ and
$\chi_1=\psi_1(0)$. Two anticipated surcharges produce three separate
intervals with a finite selling boundary.

\item \textbf{Alternating relative housing costs.}
Set $d(\tau)=0.028\cos(\pi\tau/6)$ and $\chi_1=0.04$.
With owner expenses fixed, the rental yield varies so that the effective
cost advantage alternates. All-buy intervals and finite-selling intervals
then alternate. An individual trades only when its state reaches contact.
\setcounter{enumi}{8}
\item \textbf{A negative terminal sale incentive.}
Set $d(\tau)=-0.004+0.022(1-e^{-\tau/3})$ and $\chi_1=0.04$.
Despite full terminal liquidation, the local advantage of retaining ownership
removes the sale boundary during the last 1.6079 years.
\end{enumerate}

\begin{table}[htbp]\centering\small
\caption{Contact times from the limiting ODE, in remaining years.
Intervals are intersected with $(0,20)$, and isolated times are shown in braces.}
\label{tab:components}
\begin{tabular}{cll}
\toprule
Case & Finite selling boundary $\mathcal T_1$ & All-buy times $\mathcal C_0(\infty)$\\
\midrule
1 & $(2.210507,20)$ & $\varnothing$\\
2 & $(0,20)$ & $\varnothing$\\
3 & $\varnothing$ & $\varnothing$\\
4 & $\varnothing$ & $[8.628343,20)$\\
5 & $(0,4.568620)\cup(7.799007,20)$ & $\varnothing$\\
6 & $(0,3.861939)\cup(6.810759,10.862036)$ & $\varnothing$\\
  & $\quad\cup(13.810750,20)$ & \\
7 & $(0,2.931776)\cup(13.367584,14.931776)$
  & $[7.421167,8.914710]\cup[19.421167,20)$\\
8 & $(0,20)$ & $\varnothing$\\
9 & $(1.607891,20)$ & $\varnothing$\\
10 & $(0,10)\cup(10,20)$ & $\varnothing$\\
11 & $\varnothing$ & $\{10\}$\\
\bottomrule
\end{tabular}
\end{table}

\begin{figure}[p]\centering
\includegraphics[width=\textwidth,height=.82\textheight,keepaspectratio]{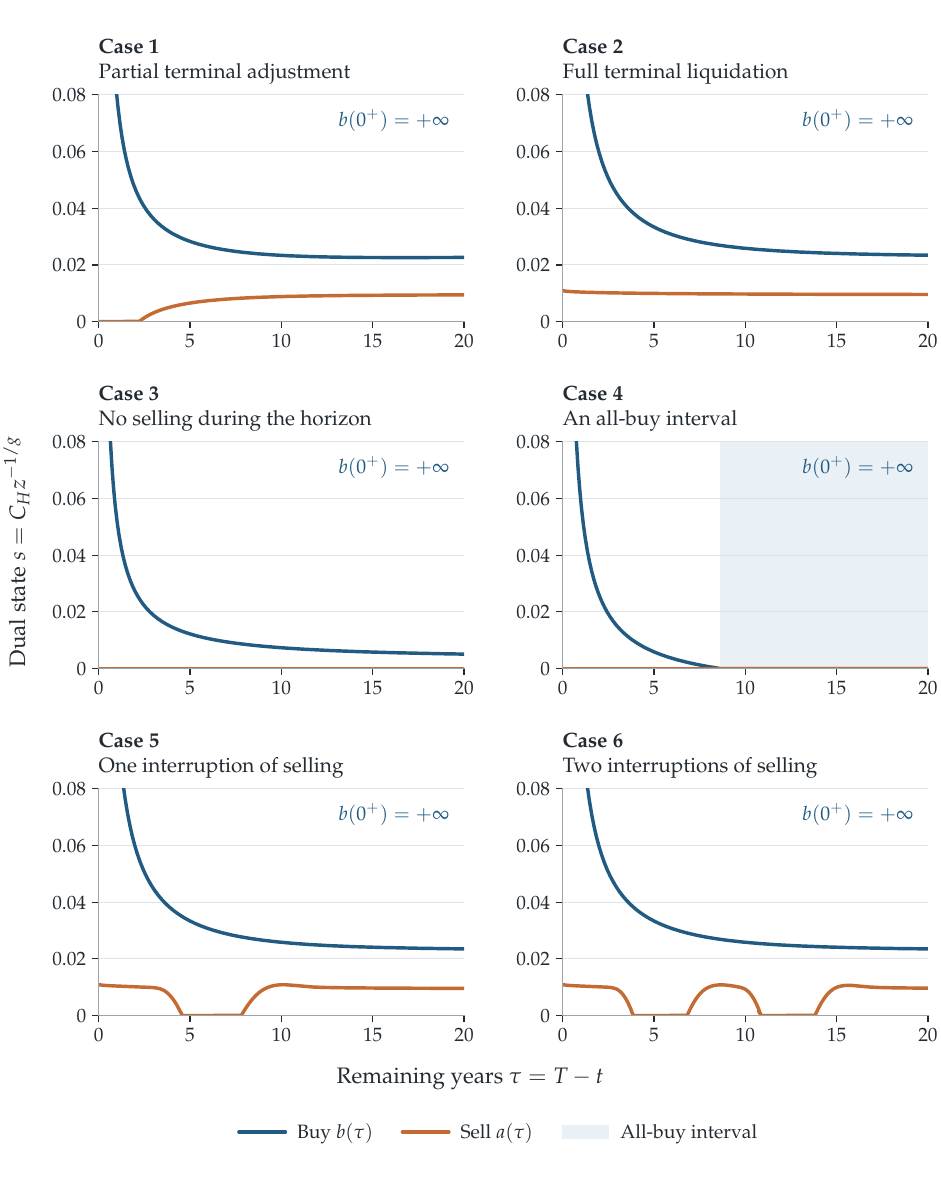}
\caption{Cases 1--6. Buying occurs above the blue curve $b$ and selling
below the orange curve $a$. The vertical scale is the transformed state
$s=C_Hz^{-1/g}$. A zero selling curve means that selling is unavailable
at every positive state. Shading denotes an all-buy interval. The
divergent terminal buy threshold is indicated explicitly and exceeds the
display range close to zero remaining time.}\label{fig:cases1}
\end{figure}

Table~\ref{tab:components} lists the active time components, and
Figure~\ref{fig:cases1} shows the dual boundaries for Cases 1--6
(Figure~\ref{fig:cases2} of the e-companion shows Cases 7--11). All entries
of the table come from the scalar equation \eqref{v,infty} and, by
Theorem~\ref{thm:timing-windows}, are the same for every admissible choice
of $\beta$, $\gamma$, $\alpha$, $H_0<H_1$, $\eta$ and $\ell$.
Case 1 is the baseline. Owning costs more than renting, $d=0.012$, so no
date is an all-buy date and the timing windows consist of the selling
window alone. Whatever the preferences, services, and income, no owner sells during
the last 2.21 years, a length given in closed form by
\eqref{eq:sale-onset-economic} below. Section~\ref{ec:robust} of the
e-companion recomputes the thresholds when short sales or borrowing are
restricted and when the ownership cost is subject to an uninsurable
increase.

Case 5 is easiest to read in calendar time, from right to left in
Figure~\ref{fig:cases1}. The surcharge peaks seven years before the horizon.
Selling regions exist until 7.80 years remain, about 0.8 years before the
peak, when the selling cost has already risen to about 8 percent of the
property value. They are absent through the peak and reappear only when 4.57
years remain, 2.4 years after the peak, when the surcharge has fallen to
about 0.6 percent. The two ends of the interruption have different causes.
The window closes when a sale after the surcharge first offers a larger
saving than a sale now, which is the running-record condition of
Lemma~\ref{lem:running-record}. It reopens when $G_1$ becomes positive
again, that is, when the selling cost declines more slowly than the net flow
gain from selling, $d-r_Z\psi_1$, accrues. Case 6 shows the same mechanism
twice. In both cases a favorable instantaneous sale incentive does not by
itself restore selling.

\subsection{Representative boundaries in household wealth}
\label{sec:wealth-cases}
\begin{figure}[tbp]\centering
\includegraphics[width=.86\textwidth]{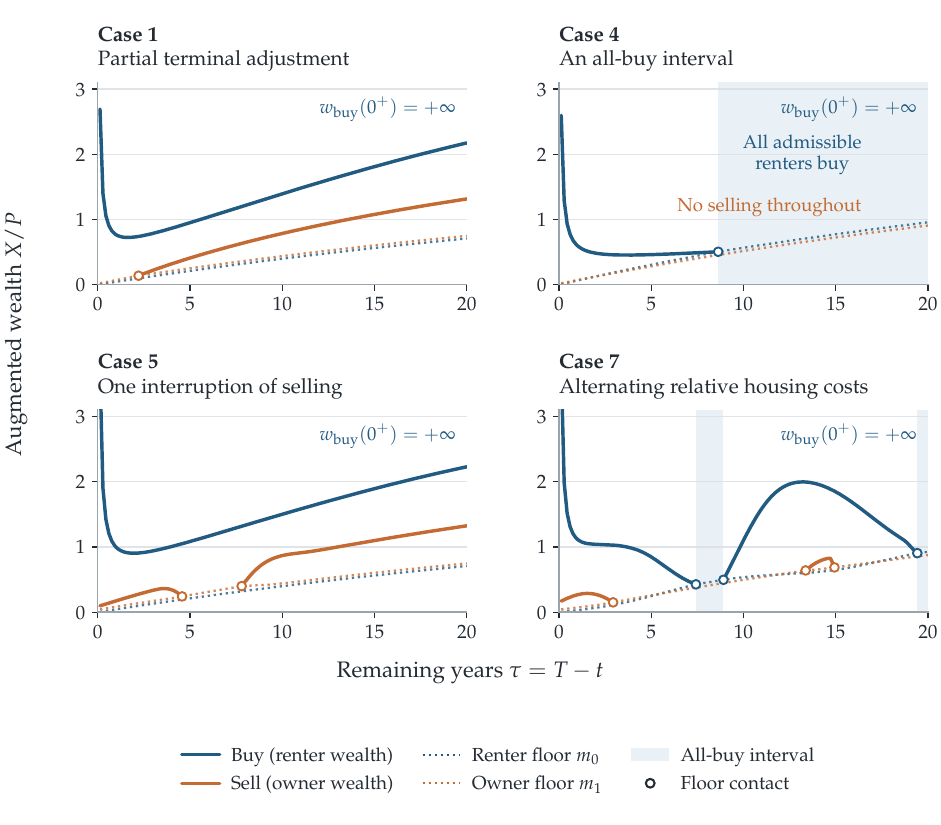}
\caption{Representative pre-switch wealth boundaries for Cases 1, 4, 5
and 7. Solid blue and orange curves are the renter's buying and owner's
selling thresholds. Dotted curves show the corresponding minimum
housing-financing costs $m_0$ and $m_1$. Shading denotes all-buy intervals.
Open circles mark limits at a financing floor, not extra finite boundaries
on an inactive slice. Missing selling curves mean no selling at any
admissible owner wealth. The vertical display ends at 3.1 property values,
and the terminal buying threshold diverges beyond it.}
\label{fig:wealth-cases}
\end{figure}

Figure~\ref{fig:wealth-cases} presents Cases 1, 4, 5 and 7 in
$X/P=(W+\mathcal H(t))/P$. These are pre-switch wealth thresholds:
the buying threshold is evaluated with the renter value, and the selling
threshold with the owner value. Deterministic labor income translates the
net-worth thresholds by exactly $-\mathcal H(t)$, and all wealth axes use
$X$. The eleven cases concern two
free boundaries and their time domains, not eleven different
boundaries. Proposition~\ref{prop:wealth-coordinates} explains why their
active time components carry over to wealth, whereas monotonicity or
curvature of the plotted curves need not carry over.

In Case 1, a renter buys only after reaching the upper wealth threshold,
whereas an owner sells at sufficiently low wealth when the sale region
exists. With ten remaining years, the two pre-switch levels are 1.3940
and 0.7852 property values. Their separation illustrates the importance
of current tenure for households facing the same prices and opportunities.
The absence of selling near the horizon is visible as a missing orange
curve, not a zero wealth threshold.
The computed buying threshold is U-shaped in remaining time, although
the corresponding transformed dual threshold $b(\tau)$ decreases over
these horizons.
A longer horizon changes both the switching incentive and the resources
needed to support the remaining consumption and housing plan. The wealth
conversion therefore matters economically: a lower trigger in $s$ does
not by itself imply a lower absolute wealth requirement.

Case 4 shows an all-buy interval produced by a favorable long-run
ownership opportunity alone, since the model has no liquidity or
down-payment restriction. Once at least 8.6283 years remain,
every renter with $X/P>m_0(t)$ buys. The finite buying boundary meets the
renter's financing floor when the all-buy interval begins. Unrestricted
financial hedging and borrowing are essential to this interpretation.
The statement does not apply to households outside the solvency domain.

Case 5 shows how the timing of a known selling-cost surcharge changes
the availability of an exit from ownership. During the interruption,
selling is unattractive even at low admissible wealth, and a finite
selling threshold reappears away from that interval at the dates of
Table~\ref{tab:components}. Empirical evidence that transfer
taxes discourage housing-related moves is provided by
\citet{HilberLyytikainen2017}. This supports the transaction-cost channel,
without identifying the complete disappearance predicted in this example.

Case 7 combines changes in the location of wealth thresholds with changes
in their existence. A predictable relative-cost cycle can create separate
opportunities to enter ownership and to exit it. The figure represents
policy regions across remaining horizons, not a realized household path.
It must be read from right to left as calendar time advances.

\subsection{Transaction costs, terminal adjustment and carrying costs}
\begin{figure}[htbp]\centering
\includegraphics[width=.86\textwidth]{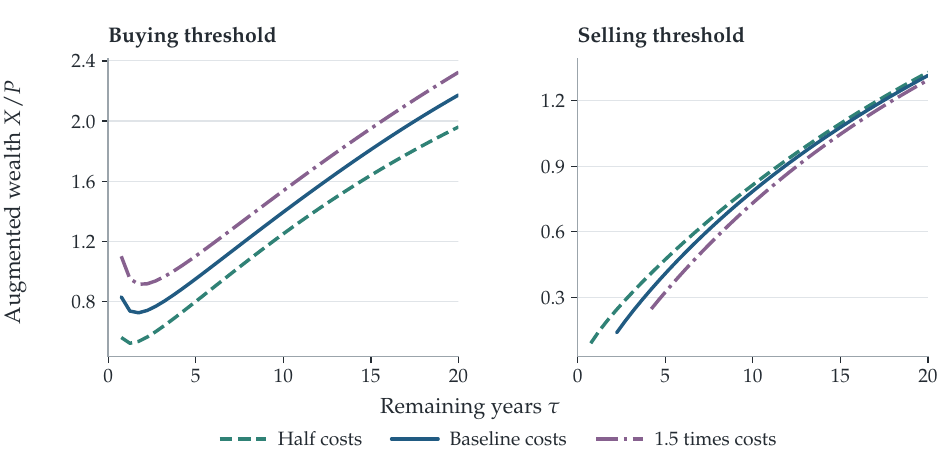}
\caption{Transaction costs and tenure persistence. Both thresholds are
pre-switch augmented wealth divided by current house price. The terminal
adjustment is $\chi_1=0.015$ in every curve.}\label{fig:costs}
\end{figure}

Figure~\ref{fig:costs} compares one half, one and one-and-a-half times
the baseline buying and selling costs, holding the terminal adjustment
fixed. In this experiment higher costs increase the pre-purchase wealth
threshold and reduce the wealth level at which an owner sells. The
resulting larger inaction interval is tenure persistence caused by the
resource cost of changing housing, even though financial risk can be
hedged continuously. This is not a mortgage or down-payment effect, since
the market allows unrestricted borrowing and hedging.

At ten remaining years, the baseline pre-purchase and pre-sale thresholds
are $X/P=1.3940$ and $0.7852$. Halving both transaction costs changes them
to $1.2495$ and $0.8142$, and raising costs by one half changes them to
$1.5360$ and $0.7317$. Thus the threshold separation increases from
$0.4352$ to $0.8043$ property values between these two experiments.
These quantities compare thresholds in augmented wealth. Actual wealth
drops by the appropriate transaction cost when a switch occurs. For a renter
with $X/P=1.30$ and ten remaining years, the lower-cost schedule changes
the current decision from waiting to buying. For an owner with $X/P=0.80$,
it changes waiting to selling. These are decisions read from the computed
thresholds at the same initial wealth and tenure, not an estimate of a
population response to a transaction-tax reform.

\begin{figure}[htbp]\centering
\includegraphics[width=.86\textwidth]{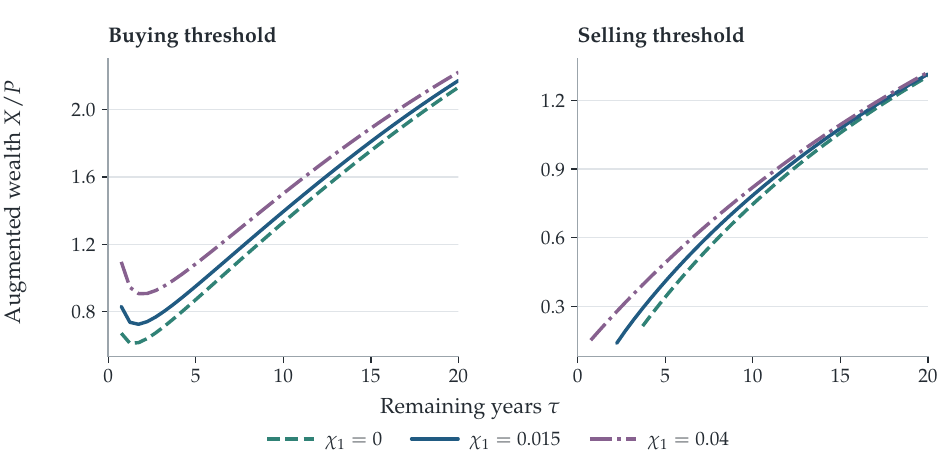}
\caption{Terminal housing adjustment. Missing portions of a selling curve
mean no finite-state sale region, as determined by
\eqref{eq:sale-onset-economic}. Buying and selling transaction costs are
held fixed.}\label{fig:terminal}
\end{figure}

Terminal accounting has a distinct effect, shown in
Figure~\ref{fig:terminal}. Write $\delta=r_Z>0$ in this constant-coefficient comparison.
For constant $d>\delta\psi_1$ and constant
costs, with no intervening lower limiting contact, the onset of selling is
\begin{equation}\label{eq:sale-onset-economic}
 \tau_{\rm sell}
 =\frac1\delta\log\frac{d/\delta-\chi_1}{d/\delta-\psi_1}.
\end{equation}
For the baseline parameters this is 3.480455 years when $\chi_1=0$,
2.210507 years when $\chi_1=0.015$, and zero under full terminal
liquidation. A small terminal adjustment makes a late discretionary sale
unattractive because the household can retain the property at the horizon.
With full liquidation, the sale cost is already embedded in the terminal
alternative, so a finite sale threshold survives arbitrarily close to $T$.
The higher purchase threshold under a larger terminal adjustment reflects
the additional ownership expense that must be recovered over the remaining
residence period.

Section~\ref{ec:numerics} of the e-companion shows how the risk-adjusted
carry moves both thresholds and reports the optimal consumption and the
owner's financial hedge.

\subsection{The value of tenure flexibility}
\begin{figure}[htbp]\centering
\includegraphics[width=.86\textwidth]{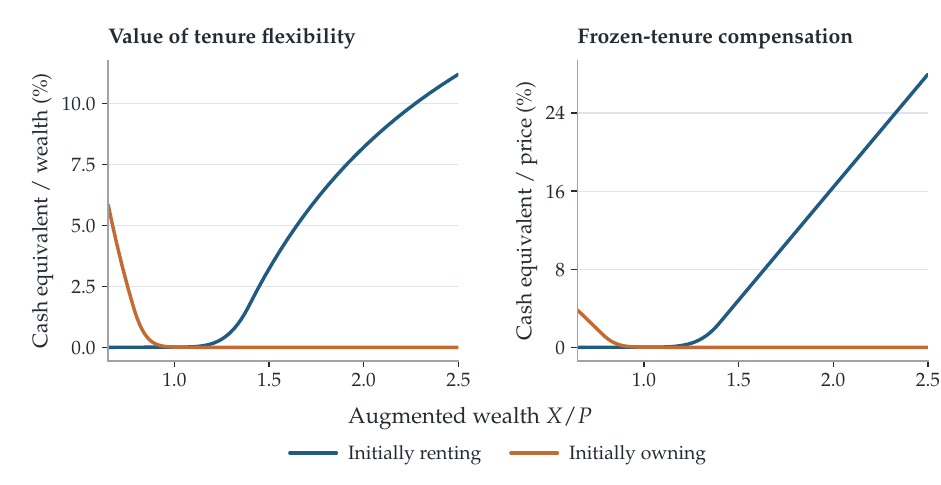}
\caption{Cash-equivalent value of tenure flexibility at ten remaining
years. The no-switch comparison is evaluated on its financeable domain.
The left panel scales compensation by wealth, and the right scales it by
current house price.}\label{fig:welfare}
\end{figure}

Figure~\ref{fig:welfare} compares optimal switching with a
household forced to retain its initial tenure for the entire remaining
horizon. Let $V_i^{\rm fixed}$ denote the corresponding no-switch value.
The cash-equivalent value of flexibility $\Delta_i(x)$ is defined by
\begin{equation}\label{eq:cash-equivalent}
 V_i^{\rm fixed}(t,x+\Delta_i(x),p)
       =\widetilde V_i(t,x,p).
\end{equation}
The fixed-tenure consumption--investment problem is solved analytically,
and the switching value uses the Gaussian integral reconstruction and the
unique multiplier. Thus the comparison holds preferences, market risks,
income treatment and initial tenure fixed. Its value reflects both future
contingent switches and the ability to change tenure immediately. The two
initial regimes have different gains because keeping an expensive housing
bundle and foregoing a superior housing bundle impose different losses.

For example, at $X/P=2.5$ an initially renting household values switching
flexibility at $0.2797$ property values, or $11.19$ percent of wealth.
At $X/P=0.65$, the corresponding compensation for an initially owning
household is $0.0380$ property values, or $5.85$ percent of wealth. At both
wealth levels the household switches immediately, and a single switch
followed by fixed tenure reproduces both numbers to the displayed
precision, so they measure the value of the immediate change. The value of
future contingent switches is the gain in the middle of the inaction
region, which is much smaller. These
asymmetries follow from the initial tenure and future adjustment option,
without adding financing constraints or changing risk preferences.

\FloatBarrier

\section{Conclusion}\label{sec:conclusion}

Time-dependent housing costs determine both the dates at which some owner
sells or every renter buys and the wealth required to change tenure. In a complete
market with unrestricted portfolios, the dates are those of a deterministic
minimum cost schedule for housing, discounted at the risk-adjusted carry,
and they do not depend on preferences, housing services, or labor income.
Because financial hedging offsets the house-price exposure of the home,
house-price coefficients enter these dates only through the carry. When owning costs at least as much as renting, the timing windows
concern selling alone. Running records describe when the windows open and close, and strong
duality connects the boundaries inside the windows to consumption and
financial investment on the natural financing domain. The windows created by
an announced cost schedule can therefore be computed before any preference
parameter is chosen. The windows are a property of the financing floors,
and two experiments in the e-companion show what this means outside the
benchmark. When the ownership cost is subject to an uninsurable permanent
increase, the dates at which some owner sells remain deterministic in our
computations. They are the dates at which the worst-case minimum cost
schedule sells immediately and the sale incentive at the current cost is
positive, and they did not move under any of the preference changes we
computed. When short sales or
borrowing are restricted, an owner with less than one house value of net
worth cannot remove its house-price exposure and sells to stay solvent, so
a selling region exists at every date and the windows disappear. Whether a
deterministic description survives unspanned house-price risk is open.

\clearpage
\setcounter{section}{0}\setcounter{figure}{0}\setcounter{table}{0}\setcounter{remark}{0}
\renewcommand{\thesection}{EC.\arabic{section}}
\renewcommand{\thefigure}{EC.\arabic{figure}}
\renewcommand{\thetable}{EC.\arabic{table}}

\begin{center}
{\Large\bfseries E-Companion\par}
\vspace{4pt}
{\large Renting, Owning, and Moving: Housing Tenure Choice with Transaction Costs\par}
\end{center}
\vspace{6pt}

\noindent
This e-companion contains auxiliary results (Section~\ref{ec:auxiliary}), the proofs of all results
(Section~\ref{sec:proofs}), the integral representation and numerical
implementation (Section~\ref{sec:integral-method}), additional contact
configurations (Section~\ref{ec:cases}), additional numerical results
(Section~\ref{ec:numerics}), and numerical experiments with portfolio
restrictions and an uninsurable cost increase (Section~\ref{ec:robust}).
Numbers without the prefix EC refer to the main text above.

\section{Auxiliary results}\label{ec:auxiliary}

This section states the budget constraint and the duality criterion of
Section~\ref{sec:model}, the comparison and existence results for
\eqref{eq:vdob}, the running-record form of the timing windows, and the
local regularity and endpoint results summarized in
Section~\ref{sec:regularity}. Their proofs are in Section~\ref{sec:proofs}.

\subsection{Budget constraint and duality}
\begin{lemma}\label{lem:budget}
For every admissible strategy $(c,\pi,\vartheta)$ starting from $(t,w,p,i)$,
\begin{align}\label{eq:budget-W}
    \E_t\bigg[
        &\int_t^T \Gamma_{t,s}
        \Big(c_s+P_s\wh a_{I_s}(s)-\ell(s)\Big)\,\dd s
        +\Gamma_{t,T}W_T 
        +\sum_{n\ge1:\,\tau_n<T}
        \Gamma_{t,\tau_n}P_{\tau_n}
        k_{\iota_{n-1},\iota_n}(\tau_n)
    \bigg]\le w.
\end{align}
Equivalently, after adding human wealth $x=w+\mathcal H(t)$,
\begin{align}\label{eq:budget-X}
    \E_t\bigg[
        &\int_t^T \Gamma_{t,s}
        \Big(c_s+P_s\wh a_{I_s}(s)\Big)\,\dd s
        +\Gamma_{t,T}W_T
        +\sum_{n\ge1:\,\tau_n<T}
        \Gamma_{t,\tau_n}P_{\tau_n}
        k_{\iota_{n-1},\iota_n}(\tau_n)
    \bigg]\le x.
\end{align}
If the deflated gain process is a uniformly integrable martingale, then the
corresponding inequality holds with equality.
\end{lemma}

The standard martingale verification argument of \citet{JeonYang} gives
the following criterion. Its attainment and intermediate financeability
conditions will be established on the natural wealth domain in
Theorem~\ref{thm:complete-verification}.

\begin{proposition}
\label{prop:duality}
For every initial state with nonempty admissible set, weak duality
\eqref{eq:duality-bound} holds. Suppose there
exist \(y^*>0\) and a dual-optimal switching strategy \(\vartheta^*\) such that, with
\[
 c_s^*=c_{I_s^*}^*(Y_s^{t,y^*}),
 \qquad
 W_T^*=P_T\chi_{I_T^*}
       +\eta^{1/\gamma_1}(Y_T^{t,y^*})^{-1/\gamma_1},
\]
the complete financial market admits a portfolio \(\pi^*\) financing
\(c^*\), all switching costs generated by \(\vartheta^*\), and \(W_T^*\) from initial
augmented wealth \(x\), with the budget equality attained and all intermediate
financeability conditions satisfied.  Then \((c^*,\pi^*,\vartheta^*)\) is primal
optimal and
\begin{equation}\label{eq:strong-duality}
    \wt V_i(t,x,p)=
    \cJ_i(t,y^*,p)+y^*x
    =\inf_{y>0}\{\cJ_i(t,y,p)+yx\}.
\end{equation}
A minimizing multiplier is characterized in general by the convex-subgradient
condition
\begin{equation}\label{eq:foc}
    x\in-\partial_y\cJ_i(t,y^*,p),
\end{equation}
which reduces to \(x=-\partial_y\cJ_i(t,y^*,p)\) at differentiability points.
\end{proposition}

\subsection{Comparison and existence}
\begin{lemma}\label{lem:bounded-comparison}
Let $u,\bar u$ be bounded continuous strong solutions, locally in
$W^{1,2}_p$ for some $p>3$, with the same time-dependent obstacles and
operator as in \eqref{eq:vdob}, sources $f\le\bar f$, and ordered initial
values $u(0,\cdot)\le\bar u(0,\cdot)$. The sources need only be locally
bounded, and no common polynomial growth bound is required. Then $u\le\bar u$.
The same comparison holds on a bounded cylinder with ordered parabolic
boundary data, for sub- and supersolutions satisfying the appropriate
inequality on their positive difference set, and for a common penalized
equation with a penalty nonincreasing in the unknown. For the scalar
reflected ODE, comparison holds in the class $AC([0,T])$.
\end{lemma}

The proof below uses the bounded difference of the two solutions, so the
exponential growth of $\mathcal R$ as $\xi\to-\infty$ causes no difficulty.
For comparison results for more general systems on unbounded sets, see
\citet{YangYan2008}.

\begin{lemma}\label{lem:dob-existence}
There exists a unique bounded strong solution
$v\in W^{1,2}_{p,\mathrm{loc}}(\Omega_T)\cap C([0,T]\times\mathbb R)$
for every $1<p<\infty$ to \eqref{eq:vdob}. More precisely,
$v\in W^{1,2}_p((0,T)\times(-m,m))$ for every $m<\infty$ and
$1<p<\infty$, including estimates up to the initial time.
\end{lemma}

\subsection{Finite upper contact}
\begin{lemma}\label{lem:finite-contact-interior}
Let $(\tau_0,\xi_0)\in\Omega_T$ be such that
$v(\tau_0,\xi_0)=\psi_1(\tau_0)$. Then
\[
    \tau_0\in\operatorname{int}\mathcal C_1(\infty).
\]
\end{lemma}

We next prove the converse implication: an interior contact time of
$v_\infty$ satisfying the pointwise necessary condition indeed produces
contact at finite spatial points.

\begin{lemma}\label{lem:recover-upper-contact}
Suppose that
$\tau_0\in\operatorname{int}\mathcal C_1(\infty)$ and
$\Lambda_1(\tau_0)<\infty$. Then there exist $\nu>0$ and $\Xi_{\tau_0}>0$ such
that
\[
    v(\tau,\xi)=\psi_1(\tau)
    \qquad
    \text{for }|\tau-\tau_0|<\nu,\quad \xi\ge \Xi_{\tau_0}.
\]
\end{lemma}

\subsection{Running-record form of the timing windows}

The lower running-record criterion yields the following equivalent form of
Theorem~\ref{thm:timing-windows}(i).
\begin{corollary}\label{cor:lower-time-characterization}
For every $\tau\in(0,T)$,
    $\tau\notin\mathcal T_0$,
or equivalently
\[
    v(\tau,\xi)=-\psi_0(\tau)
    \qquad\text{for all }\xi\in\R,
\]
if and only if
\[
\begin{aligned}
    &\begin{cases}
        F_0(\tau)+F_1(\ell_1(\tau))\ge0,
        &\ell_1(\tau)>0,\\
        F_0(\tau)\ge0,
        &\ell_1(\tau)=0,
    \end{cases}\\
    &F_0(\tau)
    =
    \max_{\ell_1(\tau)\le s\le\tau}F_0(s),
    \qquad
    F_0'(\tau)\ge0.
\end{aligned}
\]
\end{corollary}

The finite selling boundary requires strict level and record conditions.

\begin{lemma}\label{lem:upper-time-characterization}
For every $\tau\in(0,T)$, $\tau\in\mathcal T_1$ if and only if
the following conditions hold:
\[
\begin{aligned}
    &\text{\rm (i)}\quad
    \begin{cases}
        F_1(\tau)+F_0(\ell_0(\tau))>0,
        &\ell_0(\tau)>0,\\
        F_1(\tau)>0,
        &\ell_0(\tau)=0,
    \end{cases}\\
    &\text{\rm (ii)}\quad
    F_1(\tau)>F_1(s)
    \quad\text{for every }s\in[\ell_0(\tau),\tau),\\
    &\text{\rm (iii)}\quad
    F_1'(\tau)>0.
\end{aligned}
\]
\end{lemma}

\begin{remark}
At every positive remaining time, sufficiently small dual states lie in the
buying region. Outside $\mathcal T_0$, buying is optimal for every dual state.
Selling is possible at a finite state exactly on $\mathcal T_1$. The asymmetry
between the weak lower and strict upper record conditions reflects ownership's
service advantage: limiting upper contact alone does not make selling optimal.
\end{remark}

\subsection{Local estimates near the free boundaries}
We now work locally on the open sets $\mathcal T_0$ and $\mathcal T_1$.
Semicontinuity and derivative estimates give boundary non-degeneracy.

\begin{proposition}\label{prop:semicontinuity-free-boundaries}
The lower free boundary $\Xi_0$ is upper semicontinuous on
$\mathcal T_0$, and the upper free boundary $\Xi_1$ is lower
semicontinuous on $\mathcal T_1$.
\end{proposition}

\begin{lemma}\label{lem:derivative-estimate}
The following local derivative estimates hold.

\begin{enumerate}
    \item[(i)] For each $\tau_0\in\mathcal T_0$, there exist
    $\delta=\delta(\tau_0)>0$, $\eta=\eta(\tau_0)>0$, and
    $C=C(\tau_0)>0$ such that
    \[
        \lvert\partial_\tau v(\tau,\xi)+\psi_0'(\tau)\rvert
        \le C\,\partial_\xi  v(\tau,\xi)
    \]
    whenever
    \[
        \tau\in(\tau_0-\delta,\tau_0+\delta)\subset\mathcal T_0,
        \qquad
        \Xi_0(\tau)<\xi\le \Xi_0(\tau_0)+\eta.
    \]

    \item[(ii)] For each $\tau_0\in\mathcal T_1$, there exist
    $\delta=\delta(\tau_0)>0$, $\eta=\eta(\tau_0)>0$, and
    $C=C(\tau_0)>0$ such that
    \[
        \lvert\partial_\tau v(\tau,\xi)-\psi_1'(\tau)\rvert
        \le C\,\partial_\xi  v(\tau,\xi)
    \]
    whenever
    \[
        \tau\in(\tau_0-\delta,\tau_0+\delta)\subset\mathcal T_1,
        \qquad
        \Xi_1(\tau_0)-\eta\le \xi<\Xi_1(\tau).
    \]
\end{enumerate}
\end{lemma}

Using the preceding local estimates, we obtain the following
non-degeneracy estimates.

\begin{lemma}\label{lem:non-degeneracy}
For each $\tau_0\in\mathcal T_0$, there exist constants
$c_0>0$, $\eta>0$, and $\delta>0$ such that, whenever
$|\tau-\tau_0|<\delta$ and
$\Xi_0(\tau)<\xi<\Xi_0(\tau)+\eta$,
\[
    \partial_\xi v(\tau,\xi)
    \ge c_0\bigl(\xi-\Xi_0(\tau)\bigr),
    \qquad
    v(\tau,\xi)
    \ge -\psi_0(\tau)
    +c_0\bigl(\xi-\Xi_0(\tau)\bigr)^2.
\]

Similarly, for each $\tau_0\in\mathcal T_1$, there exist constants
$c_1>0$, $\eta>0$, and $\delta>0$ such that, whenever
$|\tau-\tau_0|<\delta$ and
$\Xi_1(\tau)-\eta<\xi<\Xi_1(\tau)$,
\[
    \partial_\xi v(\tau,\xi)
    \ge c_1\bigl(\Xi_1(\tau)-\xi\bigr),
    \qquad
    v(\tau,\xi)
    \le \psi_1(\tau)
    -c_1\bigl(\Xi_1(\tau)-\xi\bigr)^2.
\]
\end{lemma}

The derivative estimates control the slopes of nearby level curves, yielding
local Lipschitz continuity as in \citet{HaJeonOk2026}.

\begin{lemma}\label{lem:local-lipschitz-boundaries}
The free boundaries $\Xi_0$ and $\Xi_1$ are locally Lipschitz continuous
on $\mathcal T_0$ and $\mathcal T_1$, respectively.
\end{lemma}

In one spatial dimension, a graph Lipschitz in time is locally $C_p^1$:
its $1/2$-H\"older time modulus vanishes at small scales. The boundary Harnack
inequality of \citet[Theorem~1.2]{Kukuljan2022}, applied to the differentiated
equations, therefore improves regularity.

\begin{lemma}\label{lem:boundary-C1alpha}
For every $\alpha\in(0,1/2)$,
\[
    \Xi_0\in C_{\rm loc}^{1,\alpha}(\mathcal T_0),
    \qquad
    \Xi_1\in C_{\rm loc}^{1,\alpha}(\mathcal T_1).
\]
\end{lemma}

\begin{remark}\label{rem:higher-order}
With correspondingly smoother coefficients and obstacles, the higher-order
boundary Harnack inequality \cite[Theorem~1.3]{Kukuljan2022} and the level-curve
identities bootstrap this regularity, and smooth data give locally smooth boundaries.
\end{remark}

\subsection{Endpoint behavior}
The following endpoint limits exclude additional finite spatial boundary
pieces at the edges of the active time domains.

\begin{lemma}\label{lem:boundary-behavior-free-boundaries}
The lower free boundary satisfies
\[
    \lim_{\tau\to0^+}\Xi_0(\tau)=-\infty,
    \qquad
    \lim_{\tau\to\tau_*,\,\tau\in\mathcal T_0}
    \Xi_0(\tau)=+\infty
\]
for every $\tau_*\in\partial\mathcal T_0\cap(0,T)$. Moreover,
\[
    \lim_{\tau\to\tau^*,\,\tau\in\mathcal T_1}
    \Xi_1(\tau)=+\infty
\]
for every $\tau^*\in\partial\mathcal T_1\cap(0,T)$.
\end{lemma}

\begin{corollary}\label{cor:free-boundary-graphs}
The graphs of $\Xi_0$ and $\Xi_1$ constitute the entire free boundaries
relative to the open space-time domain $\Omega_T$. More precisely,
\[
    \partial\{v>-\psi_0\}
    =
    \{(\tau,\Xi_0(\tau)):\tau\in\mathcal T_0\},
    \qquad
    \partial\{v<\psi_1\}
    =
    \{(\tau,\Xi_1(\tau)):\tau\in\mathcal T_1\}.
\]
\end{corollary}

We finally consider the upper free boundary at the initial time when
the initial datum coincides with the upper obstacle.

\begin{lemma}\label{lem:initial-upper-boundary}
Suppose that
\[
    \psi_1(0)=\chi_1
    \qquad\text{and}\qquad
    \Lambda_1(0)<\infty.
\]
Then there exists $\delta>0$ such that
$(0,\delta)\subset\mathcal T_1$, and
\[
    \lim_{\tau\to0^+}\Xi_1(\tau)
    =
    \Lambda_1(0)
    =
    \gamma_1\log\left(\frac{C_H}{G_1(0)}\right).
\]
\end{lemma}

\begin{lemma}\label{lem:remaining-terminal-cases}
If $\chi_1<\psi_1(0)$, there is a $\delta>0$ such that
$\mathcal T_1\cap(0,\delta)=\varnothing$.
The same conclusion holds
if $\chi_1=\psi_1(0)$ and $G_1(0)<0$.
Thus, in either case, $0$ is not
a limit point of $\mathcal T_1$.

If $\chi_1=\psi_1(0)$ and $G_1(0)=0$, then either
$\mathcal T_1\cap(0,\delta)=\varnothing$ for some $\delta>0$,
or
\[
    \lim_{\tau\to0^+,\,\tau\in\mathcal T_1}
    \Xi_1(\tau)=+\infty.
\]
\end{lemma}

\section{Proofs}\label{sec:proofs}
The proofs concentrate on the contact classification and boundary analysis.
For standard martingale verification and local PDE constructions we use
\citet{JeonYang,HaJeonOk2026}, specifying the adaptations needed for housing
costs, different terminal data, and possibly disconnected active time domains.

\subsection{Reduction identities and standard verification}
Proposition~\ref{prop:neutralization} follows by substituting
$\pi=\widehat\pi-q_iP\zeta_P$ and $\mu-r\one=\Sigma\theta$ into
\eqref{eq:W-dyn}, and the remaining owner drift is $q_iP\kappa_H$.
For Lemma~\ref{lem:no-terminal-switch}, selling at $T$ gives adjusted wealth
$W_{T-}-P_Tk_{10}(T)\le W_{T-}-P_T\chi_1$, while buying gives
$W_{T-}-P_T(k_{01}(T)+\chi_1)<W_{T-}$.

Lemma~\ref{lem:budget} and Proposition~\ref{prop:duality} use the standard
deflator and Fenchel arguments of
\citet{KaratzasLehoczkyShreve1987,CoxHuang1989,JeonYang}.
Here the only budget change is the effective cost $\widehat a_i$:
between switches,
$d(\Gamma_{t,s}W_s)=\Gamma_{t,s}(-c_s-P_s\widehat a_{I_s}+\ell)\,ds+dM_s$,
and each switch contributes its deflated transaction cost. The assumed
supermartingale property gives the budget inequality, with equality for a
uniformly integrable martingale. Fenchel equality at the displayed optimizers
then gives strong duality whenever the stated attainment conditions hold.
Their model-specific financeability checks are supplied below.

Finally, Proposition~\ref{prop:smooth-fit-complete} follows directly from
the spatial $C^1$ regularity in Lemma~\ref{lem:dob-existence}: the nonnegative
gaps from the smooth obstacles have zero derivative at contact.

\subsection{Proof of Lemma~\ref{lem:bounded-comparison}}
\begin{proof}
Put $P=\partial_\tau-\mathcal L_\xi $. On $\{u>\bar u\}$, the common
obstacles imply $u>-\psi_0$, $\bar u<\psi_1$, and hence
$P(u-\bar u)\le f-\bar f\le0$ a.e. Choose $\lambda$ with
$\widetilde r+\lambda\ge1$. Since $u-\bar u$ is bounded, the positive set of
\[
 e^{-\lambda\tau}(u-\bar u)-\varepsilon e^{K\tau}(1+\xi^2)
\]
is spatially bounded, and bounded coefficients allow $K$ large enough that the
subtracted function is a strict supersolution. The strong-solution maximum
principle \citep[Chapters~2 and~7]{Lieberman} and
$\varepsilon\downarrow0$ give comparison. Thus no growth bound on the
individual sources is needed (compare \citet{YangYan2008}).
The same argument on bounded cylinders uses the prescribed boundary order.
A nonincreasing penalty adds a nonpositive term on the positive difference
set. For the reflected ODE, the positive-part chain rule gives
$(e^{2\int_0^\tau\widetilde r}(u-\bar u)_+^2)'\le0$ a.e., proving
comparison and uniqueness.
\end{proof}

\subsection{Proof of Lemma~\ref{lem:dob-existence}}
\begin{proof}
We use the bounded-cylinder penalization and exhaustion construction of
\citet{JeonYang,HaJeonOk2026,Ha_Jeon_OK}, recording the changes needed for
the nonzero initial datum and the unbounded source. On
$\Omega_T^n=(0,T)\times(-n,n)$ choose
\[
 M_n>\|\mathcal R\|_\infty+
       \max_i\|\psi_i'+\widetilde r\psi_i\|_\infty,
 \qquad g(\tau)=\chi_1e^{-K\tau},\quad
 K>\|\psi_1'/\psi_1\|_\infty.
\]
Then $-\psi_0<g\le\psi_1$. With a smooth nonincreasing cutoff
$\varphi=1$ on $(-\infty,0]$ and $\varphi=0$ on $[1,\infty)$, solve
\[
 Pv_{n,\varepsilon}=\mathcal R+
 M_n\varphi((v_{n,\varepsilon}+\psi_0)/\varepsilon)
 -M_n\varphi((\psi_1-v_{n,\varepsilon})/\varepsilon),
\]
with lateral data $g$ and initial value $\chi_1$, taking
$2\varepsilon<\min(\psi_0+\psi_1)$. Here $P=\partial_\tau-\mathcal L_\xi $.
The standard linear estimates and fixed-point construction
\citep[Theorem~7.32]{Lieberman}, \citep[Corollary~11.2]{GilbargTrudinger},
together with the ABP estimate \citep{Tso-1985}, give a solution.
The obstacles are sub- and supersolutions, so comparison bounds the solution
between them and bounds both penalties independently of $\varepsilon$.
The $W^{1,2}_p$ estimates and compact embedding therefore yield a strong
obstacle solution $v_n$ as $\varepsilon\downarrow0$.

For exhaustion, the level-set property of Sobolev derivatives gives
\[
 Pv_n=(-\psi_0'-\widetilde r\psi_0)\one_{\{v_n=-\psi_0\}}
 +(\psi_1'+\widetilde r\psi_1)\one_{\{v_n=\psi_1\}}
 +\mathcal R\one_{\{-\psi_0<v_n<\psi_1\}}.
\]
On any fixed spatial strip this is bounded independently of $n$.
Interior spatial estimates including the constant initial datum
\citep[p.~355]{Ladyzhenskaya} give uniform $W^{1,2}_p$ bounds on smaller
strips up to time zero. A diagonal limit in $n$ and $p$ yields $v$ and
its continuous spatial derivative. Testing each inequality on compact subsets
where the opposite obstacle is inactive passes it to the limit. The obstacle
bounds give global boundedness, and Lemma~\ref{lem:bounded-comparison}
gives uniqueness.
\end{proof}

\subsection{Proof of Lemma~\ref{lem:state-monotonicity}}
\begin{proof}
For $h>0$, the translate $v(\tau,\xi+h)$ has the same obstacles and initial
data and a larger source. Comparison gives $v(\tau,\xi+h)\ge v(\tau,\xi)$.
The spatial derivative is continuous by Sobolev embedding
\citep[Lemma~3.3, p.~80]{Ladyzhenskaya}, so $v_\xi\ge0$ everywhere.
Spatial difference quotients in the continuation set give
\[
 (\partial_\tau-\mathcal L_\xi )v_\xi
 =\mathcal R_\xi =\frac{C_H}{\gamma_1}e^{-\xi/\gamma_1}>0.
\]
The strong minimum principle \citep[Chapters~2 and~7]{Lieberman}, with
a time exponential if needed for the zero-order term, excludes an interior zero.
\end{proof}

\subsection{Proof of Lemma~\ref{lem:contact-location}}
\begin{proof}
Write $P=\partial_\tau-\mathcal L_\xi $ and $u_0=v+\psi_0$.
Near a lower contact, the opposite obstacle is inactive and
$Pu_0\ge G_0-C_He^{-\xi/\gamma_1}$. If $G_0(\tau_0)>0$ and
$\xi_0>\Lambda_0(\tau_0)$, the right-hand side is strictly positive
in a neighborhood, contradicting the strong minimum principle at the zero
of $u_0$. Thus $\xi_0\le\Lambda_0(\tau_0)$.
Equality is impossible: locally $\Lambda_0$ is Lipschitz, and on its right
the preceding non-strict inequality excludes lower contact. There
$u_0>0$ and $Pu_0=G_0-C_He^{-\xi/\gamma_1}>0$.
Flattening $\xi=\Lambda_0(\tau)$ gives a uniformly parabolic operator with
bounded drift. The Hopf lemma \citep[Chapters~2 and~7]{Lieberman}
would give $u_{0,\xi}>0$ at equality, contrary to spatial smooth fit.
When $G_0\le0$, the asserted bound is automatic.

At upper contact use $u_1=\psi_1-v\ge0$, for which
$Pu_1\ge C_He^{-\xi/\gamma_1}-G_1$. The strong minimum principle excludes
$G_1\le0$ and $\xi<\Lambda_1$. At equality, the continuation region to
the left of the locally Lipschitz graph $\Lambda_1$ has $Pu_1>0$,
and the Hopf lemma gives $u_{1,\xi}<0$, again contradicting smooth fit.
\end{proof}

\subsection{Proof of Lemma~\ref{lem:limit-v}}
\begin{proof}
The limit defining $v_\infty$ is well defined by the spatial
monotonicity of $v$ and the obstacle bounds. In particular,
\[
    v(\tau,\xi)\le v_\infty(\tau)
    \qquad
    \text{for all }(\tau,\xi)\in\Omega_T.
\]

For each $n\in\mathbb N$, define
\[
    v^{(n)}(\tau,\xi):=v(\tau,\xi+n).
\]
Then $v^{(n)}$ satisfies the same double-obstacle problem as $v$,
with source
\[
    \mathcal R^{(n)}(\tau,\xi)
    :=
    -C_H e^{-(\xi+n)/\gamma_1}-\mathfrak a(\tau).
\]
Clearly, $\mathcal R^{(n)}\to-\mathfrak a$ locally uniformly on
$(0,T)\times\R$. Moreover, for every fixed $(\tau,\xi)$,
\[
    v^{(n)}(\tau,\xi)
    =
    v(\tau,\xi+n)
    \longrightarrow
    v_\infty(\tau)
    \qquad\text{as }n\to\infty.
\]

The strip estimates in Lemma~\ref{lem:dob-existence} apply to the
translates with constants uniform in $n$: on every fixed spatial strip
their sources are uniformly bounded, and their initial data equal $\chi_1$.
Compact embedding therefore gives convergence along subsequences
uniformly on $[0,T]\times[-m,m]$, and weakly in $W^{1,2}_p$.
The pointwise limit identifies every such subsequence, so the whole
sequence converges uniformly. Spatial monotonicity then extends the
conclusion from integer translations to all real $\xi\to\infty$.
Passing to the limit in the variational inequalities is justified by
the same compact-support test argument as in the existence proof. Since the limiting function is independent of $\xi$,
the spatial derivative terms vanish, while
$\mathcal R^{(n)}\to-\mathfrak a$. Hence $v_\infty$ satisfies
\eqref{v,infty}.
The initial condition $v_\infty(0)=\chi_1$ follows from
$v(0,\xi)=\chi_1$, and uniqueness follows from the absolutely continuous ODE comparison in
Lemma~\ref{lem:bounded-comparison}. On each contact set the derivative
equals the derivative of that obstacle a.e., and on the continuation set
the ODE gives $v_\infty'=-\widetilde r v_\infty-\mathfrak a$.
These quantities are bounded, proving $v_\infty\in W^{1,\infty}(0,T)$.

If $\tau\in\mathcal C_0(\infty)$, then
\[
    -\psi_0(\tau)
    \le v(\tau,\xi)
    \le v_\infty(\tau)
    =
    -\psi_0(\tau)
\]
for every $\xi\in\R$, and hence
$v(\tau,\xi)=-\psi_0(\tau)$ for all $\xi\in\R$. If
$v(\tau,\xi)=\psi_1(\tau)$ for some $\xi\in\R$, then
\[
    \psi_1(\tau)
    =
    v(\tau,\xi)
    \le v_\infty(\tau)
    \le \psi_1(\tau),
\]
and therefore $\tau\in\mathcal C_1(\infty)$.
\end{proof}

\subsection{Proof of Lemma~\ref{lem:running-record}}
\begin{proof}
Set $M(s)=\exp(\int_0^s\widetilde r(q)\,dq)$ and fix $\tau$.
For the upper contact criterion write $\ell=\ell_0(\tau)$ and suppose
first that $\ell<\tau$. On $(\ell,\tau)$ there is no lower contact, so
\[
 R(s):=M(s)(\vns(s)-v_\infty(s))\ge F_1(s),\qquad
 R'\ge0,\qquad R'=0\ \hbox{a.e. on }\{R>F_1\}.
\]
Its initial value is $c=-F_0(\ell)$ if $\ell>0$, and $c=0$ otherwise.
Consequently,
\begin{equation}\label{eq:running-regulator}
 R(s)=\max\left\{c,\max_{\ell\le q\le s}F_1(q)\right\},
 \qquad \ell\le s\le\tau.
\end{equation}
For completeness, the right-hand side is absolutely continuous, nondecreasing,
dominates $F_1$, and increases only where it equals $F_1$. It has initial
value $c$ because $c\ge F_1(\ell)$. These properties characterize the
one-sided reflected equation uniquely by the positive-part comparison
argument in Lemma~\ref{lem:bounded-comparison}, proving the identity.
Upper contact at $\tau$ is equivalent to $R(\tau)=F_1(\tau)$, which is
exactly the stated level and running-record conditions. If $\ell=\tau$,
continuity gives lower contact, and the proposed level condition is impossible
because $F_1(\tau)+F_0(\tau)=-M(\tau)(\psi_0+\psi_1)(\tau)<0$.

For the lower criterion use $\ell=\ell_1(\tau)$ and
$R(s)=M(s)(v_\infty(s)-\vns(s))$. The same formula holds with
$F_0$ in place of $F_1$ and initial value $-F_1(\ell)$ or zero.
Finally, the limiting obstacle inequalities give
$F_1'=MG_1\ge0$ a.e. on $\mathcal C_1(\infty)$ and
$F_0'=-MG_0\ge0$ a.e. on $\mathcal C_0(\infty)$.
\end{proof}

\subsection{Proof of Lemma~\ref{lem:lower-contact-existence}}
\begin{proof}
Set $M_0=\max(\psi_0+\psi_1)$ and
$A_\varepsilon=2(M_0+1)/\varepsilon$. On
$Q=(\varepsilon/2,T)\times(\bar\xi-1,\bar\xi+1)$ define
\[
 \zeta(\tau,\xi)=-\psi_0(\tau)
       +A_\varepsilon(\varepsilon-\tau)_+
       +(M_0+1)(\xi-\bar\xi)^2.
\]
The added terms are at least $M_0+1$ on the parabolic boundary, so
$\zeta\ge v$ there. Bounded coefficients give, with $C$ independent of $\bar\xi$,
\[
 (\partial_\tau-\mathcal L_\xi )\zeta-\mathcal R
 \ge C_He^{-(\bar\xi+1)/\gamma_1}-\|G_0\|_\infty
           -A_\varepsilon-C(M_0+1).
\]
Choose $\bar\xi=\Xi_\varepsilon$ sufficiently negative to make this nonnegative.
On $\{v>\zeta\}$, one has $v>-\psi_0$ and hence
$(\partial_\tau-\mathcal L_\xi )v\le\mathcal R$. The comparison principle
\citep[Theorem 1.1]{YangYan2008} gives $v\le\zeta$ on that set, which is a
contradiction. At $\xi=\Xi_\varepsilon$ and $\tau\ge\varepsilon$,
$\zeta=-\psi_0$, so the equality $v=-\psi_0$ holds there. Spatial monotonicity extends
contact to all $\xi\le \Xi_\varepsilon$, and continuity includes $\tau=T$.
\end{proof}

\subsection{Proof of Proposition~\ref{prop:lower-contact-times}}
\begin{proof}
By Lemma~\ref{lem:lower-contact-existence}, $\Xi_0(\tau)>-\infty$ for
every $\tau\in(0,T)$. If
$v_\infty(\tau)>-\psi_0(\tau)$, then the convergence
$v(\tau,\xi)\to v_\infty(\tau)$ as $\xi\to\infty$ shows that the lower
contact set is bounded from above, and hence $\Xi_0(\tau)\in\R$.
Conversely, if $v_\infty(\tau)=-\psi_0(\tau)$, then
$-\psi_0(\tau)\le v(\tau,\xi)\le v_\infty(\tau)$ gives
$v(\tau,\cdot)\equiv-\psi_0(\tau)$, so that $\Xi_0(\tau)=+\infty$.
\end{proof}

\subsection{Proof of Corollary~\ref{cor:lower-time-characterization}}
\begin{proof}
Proposition~\ref{prop:lower-contact-times} and
Lemma~\ref{lem:running-record} give the equivalence with the two
level and running-record conditions. To justify the additional derivative
condition pointwise, let $\tau\in\mathcal C_0(\infty)$.
Lemma~\ref{lem:limit-v} gives $v(\tau,\xi)=-\psi_0(\tau)$ for every
$\xi\in\mathbb R$. Lemma~\ref{lem:contact-location} therefore forces
$\Lambda_0(\tau)=+\infty$, equivalently $G_0(\tau)\le0$. Hence
\[
 F_0'(\tau)=-e^{\int_0^\tau\widetilde r(q)\,dq}G_0(\tau)\ge0.
\]
This argument includes isolated contacts and does not infer a pointwise
claim from an a.e. statement on a possibly measure-zero set.
\end{proof}

\subsection{Proof of Lemma~\ref{lem:finite-contact-interior}}
\begin{proof}
Finite upper contact at $(\tau_0,\xi_0)$ and spatial monotonicity give
$v(\tau_0,\xi)=v_\infty(\tau_0)=\psi_1(\tau_0)$ for all $\xi\ge \xi_0$.
Lemma~\ref{lem:contact-location} gives $G_1(\tau_0)>0$. On a small
interval about $\tau_0$, continuity ensures $v_\infty>-\psi_0$ and
$G_1>0$. With $M(\tau)=e^{\int_0^\tau\widetilde r}$, the nonnegative
gap $M(\psi_1-v_\infty)$ has derivative $-MG_1<0$ off upper contact
and zero a.e. on contact. It is nonincreasing, and hence vanishes to the
right of $\tau_0$.

If $\tau_0$ were not an interior contact time, the same monotonicity would
force $v_\infty<\psi_1$ throughout a left neighborhood. Choose a compact
spatial interval strictly to the right of $\xi_0$ and a sufficiently short
left time cylinder $Q$. Continuity gives $v>-\psi_0$ there, while
$v\le v_\infty<\psi_1$. Consequently
\[
 (\partial_\tau-\mathcal L_\xi )(v-v_\infty)
       =-C_He^{-\xi/\gamma_1}<0\quad\hbox{in }Q.
\]
The nonpositive difference vanishes on the top interior of $Q$, contradicting
the strong maximum principle. Thus limiting upper contact also persists to
the left of $\tau_0$.
\end{proof}

\subsection{Proof of Lemma~\ref{lem:recover-upper-contact}}
\begin{proof}
Choose $\kappa,\lambda>0$ so that
$[\tau_0-\kappa,\tau_0+\kappa]\subset\mathcal C_1(\infty)$ and
$G_1\ge2\lambda$ there. Write $s_0=\tau_0-\kappa$ and fix
$\tau\in[\tau_0-\kappa/2,\tau_0+\kappa/2]$. Let
\[
 \rho_\tau'+\widetilde r\rho_\tau=-\mathfrak a-\lambda,
 \qquad \rho_\tau(\tau)=\psi_1(\tau).
\]
For $s<\tau$, integration with $M(s)=e^{\int_0^s\widetilde r}$ gives
\[
 \psi_1(s)-\rho_\tau(s)
 =M(s)^{-1}\int_s^\tau M(q)(G_1(q)-\lambda)\,dq>0.
\]
At $s_0$ this gap is bounded below by a constant $c_0>0$ independent
of $\tau$. Since $v(s_0,\xi)\to\psi_1(s_0)$, for sufficiently small
$h>0$ one has $v(s_0,\xi)\ge\rho_\tau(s_0)$ for all $\xi\ge h^{-2}$,
uniformly in $\tau$.

On $(s_0,\tau)\times(h^{-2},\infty)$ set
\[
 \phi_h(s,\xi)=e^{-\int_{s_0}^s\widetilde r(q)dq}
                   h^3(2h^{-2}-\xi)_+^2,
 \qquad \zeta_\tau=\rho_\tau-\phi_h.
\]
At the left boundary, $\phi_h\ge c/h$, whereas $\rho_\tau$ and the
obstacles are uniformly bounded. Thus $\zeta_\tau\le v$ on the entire
parabolic boundary when $h$ is small. Moreover,
$\phi_{h,s}=-\widetilde r\phi_h$, $|\phi_{h,\xi}|\le Ch$, and
$|\phi_{h,\xi\xi}|\le Ch^3$ a.e., so
\[
 (\partial_s-\mathcal L_\xi )\zeta_\tau
 \le-\mathfrak a-\lambda+C(h+h^3)
 \le-\mathfrak a-C_He^{-\xi/\gamma_1}=\mathcal R
\]
after reducing $h$ once more. On $\{\zeta_\tau>v\}$,
$v<\zeta_\tau\le\rho_\tau<\psi_1$, and hence
$(\partial_s-\mathcal L_\xi )v\ge\mathcal R$. 
Applying the comparison principle \citep[Theorem 1.1]{YangYan2008} on that
set yields $\zeta_\tau\le v$ there, a contradiction. Thus $\zeta_\tau\le v$,
and at $s=\tau$ and $\xi\ge2h^{-2}$,
$\zeta_\tau=\psi_1$, proving contact. All choices were uniform in $\tau$,
so $\nu=\kappa/2$ and $\Xi_{\tau_0}=2h^{-2}$ suffice.
\end{proof}

\subsection{Proof of Proposition~\ref{prop:upper-contact-times}}
\begin{proof}
Finite upper contact implies $G_1>0$ by Lemma~\ref{lem:contact-location}
and interior limiting contact by Lemma~\ref{lem:finite-contact-interior}.
Conversely, Lemma~\ref{lem:recover-upper-contact} makes $\Xi_1<\infty$
on the stated intersection. Lower contact exists at every positive time,
and the positive obstacle gap and continuity then give $\Xi_1> \Xi_0>-\infty$.
\end{proof}

\subsection{Proof of Lemma~\ref{lem:upper-time-characterization}}
\begin{proof}
If $\tau_0\in\mathcal T_1$, then it lies inside $\mathcal C_1(\infty)$
and $G_1(\tau_0)>0$. Thus $F_1'=e^{\int_0^\tau\widetilde r}G_1>0$
near $\tau_0$. Choose $\tau_-<\tau_0$ in that neighborhood with
$[\tau_-,\tau_0]\subset\mathcal C_1(\infty)$ and the same last lower
contact $\ell$. The weak level and record conditions at $\tau_-$, combined
with strict increase to $\tau_0$, give (i)--(iii).

Conversely, (i)--(ii) imply upper limiting contact by
Lemma~\ref{lem:running-record}. The positive obstacle gap makes
$\ell_0(\tau)=\ell_0(\tau_0)=:\ell<\tau_0$ in a neighborhood of
$\tau_0$. Condition (i) persists there by continuity. By (iii), $F_1$
is strictly increasing on a larger neighborhood. For a small $\delta>0$,
condition (ii) and compactness give
\[
 F_1(\tau_0)>\max_{\ell\le s\le\tau_0-\delta}F_1(s).
\]
This strict gap persists at nearby times, while local increase handles
$s>\tau_0-\delta$. Thus the running-record condition also holds throughout
a neighborhood. Lemma~\ref{lem:running-record} gives interior limiting
contact, and (iii) gives $G_1(\tau_0)>0$.
Proposition~\ref{prop:upper-contact-times} concludes the proof.
\end{proof}

\subsection{Proof of Proposition~\ref{prop:semicontinuity-free-boundaries}}
\begin{proof}
At $\tau_0\in\mathcal T_0$, every $\xi>\Xi_0(\tau_0)$ satisfies
$v(\tau_0,\xi)>-\psi_0(\tau_0)$. Continuity preserves this strict inequality
for nearby times, so $\Xi_0(\tau)<\xi$ there. Let $\xi\downarrow \Xi_0(\tau_0)$.
The upper case follows identically from $v(\tau_0,\xi)<\psi_1(\tau_0)$
for $\xi<\Xi_1(\tau_0)$.
\end{proof}

\subsection{Proof of Lemma~\ref{lem:derivative-estimate}}
\begin{proof}
We localize the penalized derivative comparison of
\citet[Lemmas~5.4--5.5]{HaJeonOk2026}. The required localization and the
additional terms from time-dependent diffusion coefficients are recorded
below. No global ordering of the curves $\Lambda_i$ is used.
Write
\[
 A=\widetilde\sigma^2/2,\qquad B=A+\widetilde\beta-\widetilde r,
 \qquad P=\partial_\tau-A\partial_{\xi\xi}-B\partial_\xi +\widetilde r.
\]
Fix $\tau_0\in\mathcal T_0$. Semicontinuity, continuity of $v$, and
Lemma~\ref{lem:lower-contact-existence} give a rectangle
$Q=(\tau_-,\tau_+)\times(L,R)$ about the lower boundary such that
$v(\tau,L)=-\psi_0(\tau)$, $L<\Xi_0(\tau)<R<\Xi_1(\tau)$, and
the right side stays in a compact continuation subset. This construction
also works when $\Xi_1(\tau_0)=\infty$. Decrease $L$ so that
$C_He^{-L/\gamma_1}>\max_{[\tau_-,\tau_+]}G_0$.

We first justify the uniformly regular initial trace used in penalization.
Near lower contact $u=v+\psi_0$ satisfies
\[
 Pu=(G_0-C_He^{-\xi/\gamma_1})\one_{\{u>0\}}.
\]
Lemma~\ref{lem:contact-location} makes
$C_He^{-\xi/\gamma_1}-G_0$ strictly positive locally at every finite contact.
The time change $\varkappa=\int A$ puts this into the one-dimensional
obstacle form with Lipschitz coefficients and positive right-hand coefficient.
The local optimal regularity estimate of \citet{Blanchet-et-al-2006}, with
the correction in \citet{BDMErratum2010}, gives local bounds on $v_\tau$
and $v_{\xi\xi}$. Contact interiors and continuation interiors satisfy the same
local bounds directly. A finite cover and continuity of $v_\xi$ give
$v(\tau_-,\cdot)\in W^{2,\infty}(L,R)$, after a harmless adjustment of
the rectangle. Mollify this nondecreasing trace, preserving its lower bound
and leaving it unchanged near $L,R$, to obtain $g_\varepsilon$ with
uniform $W^{2,\infty}$ bounds and $g_\varepsilon'\ge0$.

With $\varphi$ the cutoff used in the existence proof, set
\[
 \rho_\varepsilon(\tau,q)=
       (C_He^{-L/\gamma_1}-G_0(\tau))\varphi(q/\varepsilon),
 \qquad Pv_\varepsilon=\mathcal R+
                   \rho_\varepsilon(\tau,v_\varepsilon+\psi_0),
\]
using lateral data $v$ and initial data $g_\varepsilon$.
The standard penalized construction and parabolic estimates
\citep[Theorem~9.1, p.~342]{Ladyzhenskaya} give
$v_\varepsilon\ge-\psi_0$, a bounded penalty, uniform $W^{1,2}_p$
bounds, and uniform convergence of $v_\varepsilon,v_{\varepsilon,\xi}$
to $v,v_\xi$. The chosen penalty makes the equation compatible with
$v=-\psi_0$ at $L$, and at $R$ it vanishes for small $\varepsilon$.
Thus the initial and lateral derivative traces are compatible.

Define $w=v_{\varepsilon,\xi}$, $k=v_{\varepsilon,\tau}+\psi_0'$, and
$P_\varepsilon=P-\rho_{\varepsilon,q}$, evaluating the penalty and its
derivatives at $(\tau,v_\varepsilon+\psi_0)$. Spatial differentiation gives
\[
 P_\varepsilon w=\mathcal R_\xi >0.
\]
On the parabolic boundary $w\ge0$: this follows from the obstacle at $L$,
the nondecreasing initial trace, and uniform convergence to $v_\xi>0$ at $R$.
Comparison gives $w\ge0$ in $Q$, with $w\ge c>0$ at $R$ for small
$\varepsilon$.

Time difference quotients, justified by Lipschitz coefficients, give
\[
 P_\varepsilon k=A'v_{\varepsilon,\xi\xi}+B'w-\widetilde r'v_\varepsilon
       +\mathcal R_\tau+\psi_0''+\widetilde r\psi_0'
       +\rho_{\varepsilon,\tau}.
\]
Eliminating $v_{\varepsilon,\xi\xi}$ with the penalized equation yields
\begin{align*}
 (P_\varepsilon-A'/A)k={}&(B'-(A'/A)B)w
       +(-\widetilde r'+(A'/A)\widetilde r)v_\varepsilon\\
 &+\mathcal R_\tau+\psi_0''+\widetilde r\psi_0'
       -G_0'\varphi((v_\varepsilon+\psi_0)/\varepsilon)\\
 &-(A'/A)(\psi_0'+\mathcal R+\rho_\varepsilon).
\end{align*}
All terms on the right are uniformly bounded. Since $A'/A$ is bounded and
$\rho_{\varepsilon,q}\le0$, the zero-order coefficient is bounded below
uniformly in $\varepsilon$. The boundary values of $k$ are zero at $L$,
bounded at $R$ by interior regularity, and bounded initially by the equation
and the uniform $W^{2,\infty}$ norm of $g_\varepsilon$. The estimate in \citet[Theorem~7.1]{Lieberman} gives $|k|\le C$.
The equation now bounds $v_{\varepsilon,\xi\xi}$, and hence
$|P_\varepsilon k|\le C$.

On the bounded rectangle $\mathcal R_\xi \ge c_1>0$. Therefore, for $C_1$
large enough independently of $\varepsilon$,
\[
 P_\varepsilon\{C_1w\pm(\tau-\tau_-)k\}\ge0.
\]
These quantities are nonnegative initially and at $L$. At $R$ increase
$C_1$ using $w\ge c>0$ and the bound on $k$. Comparison
\citep[Corollary~7.4]{Lieberman} gives
$(\tau-\tau_-)|k|\le C_1w$.
Pass to the weak limit in $v_{\varepsilon,\tau}$ and the uniform limit
in $v_{\varepsilon,\xi}$, and restrict to a smaller time interval about
$\tau_0$. Interior continuity in the continuation set then yields (i).

For (ii), Lemma~\ref{lem:recover-upper-contact} supplies a fixed right side
$R$ in upper contact and a left side in continuation. Enlarge $R$ so that
$G_1-C_He^{-R/\gamma_1}>0$ on the interval, and use
\[
 Pv_\varepsilon=\mathcal R-
\rho_{1,\varepsilon}(\tau,\psi_1-v_\varepsilon),
 \qquad
 \rho_{1,\varepsilon}(\tau,q)=
       (G_1(\tau)-C_He^{-R/\gamma_1})\varphi(q/\varepsilon).
\]
This sign makes the upper obstacle a supersolution and gives exact
compatibility at $R$. Apply the preceding argument to
$k=v_{\varepsilon,\tau}-\psi_1'$. The derivative of the signed penalty
with respect to $v_\varepsilon$ is again nonpositive. The same comparison
proves (ii).
\end{proof}

\subsection{Proof of Lemma~\ref{lem:non-degeneracy}}
\begin{proof}
At a lower boundary let $u=v+\psi_0$. The derivative estimate and smooth
fit give $u_\tau\to0$, $u_\xi\to0$, and $u\to0$ from continuation.
There, with $A=\widetilde\sigma^2/2$ and
$B=A+\widetilde\beta-\widetilde r$,
\[
 Au_{\xi\xi}=u_\tau-Bu_\xi+\widetilde r u-G_0+C_He^{-\xi/\gamma_1}.
\]
At $\tau_0$, Lemma~\ref{lem:contact-location} makes the last two terms
strictly positive at contact. Upper semicontinuity of $\Xi_0$, local boundedness
of the boundary, and continuity of the data preserve a uniform positive
margin for nearby times. Local H\"older continuity of $v,v_\xi$ and the
derivative estimate make the other terms uniformly small near contact.
Hence $u_{\xi\xi}\ge2c_0>0$ in a common continuation strip. Integrating twice
from $\Xi_0(\tau)$ proves both lower estimates.

For the upper boundary use $u=\psi_1-v$ and
$G_1-C_He^{-\Xi_1/\gamma_1}>0$, together with lower semicontinuity of $\Xi_1$.
The same reasoning gives $v_{\xi\xi}\le-2c_1$ near contact. Integrating from
$\xi$ to $\Xi_1(\tau)$ proves the upper estimates.
\end{proof}

\subsection{Proof of Lemma~\ref{lem:local-lipschitz-boundaries}}
\begin{proof}
We use the level-curve argument of
\citet[Theorem~5.6]{HaJeonOk2026} on a compact subinterval of the relevant
active component. For small $h>0$, strict spatial monotonicity defines
continuation curves by
\[
 v(\tau,\xi_0^h(\tau))+\psi_0(\tau)=h,
 \qquad \psi_1(\tau)-v(\tau,\xi_1^h(\tau))=h.
\]
The local rectangles in Lemma~\ref{lem:derivative-estimate} keep these
curves in a fixed spatial interval. The implicit function theorem gives
\begin{equation}\label{lv curve:uniform bound}
 (\xi_0^h)'=-\frac{v_\tau+\psi_0'}{v_\xi}(\tau,\xi_0^h),
 \qquad (\xi_1^h)'=-\frac{v_\tau-\psi_1'}{v_\xi}(\tau,\xi_1^h),
 \qquad |(\xi_i^h)'|\le C.
\end{equation}
As $h\downarrow0$, $\xi_i^h(\tau)\to \Xi_i(\tau)$ pointwise. Passing to
the limit in $|\xi_i^h(\tau)-\xi_i^h(s)|\le C|\tau-s|$ proves the claim.
\end{proof}

\subsection{Proof of Lemma~\ref{lem:boundary-C1alpha}}
\begin{proof}
We apply the boundary Harnack argument used in
\citet[Theorem~5.7]{HaJeonOk2026}, checking the hypotheses at our stated
coefficient regularity. Near a finite lower boundary choose a small cylinder
$Q_{2r}$ that excludes upper contact. The graph is Lipschitz in time by
Lemma~\ref{lem:local-lipschitz-boundaries}. In one spatial dimension it is
$C_p^1$ in the sense of \citet[Definition~1.1]{Kukuljan2022}, since
\[
 \sup_{0<|\tau-s|\le r^2}
 \frac{|\Xi_0(\tau)-\Xi_0(s)|}{|\tau-s|^{1/2}}\le Lr\longrightarrow0.
\]
In continuation set $w=v_\xi$ and $k=v_\tau+\psi_0'$. Both vanish
continuously at the boundary, by smooth fit and the derivative estimate.
For $P=\partial_\tau-A\partial_{\xi\xi}-B\partial_\xi +\widetilde r$,
\[
 Pw=\mathcal R_\xi ,\qquad
 Pk=A'v_{\xi\xi}+B'v_\xi-\widetilde r'v+
             \mathcal R_\tau+\psi_0''+\widetilde r\psi_0'.
\]
Time difference quotients justify the second equation a.e. The right-hand
sides are bounded: the derivative estimate bounds $v_\tau$, and the
continuation equation then bounds $v_{\xi\xi}$. Non-degeneracy gives
$w\ge c(\xi-\Xi_0(\tau))$, which also bounds $w$ below by a constant times
parabolic distance to the graph. A common time exponential removes the
zero-order term. Thus \citet[Theorem~1.2]{Kukuljan2022} gives
\[
 \|k/w\|_{C^{\alpha,2\alpha}(Q_r\cap\{v>-\psi_0\})}\le C,
 \qquad 0<\alpha<1/2.
\]
Here the exponents are in time and space, respectively. The bounded
right-hand sides allow every parabolic exponent below one.
By \eqref{lv curve:uniform bound}, evaluating this quotient along the
uniformly Lipschitz level curves gives uniform $C^\alpha$ bounds on
$(\xi_0^h)'$. Arzel\`a--Ascoli and $\xi_0^h\to \Xi_0$ yield
$\Xi_0\in C^{1,\alpha}$ locally. For upper contact use
$k=v_\tau-\psi_1'$ and $w\ge c(\Xi_1(\tau)-\xi)$, and the same argument applies.
\end{proof}

\subsection{Proof of Lemma~\ref{lem:boundary-behavior-free-boundaries}}
\begin{proof}
If $\Xi_0(\tau_n)\ge-M$ along $\tau_n\downarrow0$, then
$v(\tau_n,-M)=-\psi_0(\tau_n)$, contradicting continuity and
$v(0,-M)=\chi_1>-\psi_0(0)$. This proves the initial lower limit.

Now fix $\tau_*\in\partial\mathcal T_0\cap(0,T)$ and a finite $M>0$.
The entire slice at $\tau_*$ is in lower contact. Lemma
\ref{lem:contact-location} therefore gives $G_0(\tau_*)\le0$.
Choose
\[
 -G_0(\tau_*)<c_1<C_He^{-M/\gamma_1}-G_0(\tau_*).
\]
Near $\tau_*$ the upper limiting obstacle is inactive, so
$(v_\infty+\psi_0)'+\widetilde r(v_\infty+\psi_0)\ge G_0$.
Integration backwards from its zero at $\tau_*$ yields
\[
 \limsup_{\tau\uparrow\tau_*}
 \frac{v_\infty(\tau)+\psi_0(\tau)}{\tau_*-\tau}
 \le-G_0(\tau_*)<c_1.
\]
Choose a spatial radius $\rho$, preliminary time radius $\bar\delta$,
and $\kappa>0$ so that
$C_He^{-\xi/\gamma_1}-G_0(\tau)-c_1\ge\kappa$ on the corresponding
rectangle about $(\tau_*,M)$. With
$B=\widetilde\sigma^2/2+\widetilde\beta-\widetilde r$, first choose
$c_2>0$ satisfying
\[
 c_2(\|\widetilde\sigma^2\|_\infty+2\rho\|B\|_\infty
                  +\rho^2\|\widetilde r\|_\infty)\le\kappa/4.
\]
Then choose $\delta\le\bar\delta$ so that
$\|\widetilde r\|_\infty c_1\delta\le\kappa/4$ and
$v_\infty+\psi_0\le c_2\rho^2$ for $|\tau-\tau_*|\le\delta$.
On $Q=(\tau_*-\delta,\tau_*+\delta)\times(M-\rho,M+\rho)$ use
\[
 \Phi(\tau,\xi)=-\psi_0(\tau)+c_1(\tau_0-\tau)_++c_2(\xi-M)^2,
 \qquad \tau_0\in(\tau_*-\delta,\tau_*).
\]
The preceding choices give
$(\partial_\tau-\mathcal L_\xi )\Phi\ge\mathcal R$ uniformly in $\tau_0$:
the strict source margin absorbs the spatial quadratic terms and the
zero-order contribution of the time term. The strict limsup estimate allows
$\tau_0$ sufficiently close to $\tau_*$ that $v\le\Phi$ initially.
The bound $v+\psi_0\le v_\infty+\psi_0\le c_2\rho^2$ gives lateral
order. 
On $\{v>\Phi\}\cap Q$, one has $v>-\psi_0$, and hence
$(\partial_\tau-\mathcal L_\xi )v\le \mathcal R$. The comparison principle \citep[Theorem 1.1]{YangYan2008} therefore yields
$v\le\Phi$ on that set, a contradiction. At $\xi=M$ and $\tau\ge\tau_0$,
$\Phi=-\psi_0$, so contact holds. Thus $\Xi_0(\tau)\ge M$ for all nearby
active times. Since $M$ is arbitrary, the internal lower limit follows.

Finally, if $\tau^*\in\partial\mathcal T_1\cap(0,T)$ and
$\Xi_1(\tau_n)\le M$ along active times $\tau_n\to\tau^*$, then
$v(\tau_n,M)=\psi_1(\tau_n)$. Continuity gives finite upper contact
at $(\tau^*,M)$. Lemmas~\ref{lem:contact-location} and
\ref{lem:finite-contact-interior} imply $G_1(\tau^*)>0$ and interior
limiting contact, so Proposition~\ref{prop:upper-contact-times} puts
$\tau^*$ inside the open set $\mathcal T_1$, a contradiction.
\end{proof}

\subsection{Proof of Corollary~\ref{cor:free-boundary-graphs}}
\begin{proof}
At active times, continuity of $\Xi_i$ and the contact half-lines show that
points strictly inside a contact half-line are interior contact points,
and the graph itself separates contact from continuation. At a boundary time
of $\mathcal T_0$, the divergence $\Xi_0\to\infty$, together with full
lower contact outside $\mathcal T_0$, makes every finite point interior
to lower contact. The same holds at times outside
$\overline{\mathcal T_0}$. Thus no other lower free boundary pieces occur.
Finite upper contact occurs only on $\mathcal T_1$, so continuity of
$\Xi_1$ there proves the upper identity directly. All boundaries are relative
to $\Omega_T$.
\end{proof}

\subsection{Proof of Lemma~\ref{lem:initial-upper-boundary}}
\begin{proof}
Continuity gives $G_1>0$ near zero. Since $v_\infty(0)=\psi_1(0)$,
the reflected ODE keeps $v_\infty=\psi_1$ on an initial interval:
the nonnegative gap $e^{\int_0^\tau\widetilde r}(\psi_1-v_\infty)$
is nonincreasing there. Proposition~\ref{prop:upper-contact-times} gives
$(0,\delta)\subset\mathcal T_1$. Write $\xi_1=\Lambda_1(0)$.
The location bound implies $\liminf_{\tau\downarrow0}\Xi_1(\tau)\ge \xi_1$.

For the reverse bound fix $\varepsilon>0$. On a short initial interval,
$G_1(\tau)-C_He^{-\xi/\gamma_1}\ge2\lambda_\varepsilon>0$ for
$\xi\ge \xi_1+\varepsilon$. Define
\[
 \Phi(\tau,\xi)=\psi_1(\tau)
             -\eta_\varepsilon(\xi_1+2\varepsilon-\xi)_+^2,
 \qquad \xi\ge \xi_1+\varepsilon.
\]
For sufficiently small $\eta_\varepsilon>0$, bounded coefficients make
$(\partial_\tau-\mathcal L_\xi )\Phi\le\mathcal R$ a.e., and outside the
quadratic's support this follows directly from the strict source margin.
Initially $\Phi\le\psi_1(0)=v(0,\cdot)$, with a strict gap
$\eta_\varepsilon\varepsilon^2$ at the left corner. Continuity preserves
lateral order on a shorter interval. On $\{\Phi>v\}$ the upper obstacle
is inactive for $v$, so half-line comparison gives $\Phi\le v$.
Thus $v=\psi_1$ for $\xi\ge \xi_1+2\varepsilon$ near time zero, giving
$\limsup_{\tau\downarrow0}\Xi_1(\tau)\le \xi_1+2\varepsilon$.
Let $\varepsilon\downarrow0$.
\end{proof}

\subsection{Proof of Lemma~\ref{lem:remaining-terminal-cases}}
\begin{proof}
The first statement follows from continuity of $v_\infty$ and $\psi_1$.
The second follows from continuity of $G_1$ and the necessary condition
$G_1>0$ for finite upper contact. 
In the third case, either $\mathcal T_1\cap(0,\delta)=\varnothing$
for some $\delta>0$, or $0$ is an accumulation point of
$\mathcal T_1$. In the latter case, since
$\Lambda_1(\tau)=\gamma_1\log(C_H/G_1(\tau))<\Xi_1(\tau)$ on
$\mathcal T_1$ and $G_1(0)=0$, we obtain
$\Xi_1(\tau)\to+\infty$ as $\tau\downarrow0$ along $\mathcal T_1$.
\end{proof}

\subsection{Proof of Proposition~\ref{prop:minimum-financing}}
\begin{proof}
Write $\underline D=\min_{t\le u}D_H(t,u)>0$ and
$\underline k=\min_{t,i\ne j}k_{ij}(t)>0$.
The running and terminal parts of any deterministic cost are uniformly
bounded in absolute value. Each switch adds at least
$\underline D\underline k$. A minimizing sequence whose cost is no larger
than the no-switch cost plus one therefore has a uniform bound on its
number of switches. For each such number, allow weakly ordered times in
the compact time simplex. The cost is continuous there. Coincident
round trips can be deleted and strictly reduce cost, and any remaining switch
at $T$ can be deleted without increasing cost by
Lemma~\ref{lem:no-terminal-switch}. Hence a minimizing admissible finite
schedule exists. The same compact representation, or the deterministic
dynamic programming principle, gives continuity in the initial time.
For example, shifting switches from $(t,t+h)$ to $t+h$, deleting any
coincident round trips, changes the running cost by $O(h)$ and the
discounted switching costs by a common modulus of continuity times the
uniform switch bound. Immediate switching gives the corresponding
one-sided comparison. Boundedness follows from the no-switch upper bound
and the uniform lower bound on running and terminal costs.

The process
$L_{t,u}=\Gamma_{t,u}P_u/(pD_H(t,u))$ is the stochastic exponential with
Brownian coefficient $\sigma_P-\theta$. Bounded coefficients make it a
positive true martingale. Under the probability measure with density
$L_{t,T}$, the expected state-price housing cost divided by $p$ is the
expectation of the deterministic expression in
\eqref{eq:minimum-housing-cost}, evaluated at the realized schedule.
Every realized schedule has cost at least $m_i(t)$, and the deterministic
optimizer attains the bound. Conditional change of measure gives the
same conclusion at a stopping time. Finally, uniformly over policies,
the conjugate utility terms in $\mathcal J_i(t,y,p)$ are $O(y^{(g-1)/g})$.
Consequently $\mathcal J_i(t,y,p)/y\to-pm_i(t)$ as $y\to\infty$.
Since $\mathcal J_i=p^{1-g}j_i(t,yp^g)$, this implies
$j_i(t,z)/z\to-m_i(t)$ and hence the asserted identity for $v_\infty$.
\end{proof}

\subsection{Proof of Lemma~\ref{lem:switching-attainment}}
\begin{proof}
The verification argument follows \citet{JeonYang}, and we specify the
reconstruction, growth, and financing-limit checks for this model.
Put $s(z)=C_Hz^{-1/g}$, $f_i=\widehat u_i-\widehat a_i z$, and
$g_i(t)=G_i(T-t)$. Reconstruct the individual values by
\begin{align}\label{eq:individual-value-reconstruction}
 (\partial_t+\mathcal L_t)j_0
 &=-f_0-z(s(z)-g_0(t))\one_{\mathcal B},\notag\\
 (\partial_t+\mathcal L_t)j_1
 &=-f_1-z(g_1(t)-s(z))\one_{\mathcal S},
 \qquad j_i(T,z)=\widehat\Phi(z)-\chi_i z.
\end{align}
Feynman--Kac applies because the sources are locally bounded with growth
$C(z+|z^{p_*}|)$, $p_*=(g-1)/g$, integrable under the lognormal semigroup.
Parabolic estimates give continuous first spatial derivatives.
The difference solves the global linear equation for $Q$ with its terminal
data, so uniqueness gives $j_0-j_1=Q$. The contact-location inequalities
make both residuals in \eqref{eq:individual-value-reconstruction}
nonnegative before the minus sign, proving \eqref{eq:VI}.

Localized switching verification now applies: first-contact switching makes
the continuation drift zero and realizes value matching at each jump.
On every sample path, $Z$ stays in a compact subset of $(0,\infty)$,
and uniform continuity of $v(T-u,\log Z_u)$ and the strictly positive obstacle
gap exclude accumulation of alternating contacts, including at $T$.
Stopping after $N$ switches bounds expected deflated costs by running,
terminal, and continuation terms dominated by lognormal power suprema.
Monotone convergence for costs and dominated convergence for those terms
give attainment and cost integrability, exactly as in the standard
verification argument.

For each fixed policy the dual payoff is convex in $y$, with linear housing
costs. Its common terminal-conjugate term is strictly convex, so the supremum
is strictly convex, and differentiability follows from the reconstructed values.
Proposition~\ref{prop:minimum-financing} and the uniform bound on conjugate
utility imply
\[
 \mathcal J_i(t,y,p)=-pm_i(t)y+O(y^{p_*}),\qquad y\uparrow\infty.
\]
Convex secant bounds yield $\mathcal J_{i,y}\to-pm_i(t)$.
The common terminal term has second derivative $c_Ty^{-1-1/g}$ with
$c_T>0$, and subtracting it leaves a convex function. Integrating this
distributional second-derivative bound to infinity gives
$\mathcal J_{i,y}\le-pm_i(t)-gc_Ty^{-1/g}$, proving the other limit.
\end{proof}

\subsection{Proof of Theorem~\ref{thm:complete-verification}}
\begin{proof}
Lemma~\ref{lem:switching-attainment} gives the unique multiplier.
Use the martingale duality construction of \citet{JeonYang} with the
following housing-specific budget and solvency checks. Freeze the optimal
switching policy at $y^*$. Its differentiable dual payoff touches
$\mathcal J_i$ from below, so their derivatives agree. Lognormal moments
and integrable switching costs justify differentiation. Thus
$x=-\mathcal J_{i,y}$ is exactly the budget equality for this policy.

After any switch at $u$, define
\begin{align}\label{eq:wealth-replication}
 X_u^*=\E_u\bigg[&\int_u^T\Gamma_{u,q}
     (c_q^*+P_q\widehat a_{I_q^*}(q))\,dq+\Gamma_{u,T}W_T^*\notag\\
 &+\sum_{n\ge1:u<\tau_n<T}\Gamma_{u,\tau_n}P_{\tau_n}
                       k_{\iota_{n-1},\iota_n}(\tau_n)\bigg].
\end{align}
Brownian martingale representation finances this integrable conditional
claim with an effective portfolio. Subtracting $q_{I_u^*}P_u\zeta_P(u)$
gives the traded portfolio, and $W_u^*=X_u^*-\mathcal H(u)$ satisfies
the original wealth dynamics. At a switch the claim loses exactly $P_uk_{ij}$.
The conditional minimum-cost bound and strictly positive terminal surplus
give $X_u^*>P_um_{I_u^*}(u)$, including after each switch. The deflated
gain process is the uniformly integrable martingale generated by the total
budget claim. Fenchel equalities, budget equality, and dual attainment now
give admissibility and optimality by Proposition~\ref{prop:duality}.
\end{proof}

\Needspace{6\baselineskip}
\subsection{Proof of Proposition~\ref{prop:wealth-coordinates}}
\begin{proof}
The derivative range and strict convexity in
Lemma~\ref{lem:switching-attainment}, together with
$\mathcal J_i(t,y,p)=p^{1-g}j_i(t,yp^g)$, prove the claimed monotone
bijection. On the buying contact region,
$j_0=j_1-k_{01}z$ identically, so $-j_{0,z}+j_{1,z}=k_{01}$,
and the analogous identity holds throughout the selling region.
Continuity of the first derivatives extends these identities to the
boundaries and proves \eqref{eq:wealth-jump-identities}. If both
boundaries are finite, $a<b$ and strict increase of $w_1$ give
\eqref{eq:wealth-threshold-order} directly. The consumption identity is the
Fenchel optimizer. For the endpoint assertion, the estimates used in
Proposition~\ref{prop:minimum-financing} give
$j_i(t,z)=-m_i(t)z+O(z^{(g-1)/g})$ locally uniformly in $t$.
Convex secant bounds at $z/2,z,2z$ then give
$j_{i,z}(t,z)=-m_i(t)+O(z^{-1/g})$ locally uniformly. Continuity of $m_i$
proves the joint time-state endpoint limit. An all-buy slice maps onto
the entire admissible renter domain, while an absent selling region
remains absent. The coordinate map depends on both time and regime,
so its order preservation at each fixed time imposes no sign on the
time derivative of a transformed boundary.

Finally $s=C_HY^{-1/g}/P$ has Brownian coefficient
$s(\theta/g-\sigma_P)$. The diffusion coefficient of
$X=Pw_i(\tau,s)$ is
$P[w_i\sigma_P+s w_{i,s}(\theta/g-\sigma_P)]$.
Equating it to $\Sigma^{\mathsf T}\pi_i+q_iP\sigma_P$
gives \eqref{eq:portfolio-recovery}. Local Sobolev It\^o formulas justify
the calculation almost everywhere, and the financing martingale in
\eqref{eq:wealth-replication} defines the integrand up to stochastic
integration equivalence.
\end{proof}

\subsection{Proof of Proposition~\ref{prop:integral-verification}}
\begin{proof}
The global strong equation is
\[
 (\partial_\tau-\mathscr L_s)\mathcal V
 =d-s+(s-G_0)\mathbf1_{\{s\ge b\}}
                 +(s-G_1)\mathbf1_{\{s\le a\}}.
\]
The convention $a=0$ makes the last term zero on absent selling slices.
Spatial smooth fit eliminates boundary local-time terms. The lognormal
semigroup and Duhamel's formula yield \eqref{eq:Volterra-value},
and linear-growth moment bounds remove localization. Value matching gives
\eqref{eq:coupled-Volterra}. Conversely, uniqueness of the double-obstacle
solution identifies any reconstructed candidate satisfying all the stated
conditions, and then its contact sets. Value matching alone does not
establish uniqueness of arbitrary integral-equation roots.
\end{proof}

\section{Integral representation and numerical implementation}\label{sec:integral-method}
\subsection{Gaussian occupation identities}
Write $\delta=\widetilde r$, $d=-\mathfrak a$, $\varsigma=\widetilde\sigma/g$,
and
\[
 b_s=-\frac{\widetilde\sigma^2/2+\widetilde\beta-\widetilde r}{g}
       +\frac{\widetilde\sigma^2}{2g^2},\qquad
 \mathscr L_s=\tfrac12\varsigma^2s^2\partial_{ss}
                    +b_ss\partial_s-\delta.
\]
The coefficient $b_s$ is a drift, not a derivative of the buying boundary.
The difference $\mathcal V=h_0-h_1$ has generator $\mathscr L_s$, source
$d(\tau)-s$, initial value $\chi_1$ and obstacles $-\psi_0,\psi_1$.
For $0\le u<\tau$, set
\[
 D_{\tau u}=e^{-\int_u^\tau\delta(q)dq},\quad
 \mathsf v_{\tau u}=\int_u^\tau\varsigma(q)^2dq,\quad
 \mathsf m_{\tau u}=\int_u^\tau(b_s(q)-\varsigma(q)^2/2)dq,\quad
 E_{\tau u}=D_{\tau u}e^{\int_u^\tau b_s(q)dq}.
\]
With $\mathsf N$ the standard normal distribution function, the discounted
upper-tail probability and first moment are
\begin{align}\label{eq:Gaussian-moments}
 P^+_{\tau u}(s,B)&=D_{\tau u}\mathsf N\!\left(
    \frac{\log(s/B)+\mathsf m_{\tau u}}{\sqrt{\mathsf v_{\tau u}}}\right),\notag\\
 M^+_{\tau u}(s,B)&=sE_{\tau u}\mathsf N\!\left(
    \frac{\log(s/B)+\mathsf m_{\tau u}+\mathsf v_{\tau u}}{\sqrt{\mathsf v_{\tau u}}}\right).
\end{align}
The lower tails are $P^-_{\tau u}=D_{\tau u}-P^+_{\tau u}$ and
$M^-_{\tau u}=sE_{\tau u}-M^+_{\tau u}$. Here $D_{\tau u}$ is the
discount factor, $\mathsf m_{\tau u},\mathsf v_{\tau u}$ are the
log-increment mean and variance, and $M^\pm_{\tau u}$ are first moments,
and $\mathcal V$ is the obstacle value. At $B=0$, the upper tails are
$D_{\tau u},sE_{\tau u}$ and the lower tails are zero,
and at $B=\infty$ the conventions are reversed. Define
\begin{align}\label{eq:Volterra-kernel}
 K_{\tau u}(s)={}&M^+_{\tau u}(s,b(u))-G_0(u)P^+_{\tau u}(s,b(u))\notag\\
 &+M^-_{\tau u}(s,a(u))-G_1(u)P^-_{\tau u}(s,a(u)).
\end{align}
Then
\begin{equation}\label{eq:Volterra-value}
 \mathcal V(\tau,s)=D_{\tau0}\chi_1+
 \int_0^\tau\{D_{\tau u}d(u)-sE_{\tau u}+K_{\tau u}(s)\}\,du,
\end{equation}
and the coupled Volterra equations are
\begin{equation}\label{eq:coupled-Volterra}
 \mathcal V(\tau,b(\tau))=-\psi_0(\tau)\quad(\tau\in\mathcal T_0),
 \qquad
 \mathcal V(\tau,a(\tau))=\psi_1(\tau)\quad(\tau\in\mathcal T_1).
\end{equation}
Spatial smooth fit gives $\mathcal V_s=0$ at finite positive boundaries.

\begin{proposition}\label{prop:integral-verification}
The optimal boundaries satisfy
\eqref{eq:Volterra-value}--\eqref{eq:coupled-Volterra}.
Conversely, continuous positive curves on the classified active components,
with the prescribed endpoint conventions, identify the optimal boundaries
if their reconstructed $\mathcal V$ satisfies the obstacle inequalities,
the variational inequalities and the corresponding contact classification.
\end{proposition}

The same moments reconstruct the individual values. With
$e_i=gA_i/((1-g)C_H)$ and $e_T=g\eta^{1/g}/((1-g)C_H)$,
\begin{align}\label{eq:h-integral}
 h_0(\tau,s)={}&e_TsE_{\tau0}
 +\int_0^\tau\{e_0sE_{\tau u}-\widehat a_0(T-u)D_{\tau u}
                  +(M^+-G_0P^+)(s,b(u))\}\,du,\notag\\
 h_1(\tau,s)={}&e_TsE_{\tau0}-\chi_1D_{\tau0}
 +\int_0^\tau\{e_1sE_{\tau u}-\widehat a_1(T-u)D_{\tau u}
                  +(G_1P^--M^-)(s,a(u))\}\,du.
\end{align}
In the last terms all moments have indices $\tau u$, and $G_i=G_i(u)$.
The identities follow from the individual reconstruction in
\eqref{eq:individual-value-reconstruction}. Their derivatives give the
wealth and portfolio formulas. This specializes the occupation-integral method of \citet{JeonYang}
to the housing normalization and the active time components characterized
here, using the time-dependent obstacle framework of \citet{HaJeonOk2026}.

\subsection{Computation and accuracy}
We use the limiting obstacle ODE to insert contact endpoints and zero-incentive
times into the grid. A positive Gaussian-integral dynamic program initializes
the branch, and four-point product integration solves the coupled boundary
equations. A square-root time transformation resolves the short-time kernel,
and positive buying thresholds are interpolated reciprocally. Value matching is
fitted with a vanishing scaled smooth-fit residual. At an all-buy slice,
the exact constant value allows a Duhamel restart. No failed root is
replaced by a previous boundary or by zero. The implementation uses
constant financial and house-price coefficients and the deterministic
housing-cost paths specified in the main text, but the representation above
also permits deterministic time-varying diffusion coefficients.

The production and comparison grids have nominal resolutions 1,200 and
600, with exact event nodes added. Independent verification uses an
implicit monotone PDE scheme with 8,000 time steps and log-state spacing
0.002. Boundary comparisons sample $\tau\in[0.5,20]$, omit 0.15-year
neighborhoods of internal endpoints, and retain finite positive thresholds
in $s\in[10^{-4},0.15]$. Value comparisons include contact slices at
$\tau=1,5,10,15,20$ and $s\in[10^{-5},0.15]$.

\begin{table}[htbp]\centering\small
\caption{Absolute numerical discrepancies. VM denotes unscaled Volterra
value matching. Boundary discrepancies are in $s$, and value discrepancies
are in $Q/z$.}\label{tab:verification}
\begin{tabular}{lrrrr}
\toprule
Case & VM residual & Grid change & PDE boundary & PDE value\\
\midrule
1 & 1.30e-08 & 1.01e-04 & 2.91e-04 & 2.45e-05\\
2 & 2.53e-09 & 1.77e-05 & 1.93e-04 & 1.01e-05\\
3 & 3.97e-10 & 1.88e-05 & 2.04e-04 & 5.47e-06\\
4 & 3.19e-09 & 2.35e-05 & 1.90e-04 & 6.21e-06\\
5 & 8.81e-09 & 6.83e-05 & 1.93e-04 & 2.58e-05\\
6 & 7.09e-06 & 5.36e-05 & 1.74e-04 & 6.73e-05\\
7 & 4.80e-06 & 6.15e-05 & 1.53e-04 & 4.86e-05\\
8 & 1.77e-09 & 1.80e-05 & 2.77e-04 & 1.98e-05\\
9 & 6.46e-09 & 3.00e-05 & 1.13e-04 & 2.65e-05\\
10 & 5.05e-09 & 1.92e-05 & 1.56e-04 & 4.16e-05\\
11 & 3.56e-08 & 6.70e-06 & 1.44e-04 & 1.18e-05\\
\bottomrule
\end{tabular}

\end{table}

Table~\ref{tab:verification} reports the verification discrepancies.
The largest VM residual is $7.09\times10^{-6}$ and the largest sampled
PDE value difference and obstacle violation are $6.73\times10^{-5}$.
At ten remaining years in the baseline, buy and sell wealth-jump errors
are $4.83\times10^{-8}$ and $6.42\times10^{-9}$ property values.
Across the four wealth panels at the five verification horizons, the
largest wealth-jump discrepancy is $2.80\times10^{-5}$, and evaluating wealth
at $s=10^{-8}$ differs from the independently computed financing floor by
at most $6.02\times10^{-7}$. These are measured discretization discrepancies,
not certified bounds on continuous-problem errors.

\section{Additional contact configurations}\label{ec:cases}

This section describes Cases 8, 10 and 11 of Section~\ref{sec:numerical},
which isolate critical or equality configurations, and shows the dual
boundaries for Cases 7--11. Notation and parameters are as in
Section~\ref{sec:numerical}, and the contact times are listed in
Table~\ref{tab:components}.

\begin{enumerate}[leftmargin=*,itemsep=4pt]
\setcounter{enumi}{7}
\item \textbf{Zero terminal sale incentive.}
Set $d(\tau)=0.025(0.04)+0.012(1-e^{-\tau/2})$ and $\chi_1=0.04$.
Here $G_1(0)=0$ but $G_1(\tau)>0$ for $\tau>0$. Selling is possible
at every positive remaining time, with its wealth region shrinking to the
owner's financing floor as $\tau\downarrow0$.
\setcounter{enumi}{9}
\item \textbf{An isolated loss of selling.}
Set $d(\tau)=0.025(0.04)+0.012(1-\mathcal G_{10,1.2}(\tau))$ and
$\chi_1=0.04$. Although $v_\infty=\psi_1$ throughout, $G_1(10)=0$
and $G_1>0$ elsewhere. Selling disappears only at $\tau=10$, illustrating
why limiting contact alone is insufficient.

\item \textbf{An isolated all-buy time.}
Prescribe $v_\infty(\tau)=-0.05+0.065\cos^2(\pi\tau/20)$ and set
$d=v_\infty'+0.025v_\infty$, with $\chi_1=0.015$.
The limiting value touches the lower obstacle only at $\tau=10$,
and at that time every admissible renter buys.
\end{enumerate}

\begin{figure}[p]\centering
\includegraphics[width=\textwidth,height=.82\textheight,keepaspectratio]{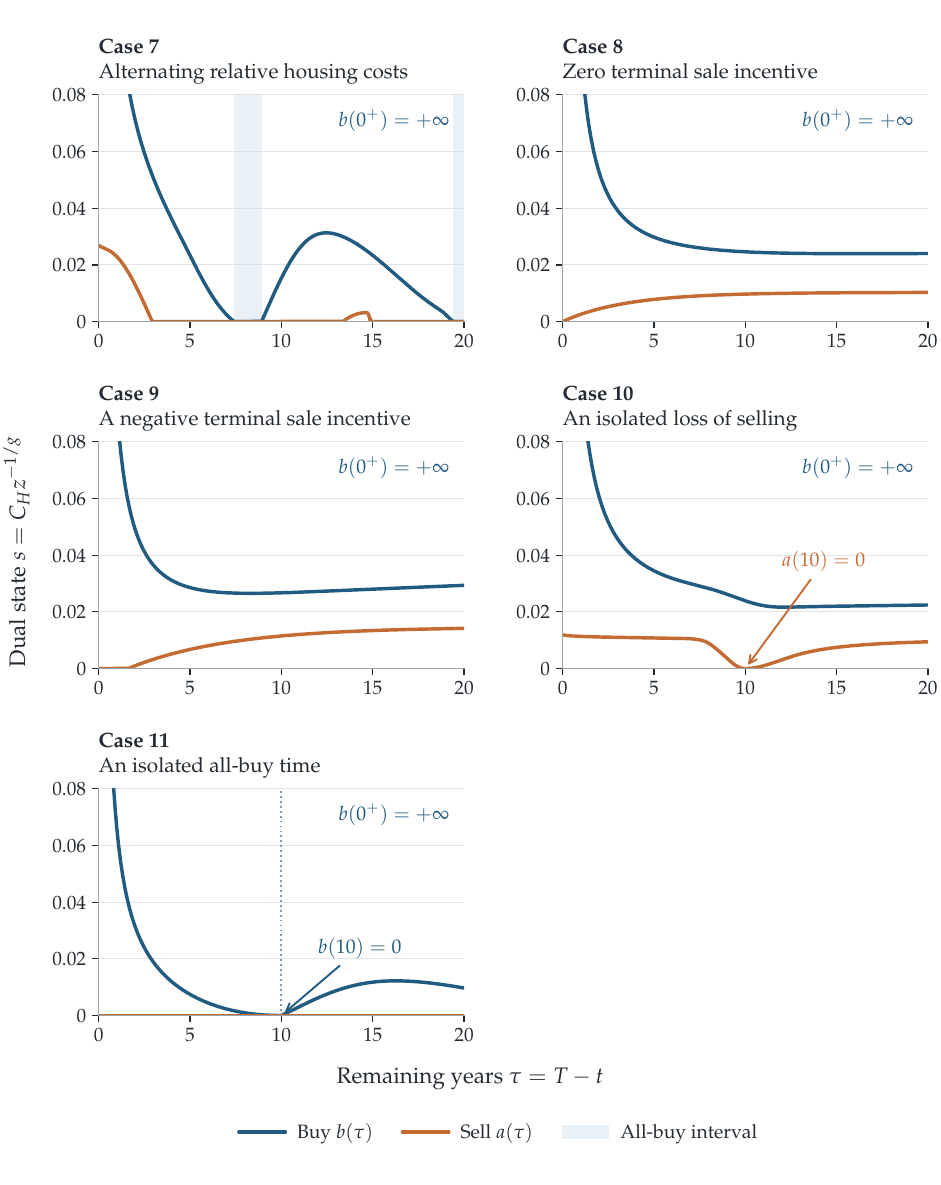}
\caption{Cases 7--11. The same conventions as in Figure~\ref{fig:cases1}
apply. Cases 10 and 11 distinguish an isolated missing sell boundary from
an isolated all-buy slice. Both phenomena occur with smooth coefficients
and strictly positive transaction costs.}\label{fig:cases2}
\end{figure}

In Case 10, $v_\infty=\psi_1$ throughout the horizon, but $G_1(10)=0$
removes the finite sell boundary at that time. In Case 11, a single contact
of the limiting value with the lower obstacle forces the entire finite-state
value onto that obstacle.

\section{Additional numerical results}\label{ec:numerics}

\subsection{Risk-adjusted housing carry}
\begin{figure}[htbp]\centering
\includegraphics[width=.86\textwidth]{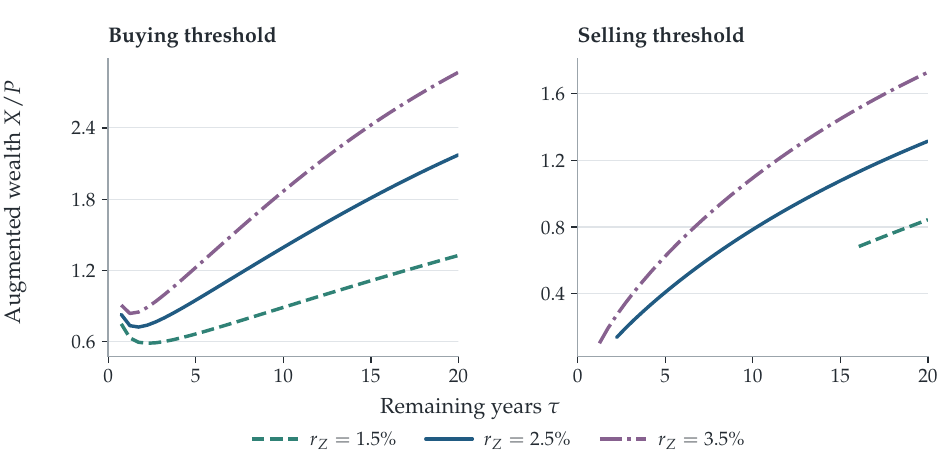}
\caption{Risk-adjusted housing carry. The changes in $r_Z$, $\beta_Z$
and the normalized-state drift are all recomputed from the primitive
house-price drift, and no coefficient is held artificially fixed.}\label{fig:carry}
\end{figure}

For a fixed normalized generator and fixed obstacles, increasing the
effective ownership premium $d$ increases $\mathcal V$ by comparison.
Consequently the high-$s$ buying region contracts and the low-$s$ selling
region expands: both $b$ and $a$ weakly increase. This is an exact
comparative-static statement in the normalized state. Changes in primitive
house-price drift also change the normalized generator, and must be
computed jointly. Figure~\ref{fig:carry} does this by varying $\mu_P$
while holding financial risk prices, volatilities, rental yield and owner
cash expenses fixed. Increasing $r_Z=-\kappa_H$ from 1.5 to 3.5 percent
raises the effective owner cost. In the displayed experiments, a household
requires more wealth before buying and sells at a higher wealth level.
Thus spanning does not make expected housing carry economically irrelevant.

At ten remaining years, the purchase threshold rises from $0.8893$ to
$1.8672$ property values across this carrying-cost comparison. With the
lower carrying cost there is no sale region at that horizon. With the
higher carrying cost the pre-sale threshold is $1.0941$. The difference
therefore affects both the location of a decision threshold and whether
selling is ever optimal at a given remaining horizon.

\subsection{Consumption and the financial hedge}
\begin{figure}[htbp]\centering
\includegraphics[width=.86\textwidth]{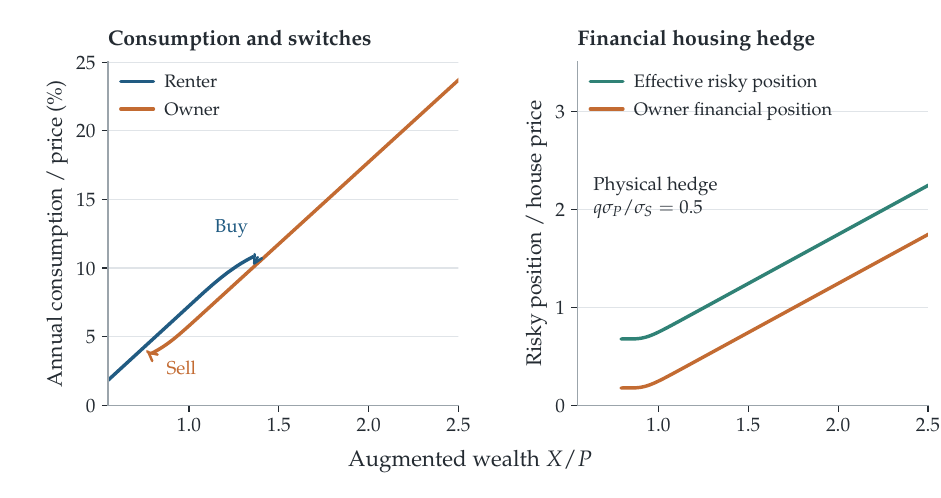}
\caption{Optimal consumption and the owner's financial hedge at ten
remaining years. Consumption curves are displayed in their respective
continuation regions. Arrows indicate the tenure change and associated
wealth loss.}\label{fig:policies}
\end{figure}

The left panel of Figure~\ref{fig:policies} shows consumption in the
two continuation regimes at ten remaining years. A purchase changes
consumption by the exact factor
$A_1/A_0=0.894427$, while a sale changes it by the reciprocal factor.
The household exchanges part of perishable consumption for the larger
housing-service bundle at a purchase, and pays the transaction cost through
the wealth jump. These changes are evaluated at the same continuous shadow
price. They are not jumps in the pricing kernel.

The right panel separates the effective risky exposure from the owner's
actual financial investment. Their difference is exactly
$q_1\sigma_P/\sigma_S=0.5$ property values in the one-stock benchmark.
The physical home supplies this Brownian exposure, so the financial position
offsets it. The figure illustrates why the owner-exposure neutralization
must be interpreted as a portfolio reparametrization, not as an absence of
housing hedging or housing costs.

\section{Portfolio constraints and uninsurable cost risk}\label{ec:robust}

Theorem~\ref{thm:timing-windows} rests on two features of the benchmark
market. Positions are unrestricted, so an owner can offset the house-price
exposure of the home at every wealth level, and every cost schedule is known
in advance. This section removes each feature in turn and recomputes the
thresholds. The first experiment imposes a short-sale or a borrowing
restriction. The second makes the ownership cost subject to a permanent
increase at a random date that cannot be insured. Both experiments keep the
preferences, services, and market prices of Table~\ref{tab:parameters} and
use the baseline and Case 5 cost schedules. Labor income is zero, so $X=W$,
and the state is $z=W/P$.

\subsection{Primal formulation and numerical method}\label{ec:primal}

With $\gamma_1=2$, flow utility is $-A_i/u$ for $u=c/P$ and
$A_i=H_i^{\gamma_1-\gamma}$, terminal utility is $-\eta/(z-\chi_i)$, and
$V_i(t,W,P)=P^{1-\gamma_1}v_i(t,z)$. Let $\varpi=\pi/P$ and let
$\nu=\sigma_S\varpi+(q_i-z)\sigma_P$ be the diffusion coefficient of $z$.
Between switches the normalized value satisfies, in remaining time $\tau$,
\begin{equation}\label{eq:primal-hjb}
 \partial_\tau v_i=\sup_{u>0,\ \nu\in\mathcal N_i(z)}
 \Big\{(r_Zz-\widehat a_i+\varrho\,\nu-u)\,\partial_zv_i
 +\tfrac12\nu^2\partial_{zz}v_i-\frac{A_i}{u}\Big\}-\beta_Zv_i,
 \qquad \varrho=\theta-\gamma_1\sigma_P,
\end{equation}
with $r_Z$ and $\beta_Z$ as in Section~\ref{sec:numerical}, together with
the switching conditions $v_0(\tau,z)\ge v_1(\tau,z-k_{01})$ and
$v_1(\tau,z)\ge v_0(\tau,z-k_{10})$ whenever the destination is
financeable. In the benchmark $\mathcal N_i(z)=\mathbb R$, and the
supremum over $\nu$ reproduces the operator of Section~\ref{sec:model}
with $\sigma_Z=|\varrho|$.

Let $f_i(\tau)$ denote the financing floor of tenure $i$, let $y=z-f_i$ be
the surplus over the floor, $\xi=\log y$, and $h=-yv_i$. The function $h$ is
bounded, and at the floor it tends to the Merton constant of the household
that follows the minimum cost schedule with zero surplus, which is known in
closed form. The equation for $h$ is discretized on a uniform grid in $\xi$
between $y=10^{-5}$ and $y=30$ with implicit Euler steps in $\tau$. At each
step the controls are updated by policy iteration. First derivatives use
the one-sided difference consistent with the sign of the drift, blended
with the central difference wherever the resulting matrix remains
monotone, and the switching conditions enter through an active set. The
scaled exposure $|\nu|/y$ is capped at $M=1$, ten times the unrestricted
optimum $|\varrho|/2$. A portfolio restriction takes precedence over the
cap. Where an owner with $z<q_1$ is required to carry a scaled exposure
larger than $M$, the control set is the restriction alone and the exposure
is the smallest admissible one, so the cap never empties the control set
and never forces a sale, and applying the same rule or a fixed cap to the
renter changes no threshold. Larger caps let the discrete scheme exploit the discrete
non-concavity that a switching condition creates next to a free boundary,
and the discrete value of a renter then falls below the value of renting
forever, which is a feasible policy. With $M\le2$ this lower bound holds to
a relative accuracy of $8\times10^{-4}$ at every date at which we store the
value functions, and $3\times10^{-4}$ at ten remaining years. At that date
the cap binds at no continuation node with $z\le8$ in the baseline, and at
every date it binds at no more than five such nodes, all adjacent to a free
boundary.

In the benchmark the primal solver returns buying and selling thresholds of
$1.3929$ and $0.7853$ at ten remaining years, against $1.3940$ and $0.7852$
from the integral equations, and it agrees with the dual finite differences
of Section~\ref{sec:integral-method} within $0.003$ at every remaining time
of at least one year in both cost schedules, and within $0.005$ at half a
year. The grid is $\Delta\xi=0.01$ and
$\Delta\tau=0.0025$ unless stated otherwise.

\subsection{Separate short-sale and borrowing restrictions}\label{ec:pc}

The first restriction forbids short positions in the risky security,
$\pi\ge0$, and leaves borrowing in the risk-free account free. The second
forbids borrowing in the risk-free account, $W-q_iP-\pi\ge0$, and leaves
short sales free. The two are imposed separately, and neither is a mortgage
contract with a loan-to-value limit. In terms of $\nu$,
\begin{equation}\label{eq:pc-sets}
 \mathcal N_i^{\rm NS}(z)=[\sigma_P(q_i-z),\infty),
 \qquad
 \mathcal N_i^{\rm NB}(z)=(-\infty,(\sigma_S-\sigma_P)(z-q_i)].
\end{equation}
Under either restriction an owner with $z<q_1$ cannot make the diffusion
coefficient of $z$ vanish. The only sure way to remain solvent is to sell
before the surplus over the renter's floor plus the sale cost is exhausted,
so the financeable interior of the owner's domain is
$z>m_0(\tau)+k_{10}(\tau)$ at every date, including the dates at which the
unrestricted floor $m_1$ is lower. At this floor the restricted owner sells.

\begin{table}[htbp]\centering\small
\caption{Thresholds $X/P$ under separate portfolio restrictions. NS: no
short sale. NB: no borrowing. At two remaining years in the baseline and at
seven in Case 5 no unrestricted owner sells.}\label{tab:pc}
\begin{tabular}{llrrrrrr}
\toprule
 & & \multicolumn{2}{c}{Unrestricted} & \multicolumn{2}{c}{NS} & \multicolumn{2}{c}{NB}\\
Schedule & $\tau$ & Buy & Sell & Buy & Sell & Buy & Sell\\
\midrule
Baseline & 10 & 1.393 & 0.785 & 1.556 & 0.825 & 1.887 & 1.073\\
Baseline & 2  & 0.729 & none  & 0.801 & 0.333 & 0.994 & 0.501\\
Case 5   & 10 & 1.502 & 0.868 & 1.563 & 0.903 & 1.879 & 1.220\\
Case 5   & 7  & 1.250 & none  & 1.327 & 0.506 & 1.638 & 0.657\\
\bottomrule
\end{tabular}
\end{table}

\begin{figure}[htbp]\centering
\includegraphics[width=.9\textwidth]{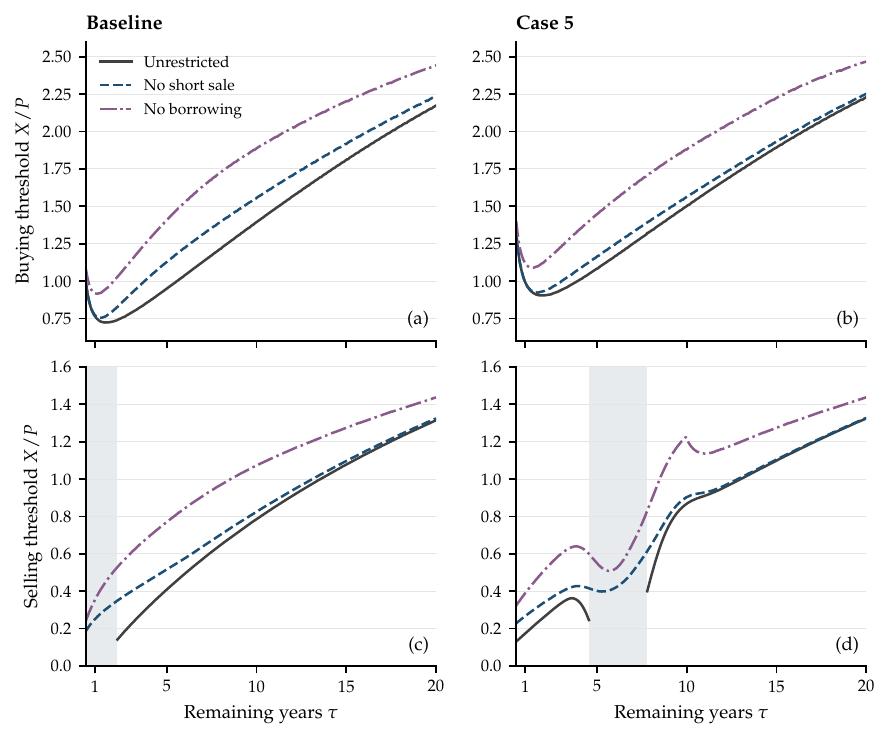}
\caption{Buying and selling thresholds under separate portfolio
restrictions. Shaded intervals are the dates without an unrestricted
selling region. Curves are shown from half a year before the horizon.}
\label{fig:pc}
\end{figure}

Table~\ref{tab:pc} and Figure~\ref{fig:pc} report the thresholds. Both
restrictions raise the buying threshold and the selling threshold, and the
borrowing restriction has the larger effect. At ten remaining years in the
baseline the buying threshold rises from $1.393$ to $1.556$ without short
sales and to $1.887$ without borrowing, and the selling threshold rises from
$0.785$ to $0.825$ and $1.073$.

The change that matters for Theorem~\ref{thm:timing-windows} concerns the
dates. In the benchmark no owner sells during the last $2.21$ years in the
baseline and between $4.57$ and $7.80$ remaining years in Case 5. Under
either restriction a selling region exists at every date. At two remaining
years the restricted selling thresholds are $0.333$ and $0.501$ property
values against a floor of $0.128$, and at seven remaining years in Case 5
they are $0.506$ and $0.657$ against a floor of $0.379$. The reason is the
floor. A restricted owner with less than one house value of net worth
carries an exposure it cannot remove, and the only sure protection of
solvency is a sale before wealth reaches $m_0+k_{10}$. This sale is made to
stay solvent, not to save on housing costs, and the deterministic cost
comparison that decides the dates in the benchmark does not apply to it.
The windows of Theorem~\ref{thm:timing-windows} are therefore a property
of a market in which the house-price exposure can be offset at every
wealth level.

An exposure cap on its own would force a sale whenever
$y<\sigma_P(1-f)/(M+\sigma_P)$, where $f=m_0+k_{10}$ is the restricted
floor, because no exposure within the cap satisfies the restriction there.
Since the restriction takes precedence (Section~\ref{ec:primal}), no such
band exists in the computation, and the computed selling thresholds lie
above it by a factor of at least $1.5$ at every date from half a year before
the horizon. Table~\ref{tab:pc-verify} reports the sensitivity of the
thresholds to the grid and to the cap at ten remaining years and at the
dates at which the unrestricted owner does not sell. Halving both steps
changes the thresholds by at most $0.005$ and doubling them by at most
$0.018$. Between $M=0.5$ and $M=1$ no selling threshold changes by more
than $0.001$ and no buying threshold by more than $0.006$. Raising the cap
to $M=2$ moves the borrowing-restricted buying threshold by $0.021$ at
seven remaining years in Case 5, which is the artifact described in
Section~\ref{ec:primal}, and leaves the selling thresholds within
$0.001$.

\begin{table}[htbp]\centering\small
\caption{Grid and cap sensitivity of the restricted thresholds $X/P$ at
ten remaining years and at the dates without an unrestricted selling
region.}\label{tab:pc-verify}
\begin{tabular}{llcccc}
\toprule
Schedule & $\tau$ & $(\Delta\xi,\Delta\tau)$ & $M$ & NS buy / sell & NB buy / sell\\
\midrule
Baseline & 10 & $(0.02,\ 0.005)$ & 1 & 1.558 / 0.826 & 1.901 / 1.079\\
 &  & $(0.01,\ 0.0025)$ & 1 & 1.556 / 0.825 & 1.887 / 1.073\\
 &  & $(0.005,\ 0.00125)$ & 1 & 1.557 / 0.825 & 1.884 / 1.070\\
 &  & $(0.01,\ 0.0025)$ & 0.5 & 1.556 / 0.825 & 1.887 / 1.073\\
 &  & $(0.01,\ 0.0025)$ & 2 & 1.556 / 0.825 & 1.884 / 1.073\\
\addlinespace
 & 2 & $(0.02,\ 0.005)$ & 1 & 0.801 / 0.334 & 1.004 / 0.502\\
 &  & $(0.01,\ 0.0025)$ & 1 & 0.801 / 0.333 & 0.994 / 0.501\\
 &  & $(0.005,\ 0.00125)$ & 1 & 0.800 / 0.333 & 0.990 / 0.501\\
 &  & $(0.01,\ 0.0025)$ & 0.5 & 0.801 / 0.333 & 0.988 / 0.501\\
 &  & $(0.01,\ 0.0025)$ & 2 & 0.801 / 0.333 & 0.996 / 0.501\\
\addlinespace
Case 5 & 7 & $(0.02,\ 0.005)$ & 1 & 1.327 / 0.506 & 1.656 / 0.656\\
 &  & $(0.01,\ 0.0025)$ & 1 & 1.327 / 0.506 & 1.638 / 0.657\\
 &  & $(0.005,\ 0.00125)$ & 1 & 1.327 / 0.506 & 1.633 / 0.657\\
 &  & $(0.01,\ 0.0025)$ & 0.5 & 1.327 / 0.506 & 1.638 / 0.658\\
 &  & $(0.01,\ 0.0025)$ & 2 & 1.327 / 0.506 & 1.617 / 0.656\\
\bottomrule
\end{tabular}
\end{table}

\subsection{An uninsurable increase in the ownership cost}\label{ec:event}

Let $N$ be a Poisson process with intensity $\lambda$, independent of $B$.
At its first jump the owner's flow cost rises permanently from $a_1$ to
$a_1+\Delta$. The jump is not traded, so the market is incomplete, while
Brownian risk remains spanned and positions remain unrestricted. Let
$v_i^{\rm post}$ denote the value after the jump, which is the benchmark
problem with cost $a_1+\Delta$, and $v_i^{\rm pre}$ the value before it.
Admissibility requires solvency in every state, and the jump can occur at
any time, so the floors before the jump are those of the high-cost
schedule, $m_0$ and $m_1^{\rm high}$, computed from
\eqref{eq:minimum-housing-cost} with $a_1+\Delta$. The pre-event value
solves \eqref{eq:primal-hjb} with the additional term
$\lambda(v_i^{\rm post}-v_i^{\rm pre})$ and the same switching conditions.
Both problems are solved on the same grid, with the post-event solution
entering the pre-event equation as a source. The boundary value at the
owner's floor differs from the benchmark. Where the high-cost schedule holds
the home, the pre-event owner pays $\Delta$ less than the schedule that
defines its floor, so with zero consumption its surplus grows at the floor,
and $h=-yv_1^{\rm pre}$ tends to zero like
$(\lambda\widetilde K_1/\Delta)\,y\log(1/y)$ instead of tending to a
Merton constant. Where the high-cost schedule sells immediately, the
pre-event owner sells at the floor. The computation uses these values.
Replacing the boundary value in the hold regime by the Merton constant of
the high-cost schedule leaves every date in Table~\ref{tab:event} and
every threshold unchanged to at least four digits, because the surplus
drifts away from the floor there and the boundary value does not enter the
interior.

Theorem~\ref{thm:timing-windows} does not apply before the jump, and the
question is whether the dates at which some owner sells are still decided
by deterministic data. Two deterministic sets are natural candidates. The
first is the contact set of the high-cost schedule, that is, the dates at
which $m_1^{\rm high}=m_0+k_{10}$ in the interior. The second is the set of
dates at which the sale incentive $G_1$ of Lemma~\ref{lem:contact-location},
computed with the current cost $a_1$, is strictly positive. Their
intersection is denoted by $\mathcal D$. Near the worst-case floor an owner
who holds the home sees its surplus over that floor shrink at the rate
$G_1$ computed with the current cost, and it must sell before the surplus
is exhausted, since the jump may arrive at any moment. The set
$\mathcal D$ therefore collects the dates at which solvency forces a sale
near the floor.

\begin{table}[htbp]\centering\small
\caption{Dates without a selling region, in remaining years, with an
uninsurable cost increase $\Delta=0.01$ arriving at rate $\lambda=0.1$.
Numerical intervals are the dates at which the selling region extends
less than $0.0001$ property values above the floor, read on a grid of
$0.05$ years. The deterministic intervals are the complements of the
contact set of the relevant schedule for the two benchmark problems and the
complement of $\mathcal D$ for the pre-event problem. The rows below the
rule change one parameter at a time.}\label{tab:event}
\begin{tabular}{llll}
\toprule
Schedule & Problem & No selling (numerical) & No selling (deterministic)\\
\midrule
Baseline & current cost, no event & $(0,\,2.20)$ & $(0,\,2.21)$\\
 & before the event & $(0,\,1.15)$ & $(0,\,1.17)$\\
 & after the event & $(0,\,1.15)$ & $(0,\,1.17)$\\
Case 5 & current cost, no event & $(4.60,\,7.75)$ & $(4.57,\,7.80)$\\
 & before the event & $(4.60,\,6.85)$ & $(4.57,\,6.86)$\\
 & after the event & $(5.25,\,6.85)$ & $(5.22,\,6.86)$\\
\midrule
Baseline, $\gamma=5$, $\alpha=1/4$ & before the event & $(0,\,1.15)$ & $(0,\,1.17)$\\
Baseline, $H_0=0.9$ & before the event & $(0,\,1.15)$ & $(0,\,1.17)$\\
Baseline, $\beta=0.05$ & before the event & $(0,\,1.15)$ & $(0,\,1.17)$\\
Baseline, $\eta=2$ & before the event & $(0,\,1.15)$ & $(0,\,1.17)$\\
Baseline, $\lambda=0.3$ & before the event & $(0,\,1.15)$ & $(0,\,1.17)$\\
Baseline, $\Delta=0.02$, $\lambda=0.2$ & before the event & $(0,\,0.80)$ & $(0,\,0.80)$\\
Case 5, $\gamma=5$, $\alpha=1/4$ & before the event & $(4.60,\,6.85)$ & $(4.57,\,6.86)$\\
Case 5, $H_0=0.9$ & before the event & $(4.60,\,6.85)$ & $(4.57,\,6.86)$\\
\bottomrule
\end{tabular}
\end{table}

\begin{figure}[htbp]\centering
\includegraphics[width=.9\textwidth]{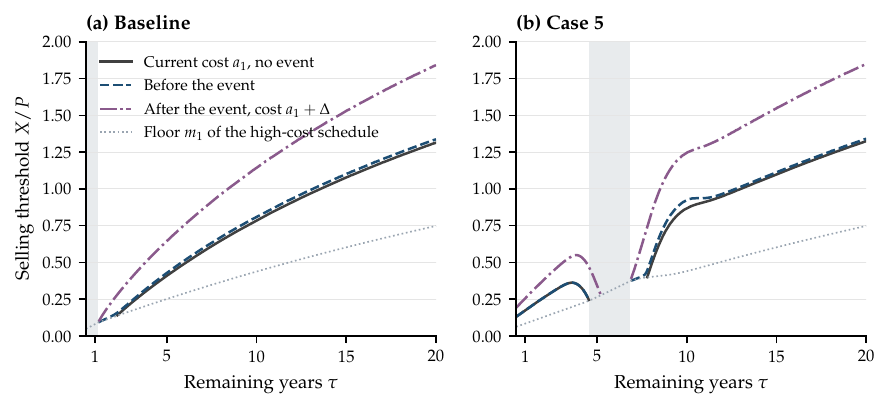}
\caption{Selling thresholds with an uninsurable cost increase,
$\Delta=0.01$ and $\lambda=0.1$. Shaded intervals are the dates outside
$\mathcal D$, at which the worst-case schedule does not sell now or the
current-cost sale incentive is not positive.}\label{fig:event}
\end{figure}

Table~\ref{tab:event} and Figure~\ref{fig:event} report the results for
$\Delta=0.01$ and $\lambda=0.1$. In the baseline no owner sells during the last
$2.21$ years at the current cost and during the last $1.17$ years at the
high cost. Before the event the selling region appears at $1.17$ remaining
years, the date at which the high-cost schedule starts to sell immediately,
and not at $2.21$. In Case 5 the current-cost schedule has no selling
region between $4.57$ and $7.80$ remaining years and the high-cost schedule
has none between $5.22$ and $6.86$. Before the event there is no selling
region between $4.57$ and $6.86$. The window closes at $4.57$, where the
current-cost sale incentive changes sign, and it reopens at $6.86$, where
the high-cost schedule resumes selling immediately. In both schedules the
dates without a pre-event selling region coincide with the complement of
$\mathcal D$ at the resolution of the computation. They do not change when the
lower truncation is moved from $10^{-5}$ to $10^{-6}$, when the cap is
lowered to $0.5$, when both grid steps are halved, or when the detection
level of Table~\ref{tab:event} is moved between $10^{-5}$ and $10^{-4}$,
while at $10^{-3}$ the pre-event opening date moves by one grid step
because the pre-event region grows slowly from the floor. They do not change
when $\gamma$ and $\alpha$ are moved to $5$ and $1/4$, when $H_0$ is raised
to $0.9$, when $\beta$ is raised to $0.05$, or when $\eta$ is raised to
$2$, although the thresholds inside the windows change with each of these
parameters. All variants keep the monetary curvature $\gamma_1=2$, which the
primal solver assumes, so the change in $\gamma$ and $\alpha$ moves the
service coefficients $A_i$ only, and the independence is an observation for
these combinations. The arrival rate $\lambda$ enters $\mathcal D$ nowhere, and
raising it to $0.3$ leaves the dates unchanged as well, whereas a larger
increase $\Delta=0.02$ moves the opening date to $0.80$ remaining years
together with the contact set of the high-cost schedule. The selling
region before the event grows more slowly from the floor than the regions
of the two benchmark problems. One tenth of a year after a window opens it
extends about $0.001$ property values above the floor, against $0.01$ to
$0.04$ for the post-event region.

These computations suggest that the preference-free description of the
dates survives uninsurable cost risk in a modified form. The worst-case
schedule supplies the floor, the current cost supplies the incentive, and
preferences and the arrival rate determine how far the selling region
extends above the floor but not when it exists. The first part can be
proved. Admissibility in every state puts the pre-event floors at those
of the high-cost schedule, the pre-event value is nonpositive above the
floor and dominates the value of replicating the high-cost schedule, so
the pre-event dual value of tenure $i$ lies between
$-m_i^{\rm high}z-2(\widetilde K_iz)^{1/2}$ and $-m_i^{\rm high}z$ for a
Merton constant $\widetilde K_i$, and the limit of the dual value
difference at the floor is $m_1^{\rm high}-m_0$. For the second part,
since the jump is not traded, the dual equation is semilinear, with an
infimum over the price of jump risk in which the dual state jumps at the
event to the post-event state with the same wealth. On the contact set
the jump terms of the two tenures differ by a nonpositive amount that
vanishes where the post-event owner also sells, so finite contact requires
the current-cost condition $G_1>0$ of Lemma~\ref{lem:contact-location}.
What is missing for a theorem is the spatial monotonicity of the value
difference and the barrier near the floor for the semilinear equation,
which we leave open. The statement is therefore a numerical observation,
and the portfolio restrictions of Section~\ref{ec:pc} show that it fails
as soon as the owner cannot remove its house-price exposure.

\end{document}